%% file: alex.tex
\documentclass[10pt]{article}
\renewcommand{\include}{\input}

\usepackage{amsmath,amsthm,verbatim,amssymb,amsfonts,amscd,diagrams,graphics}
\theoremstyle{plain}
\newtheorem{theorem}{Theorem}[section]
\newtheorem{corollary}[theorem]{Corollary}
\newtheorem{lemma}[theorem]{Lemma}
\newtheorem{proposition}[theorem]{Proposition}

\theoremstyle{definition}

\theoremstyle{remark}
\newtheorem{remark}[theorem]{Remark}

\newarrow{ul}---->
\newarrow{Backwards}<----   

\newcommand{\hra}{\hookrightarrow} 
\newcommand{\Td}[1]{\Tilde{#1}}
\newcommand{\td}[1]{\tilde{#1}}
\newcommand{\into}{\hookrightarrow}
\newcommand{\Z}{\mathbb{Z}}
\newcommand{\Q}{\mathbb{Q}}
\newcommand{\R}{\mathbb{R}}
\newcommand{\bd}{\partial}

\newcommand{\mv}{\mathversion{bold}}

\newcommand{\vg}{\varGamma}
\newcommand{\blm}[2]{\langle #1 , #2 \rangle}
\newcommand{\bl}[2]{( #1 , #2 )}
\newcommand{\vs}{\varSigma}

\newcommand{\mf}{\mathfrak}

\begin{document}

\title{Alexander polynomials of non-locally-flat knots}\label{S: nobundle}
\author{Greg Friedman\\Yale University\\Dept. of Mathematics\\10 Hillhouse Ave\\PO Box 208283\\New Haven, CT 06520\\friedman@math.yale.edu\\Tel. 203-432-6473  Fax:  203-432-7316}
\date{October 1, 2002}
\maketitle

\begin{abstract}
We generalize the classical study of Alexander polynomials of smooth or PL locally-flat knots to PL knots that are not necessarily locally-flat. We introduce three families of generalized Alexander polynomials and study their properties. For knots with point singularities, we obtain a classification of these polynomials that is complete except for one special low-dimensional case. This classification extends existing classifications for PL locally-flat knots.  For knots with higher-dimensional singularities, we further extend the necessary conditions on the invariants. We also construct several varieties of singular knots to demonstrate realizability of certain families of polynomials as generalized Alexander polynomials. These constructions, of independent interest, generalize known knot constructions such as frame spinning and twist spinning. 

\end{abstract}

\tableofcontents
 
\include{intro.fix}
\include{polyalg}

\include{alexa.fix}

\include{alexb}

\include{alexc}

\include{alexd.fix}
\include{sstrata.fix}
\include{sstratab}

\include{sstratac}

\include{sstratad}

\bibliographystyle{amsplain}
\bibliography{bib}

Several diagrams in this paper were typeset using the
\TeX\, commutative
diagrams package by Paul Taylor.   

\end{document}

%% file: intro.fix.tex

\section{Introduction}

\paragraph{Background.} One of the central motivations for studying knots and
their invariants, including Alexander polynomials, is the central role that
knots play in the understanding of the geometry of subvarieties of real
codimension two. Thus, for a subpseudomanifold, $X^{n-2}$, of a manifold $W^n$,
e.g. for a complex divisor of a complex manifold, the local geometry of $X$ in
$W$ is classically described in terms of link pairs. Therefore, there is a
large
classical topological and algebraic geometric literature   
(e.g. \cite{miln}, \cite{lib3}, \cite{lib2}, \cite{lib1}) which studies
the
topology of the non-singular knots that arise as the link pairs of isolated
singular points. However, in general, the singularities that arise
naturally in
(both high- and low-dimensional) topological and algebraic geometric
situations (see, e.g., \cite{CS2}, \cite{CS3}, \cite{Homer})
cannot be assumed to be
isolated, and the corresponding link pairs of points of the singularities will
consist of knotted sphere pairs which are themselves singular embeddings.

In the non-singular case, more specifically, the focus of knot theory
historically has been the study of smooth or locally-flat codimension two
knots, that is embeddings of $S^{n-2}$ in $S^n$ which are differentiable
or piecewise-linear such that the neighborhood pair of any image point is
PL-homeomorphic to an unknotted ball pair $D^{n-2}\subset D^n$. Furthermore,
much effort has gone into the study of invariants of knots, algebraic
objects which can be assigned to knots and which are identical for
equivalent knots. Prominent among these invariants are the Alexander
polynomials which are  elements, up to similarity class, of the ring of
integral Laurent polynomials $\Lambda:=\Z[\Z]\cong \Z[t,t^{-1}]$. They can
be
defined in may ways, one of which is as follows: By Alexander duality, the
knot complement
$C:=S^n-S^{n-2}$ is a homology circle and hence possesses an infinite
cyclic cover, $\td C$. The homology of $\td C$ with rational coefficients,
$H_i(\td C;\Q)$, has the structure of a module over
$\Gamma:=\Q[t,t^{-1}]$, where the action of $t$ is given by the covering
translation. These modules can be shown to be $\Gamma$-torsion modules
which, since $\Gamma$ is a principal ideal domain, possess square
presentation matrices. The determinants of these matrices are elements of
$\Gamma$ which can be ``normalized'' by ``clearing denominators'' to
elements of $\Lambda$ whose coefficients are relatively prime,
collectively. These normalized determinants are the Alexander polynomials. An
equivalent approach would be to begin with the homology
modules of $C$ taken with a local coefficient system $\vg$, which is  given
by stalk $\Gamma$ and
action determined by factoring the fundamental group to the group of
covering translations.

In \cite{L66}, Levine completely characterized the Alexander polynomials of
PL-locally-flat knots. (Some of these results were known
somewhat earlier for low dimensions; see \cite{FC}, \cite{Seif},
\cite{Lev}, \cite{Kin}.) If
we represent the polynomial corresponding
to
the homology group in dimension $i$ by $p_i$, he showed that the following
conditions are necessary and sufficient for the collection $\{p_i\}$,
$0<i<n-1$, to be the Alexander polynomials of such a knot $S^{n-2}\subset
S^n$ (in the other
dimensions the polynomials are trivial):
\begin{enumerate}
\item $p_i(t)\sim p_{n-i-1}(t^{-1})$ (``$\sim$'' denotes similarity in
$\Lambda$),
\item $p_i(1)=\pm 1$,
\item if $n=2q+1$, $q$ even, then $p(-1)$ is an odd square.
\end{enumerate}

In this paper, we study the generalization of these invariants
and
their properties to various classes of knots which are not necessarily
locally-flat, that is \emph{non-locally-flat} (we also sometimes refer to knots which are
definitely not locally-flat as \emph{singular}). For
knots with only point singularities, we establish necessary and
sufficient
conditions generalizing those of Levine and which form a complete
characterization in all dimensions save $n=5$ (and even for $n=5$, we come
close to a complete characterization; see below). For knots with more general
singularities, we further generalize the necessary conditions, and we
study several methods for realizing given sets of polynomials,
including a construction of independent interest in the study of smooth
knots which generalizes twist spinning \cite{Z65}, superspinning
\cite{C70}, and frame spinning \cite{Ro89}. 

We now outline our results section by section.

\paragraph{Section \ref{S: poly alg}: Polynomial algebra.} This is a
preliminary section in which
we develop some fundamental results of what we call \emph{polynomial
algebra}
by analogy with homological algebra. In particular, to each torsion
$\Gamma$-module  there is
associated an element of $\Gamma$, up to similarity class. 
This is the determinant of the presentation matrix of the module or,
equivalently, the product of its torsion coefficients (recall that
$\Gamma=\Q[\Z]=\Q[t,t^{-1}]$ is a principal ideal domain). We develop some 
relationships among the polynomials associated to the modules of an exact
sequence.

\paragraph{Section \ref{S: disk knots}: Sphere knots with point
singularities
and locally-flat disk knots.} In this section, we first show that the
study
of the complement of a
knot with a
point singularity is homologically equivalent to the study of the complement
of a locally-flat
proper disk knot whose boundary sphere knot is the link pair of the singular
point. In fact, by a technique of Milnor and Fox \cite{M66}, the same is true of a
sphere knot with any finite number of point singularities, and the boundary sphere
knot will be the knot sum of the link pair knots of all the singular points. In
this context, we define a
family of three polynomials: $\lambda_i$, corresponding to the homology module of
the cover of the disk knot complement, $C$; $\nu_i$, corresponding to the boundary
sphere knot complement, $X$; and $\mu_i$, corresponding to the relative
homology of the
cover of the pair $(C,X)$. Furthermore, there is a natural factorization of these
polynomials: $\nu_i\sim a_ib_i$, $\lambda_i\sim b_i c_i$, and $\mu_i\sim
c_ia_{i-1}$. With this notation we prove the following theorem:

\begin{theorem}[Theorem \ref{T: alexa}]\label{T: point sing**}
For $n\neq 5$ and $0<i<n-1$, $0<j<n-2$, the following conditions are necessary and
sufficient
for $\lambda_i$, $\mu_i$, and $\nu_j$ to be the polynomials associated to
the
$\Gamma$-modules $H_i(\td C;\Q)$, $H_i(\td C,\td X;\Q)$, and $H_j(\td X;\Q)$ of a
locally-flat proper disk knot $D^{n-2}\subset D^n$:  There exist polynomials
$a_i(t)$, $b_i(t)$, and $c_i(t)$, primitive in $\Lambda$, such that
\begin{enumerate}
\item\begin{enumerate}
\item $\nu_i\sim a_i b_i$
\item $\lambda_i\sim b_ic_i$
\item $\mu_i\sim c_i a_{i-1}$
\end{enumerate}
\item \label{I: dual*}\begin{enumerate}
\item $c_i(t)\sim c_{n-i-1}(t^{-1})$
\item \label{I: duala*} $a_i(t)\sim b_{n-i-2}(t^{-1})$
\end{enumerate}
\item \label{I: disk*} $a_i(1)=\pm1, b_i(1)=\pm1 , c_i(1)=\pm 1$,
$a_0(t)=1$. 
\item \label{I: mddisk*} If $n=2q+1$ and $q$ is even, then there exist an integer
$\rho$ and an integer $\omega\geq 0$ such that $\frac{(1-t)^{\omega}\rho}{\pm
c_q(t)}$ is
the discriminant of a skew Hermitian form $A\times A\to Q(\Lambda)/\Lambda$ on
a finitely-generated $\Lambda$-module, $A$, on which multiplication by $t-1$ is
an isomorphism (or equivalently, $c_q(t)=$det$[M(t)]$, where $M(t)=
(-1)^{q+1}(R^{-1})'\tau R t-\tau'$ for integer matrices $\tau$ and $R$ such
that $R$ has non-zero determinant 
and $(R^{-1})'\tau R$ is an integer matrix (here $'$ indicates transpose); see Section
\ref{S: middim disc}
for more details). \end{enumerate}

For a locally-flat proper disk knot $D^{3}\subset D^5$ (the case $n=5$),
these
conditions are all necessary. Furthermore, we can construct knots which
satisfies both these conditions and the added, perhaps unnecessary,
condition that $|c_2(-1)|$ be an odd square.
\end{theorem}  

In Section \ref{S: nec cond disk},
we prove the necessity of the duality and
normalization conditions
\eqref{I: dual*} and \eqref{I:
disk*} by a generalization of Levine's technique in \cite{L66} by (i) constructing
an appropriately generalized Seifert surface, (ii) using the surface to
construct the cover by a
cut-and-paste procedure, (iii) deducing from the homology modules of the pieces
of the construction the form of the presentation matrices of the desired
modules, and (iv) exploiting the properties of an integer linking pairing
between the homology of the Seifert surface and that of its complement to show
that these matrices have
the requisite properties to induce those claimed for the
polynomials. 

In Section \ref{S: Realization},  we prove the
sufficiency of these conditions, modulo condition  \eqref{I: mddisk*}, by 
employing various explicit constructions using surgery and relative
surgery. In particular, we show complete sufficiency under the added
(unnecessary) condition that if $n=2q+1$ and $q$ is even, then $|c_q(-1)|$
is an odd square.

Section \ref{S: middim disc} contains a study of the additional issues which
are involved in characterizing the ``middle dimension polynomial'', $c_q(t)$,
for $n=2q+1$ and $q$ even. The necessity of condition
\eqref{I: mddisk*} is a consequence of the existence of a skew Hermitian
form on
the module \emph{ker}$(\bd_*:H_q(\td C,\td X;\Q)\to H_{q-1}(\td X;\Q))$ which
we deduce from the Blanchfield pairing \cite{B57}. 
The realization of a  given polynomial $c_q(t)$, for $n\neq 5$, is deduced as a
consequence of the following more general theorem:

\begin{proposition}[Proposition \ref{P: pairing}]
Let $A$ be a finitely generated $\Z$-torsion free $\Lambda$-module on which
multiplication by $t-1$ is an automorphism and on which there is a 
non-degenerate
conjugate linear $(-1)^{q+1}$-Hermitian pairing $\blm{\,}{}:A\times A \to
Q(\Lambda)/\Lambda$. Then there exists a disk knot $D^{n-2}\subset D^{n}$,
$n=2q+1$, $q>2$, such that:
\begin{enumerate}
\item $H_q(\td C)\cong A$,
\item $H_i(\td C)=0$, $0<i<n-1$, $i\neq q$,
\item $H_i(\td X)=0$, $0<i<n-2$, $i\neq q-1$,
\item $H_{q-1}(\td X)=0$ is a $\Z$-torsion module,
\item $H_i(\td C, \td X)=0$, $0<i<n-1$, $i\neq q$,

\item the form on $H_q(\td C)$ is given by $\blm{\,}{}$. (Note that
$H_q(\td
X)=0$ implies that $H_q(\td C)\cong A \cong$\emph{ker}$(\bd_*)$ in the long exact
sequence).
\end{enumerate}
\end{proposition}

The
impediment to a complete characterization in dimension $n=5$ is a consequence
of
special difficulties associated with low-dimensional surgery and is related to
an open problem of Levine's in the study of pairings on low-dimensional
locally-flat sphere knots \cite{L77}.

\paragraph{Section \ref{S: sing}: Knots with more general singularities.}
In
this section, we consider the case of a sphere knot $K=S^{n-2}\subset
S^n$ with singular set $\Sigma$, which need no longer consist solely of
isolated points. It remains useful to study not the actual knot complement,
$S^n-K$,
but
the homotopy equivalent complement of the locally-flat restriction of the knot
to $S^n-N(\Sigma)$, where $N(\Sigma)$ is an open regular neighborhood of
$\Sigma$. Then, we again obtain a boundary ``knot'' which is the complement in
$\bd \bar N(\Sigma)$ of its intersection with the knot $K$. Accordingly,
we can again define three sets of polynomials (corresponding to the boundary,
absolute, and relative homology modules of the covers) which again have natural
factorizations $\nu_i\sim a_ib_i$, $\lambda_i\sim b_i c_i$, $\mu_i\sim c_i
a_{i-1}$. In this setting, by further generalizing the above techniques and by
employing a number of homological algebra computations, we show in
Section \ref{S: sstrata} that the
necessary conditions of  Theorem \ref{T: point sing**} generalize as follows:

\begin{theorem}[Theorem \ref{T: sstrata}]
Let $\nu_j(t)$, $\lambda_i(t)$, and $\mu_i(t)$, $0<j<n-2$ and $0<i<n-1$,
denote the Alexander polynomials corresponding to
$H_j(\td X)$, $H_i(\td C)$, and $H_i(\td C, \td X)$, respectively, of a
knotted $S^{n-2}\subset S^n$. We can assume these
polynomials to be primitive in $\Lambda$. Then, there exist polynomials
$a_i(t)$, $b_i(t)$, and $c_i(t)$, primitive in
$\Lambda$, such that
                \begin{enumerate}

                \item $\nu_j(t)\sim a_j(t)b_j(t)$,

                \item $\lambda_i(t)\sim b_i(t)c_i(t)$,

                \item $\mu_i(t)\sim c_i(t)a_{i-1}(t)$,

\item   $a_i(t)\sim b_{n-2-i}(t^{-1})(t-1)^{\td B_i}$,
\item $c_i(t)\sim c_{n-1-i}(t^{-1})$,

\item $b_i(1)=\pm 1$,
\item $c_i(1)=\pm 1$,

\item if $n=2q+1$, then $c_q(t)$ is the determinant of a matrix of the
form $(R^{-1})'\tau R t-(-1)^{q+1}\tau'$ where $\tau$ and $R$ are
matrices such that $R$ has non-zero determinant.
 \end{enumerate}

Furthermore, if $\mf H_q =
H_q(\td{C};\Q)/$\emph{ker}$(H_q(\td{C};\Q) \to H_q(\td C,\td X;\Q ))$ and
$n=2q+1$, there is a
$(-1)^{q+1}$-Hermitian pairing
$\blm{\,}{}: \mathfrak H_q \times \mathfrak H_q \to Q(\Gamma)/\Gamma$
which has
a matrix representative $\frac{t-1}{(R^{-1})'\tau -(-1)^{q+1}
\tau'tR^{-1}}$
with respect to an appropriate basis.
\end{theorem}

In this setting of general singularities, realization of polynomials is more
difficult because the allowable set of polynomials will depend subtly on the
properties of the singular set, its link pairs, and its embedding. However,
in Section \ref{S: construction*}, we
employ several constructions available for creating locally-flat knots
including the \emph{frame spinning} of Roseman \cite{Ro89} and our own
generalization to \emph{frame twist-spinning}. Together, these include as special
cases the superspinning of Cappell \cite{C70} and the twist spinning of Zeeman
\cite{Z65}. By adapting these techniques and generalizing them to knots with
singularities, it is possible to construct singular knots and to obtain some
realization results here as well. In particular, for any manifold $M$ which
can be embedded with framing in $S^{n-2}$, we construct classes knots
$S^{n-2}\subset
S^n$ whose singular sets are $M$. 

Furthermore, we calculate the polynomials of
the knots so constructed based upon the polynomials of the knots being spun and
the homology properties of the manifolds, $M$, they are being spun about. Let
$\lambda^{\sigma}_i$, $\mu^{\sigma}_i$, and $\nu^{\sigma}_i$ denote the
polynomials of a frame spun knot; $\lambda^{\tau}_i$, $\mu^{\tau}_i$, and
$\nu^{\tau}_i$ the polynomials of a frame twist-spun knot; and
$\lambda_i$, $\mu_i$, and $\nu_i$ the polynomials of the knot $K$ being spun.
Let $\Sigma$ be the singular set of the knot $K$. Denpte the  Betti numbers of $\Sigma$ by 
$\mf b_i$, let $B_i$ be the
$i$th Betti number of $M^k$, and let $\td \beta_i$ be the reduced Betti
number of $M\times \Sigma$. Suppose that  $H_j(M^k;
\vg|_{M^k})\cong \Gamma^{\mathfrak B_j}\oplus\oplus_l\Gamma/(\zeta_{jl})$
(see Section \ref{S: twist spin} for the definition of this local coefficient
system on $M$) and that the torsion coefficients of the boundary knot of $K$
which are relatively prime to $t-1$ are denoted by $\nu_{il}$, so that
$\nu_i=(t-1)^{\mf b_i}\prod_l\nu_{il}$. Similarly, let
$\lambda_i=\prod_l\lambda_{il}$ and $\mu_i=(t-1)^{\td{\mf
b}_i}\prod_l\mu_{il}$. 
Then, we show in Sections \ref{S: frame spin} and \ref{S: twist spin}
that:
\begin{align*}
\lambda_i^{\sigma}(t)=&\prod_{l=1}^{m-2} [\lambda_l(t)]^{B_{i-l}}\\
\mu_i^{\sigma}(t)= &(t-1)^{\td B_{i-1}}\prod_{l=0}^{m-2}
[\mu_l(t)]^{B_{i-l}}\\
\nu_i^{\sigma}(t)=&\prod_{l=0}^{m-3}[\nu_l(t)]^{B_{i-l}}\\
\lambda_j^{\tau}(t)=&
\prod_{\overset{r+s=j}{s>0}}\left((\prod_l\lambda_{sl}^{\mathfrak
B_r}\cdot\prod_{i,l}d(\zeta_{ri},\lambda_{sl})\right) \cdot
\prod_{\overset{r+s=j-1}{s>0}}\left(\prod_{i,l}d(\zeta_{ri},\lambda_{sl})\right)\\
\mu_j^{\tau}(t)=&
(t-1)^{\td \beta_{j-1}}
\prod_{\overset{r+s=n-j-1}{s>0}}\left((\prod_l\mu_{m-s-1,l}^{\mathfrak
B_r}\cdot
\prod_{i,l}d(\bar\zeta_{ri},\mu_{m-s-1,l})\right) \\
&\qquad\qquad\cdot
\prod_{\overset{r+s=n-j-2}{s>0}}\left(\prod_{i,l}d(\bar\zeta_{ri},\mu_{m-s-1,l})\right)\\
\nu^{\tau}_j(t)=&(t-1)^{
\beta_j}\prod_{r+s=j}\left(\prod_l\nu_{sl}^{\mathfrak B_r} \cdot
\prod_{i,l}d(\zeta_{ri},\nu_{sl})\right)
\cdot  \prod_{r+s=j-1}\left(\prod_{i,l}d(\zeta_{ri},\nu_{sl})\right).
\end{align*}

From these formulas, we then deduce the following realization theorems: 

\begin{proposition}[Proposition \ref{P: frame spin real}]
Let $M^k$ be a manifold which embeds in $S^{n-2}$ with trivial normal bundle
with framing $\phi$ and
such that $n-k>3$. Let $\Sigma$ be a single point. Let $B_i$ denote the $i$th
Betti number of $M$,
and let $\td{\mathfrak b}_i$ and
$\td \beta_i$ denote the $i$th reduced Betti numbers of $\Sigma$ and $M\times
\Sigma$, respectively.
Suppose that we are given any set of polynomials, $a_i(t)$, $b_i(t)$,
$c_j(t)$ and
$c'_l(t)$,
$0<i<n-k-2$, $0<j<n-k-1$, and $0<l<n-1$, which satisfy:
\begin{enumerate}
\item   $a_i(t)\sim b_{n-k-2-i}(t^{-1})$,
\item $c_i(t)\sim c_{n-k-1-i}(t^{-1})$,
\item $c'_i(t)\sim c'_{n-1-i}(t^{-1})$,
\item $b_i(1)=\pm 1$,
\item $c_i(1)=\pm 1$,
\item $c'_i(1)=\pm 1$,
\item if $n-k=2p+1$, $p$ even, $p\neq 2$, then $c_p(t)$ is the
determinant of a matrix of
the form $(R^{-1})'\tau R 
t-(-1)^{q+1}\tau'$ where $\tau$ and $R$ are integer matrices such that $R$ has
non-zero determinant and
$(R^{-1})'\tau R$ is an integer matrix; if $n-k=2p+1$, $p$ even, $p= 2$,
then
$|c_p(-1)|$ is an odd square,
\item if $n=2q+1$, $q$ even, then $|c'_q(-1)|$ is an odd square.
\end{enumerate}
Then there exists a knotted $S^{n-2}\subset S^n$ with singular set $M$ and
Alexander subpolynomials $a_i^{\sigma}(t)$,  $b_i^{\sigma}(t)$,
and $c_i^{\sigma}(t)$ satisfying
\begin{align*}
a_i^{\sigma}(t)&\sim (t-1)^{\td \beta_{i}} \prod_{l=1}^{m-2}
\left[a_l(t)\right]^{B_{i-l}}\\
b_i^{\sigma}(t)&\sim\prod_{l=1}^{m-2} [b_l(t)]^{B_{i-l}}\\
c_i^{\sigma}(t)&\sim c'_i(t) \prod_{l=1}^{m-2} [c_l(t)]^{B_{i-l}}. \\
\end{align*}
\end{proposition}

\begin{theorem}[Theorem \ref{T: twist spin real}]
Let $M^k$, $n-k>3$, be a manifold which embeds in $S^{n-2}$ with trivial normal bundle
with framing $\phi$. Given a map $\tau: M\to S^1$, let $\mathfrak B_i$ be the rank of the
free part and $\zeta_{il}$ be the torsion invariants of the $\Gamma$-modules
$H_i(M;\vg|_M)$. (These modules
are independent of the knot being spun in the construction.) If $\gamma\in \Gamma$, then
let
$\bar \gamma\in \Gamma$ be such that $\bar \gamma(t)=\gamma(t^{-1})$. If $K$ is a knot
$S^{m-2}\subset S^m$ with Alexander invariants $\lambda_{il}$, $\mu_{il}$, and $\nu_{il}$
and with singular set $\Sigma$ with reduced Betti numbers $\td{\mathfrak b}_i$, then
there exists a frame twist-spun knot $\sigma^{\phi,\tau}_M(K)$ with singular set $M\times
\Sigma$ (whose reduced Betti numbers we denote $\td \beta_i$) and with Alexander
polynomials given for $j>0$ by:
\begin{align*}
\lambda_j^{\tau}(K)\sim& 
\prod_{\overset{r+s=j}{s>0}}\left((\prod_l\lambda_{sl}^{\mathfrak
B_r}\cdot\prod_{i,l}d(\zeta_{ri},\lambda_{sl})\right) \cdot
\prod_{\overset{r+s=j-1}{s>0}}\left(\prod_{i,l}d(\zeta_{ri},\lambda_{sl})\right)\\
\mu_j^{\tau}(K)\sim&   
(t-1)^{\td \beta_{j-1}}
\prod_{\overset{r+s=n-j-1}{s>0}}\left((\prod_l\mu_{m-s-1,l}^{\mathfrak B_r}\cdot
\prod_{i,l}d(\bar\zeta_{ri},\mu_{m-s-1,l})\right) \\
&\qquad\qquad\cdot
\prod_{\overset{r+s=n-j-2}{s>0}}\left(\prod_{i,l}d(\bar\zeta_{ri},\mu_{m-s-1,l})\right)\\
\nu^{\tau}_j(t)\sim &
(t-1)^{\td \beta_j}\prod_{r+s=j}\left(\prod_l\nu_{sl}^{\mathfrak B_r} \cdot
\prod_{i,l}d(\zeta_{ri},\nu_{sl})\right)
\cdot  \prod_{r+s=j-1}\left(\prod_{i,l}d(\zeta_{ri},\nu_{sl})\right).
\end{align*}
In particular, by frame twist-spinning knots with a single point as their singular set,
we obtain knots with $M$ as their singular sets.
\end{theorem}

\begin{remark}
In fact, we can create a knot with a single point as its singular set
and with (nearly) any given set of allowable invariants by the results of
Section
\ref{S: disk knots}. Putting this together with the above
theorem, we know exactly what kinds of polynomials can be realized as
those of frame twist-spun knots with singular set $M$, modulo our ability
to compute the homology $H_j(M;\vg|_M)$ and our  difficulty with
the polynomial $c_2(t)$ of a disk knot $D^3\subset D^5$. 
\end{remark}

Finally, we
form singular knots by the suspension of locally-flat or singular knots
and compute their polynomials ($\lambda_i^{\vs}$, $\mu_i^{\vs}$, and
$\nu_i^{\vs}$) from those of the original the
knots  ($\lambda_i$, $\mu_i$, and
$\nu_i$). This is done in Section \ref{S: suspend}, where we obtain the following result:

\begin{proposition}[Proposition \ref{P: suspend real}]
With the notation as above,
\begin{enumerate}
\item $\lambda_i^{\vs}\sim \lambda_i\sim b_ic_i$
\item $\mu_i^{\vs}\sim \mu_{i-1}\sim c_{i-1}a_{i-2}$
\item $\nu_i^{\vs}\sim \lambda_i\mu_i \sim a_{i-1}b_{i}c_i^2$.
\end{enumerate}
\end{proposition}

This work originally appeared as part of the author's dissertation \cite{GBF}. In further papers, we study the intersection homology analogues of Alexander polynomials for non-locally-flat knots  (see \cite{GBF}, \cite{GBF2}). I thank my advisor, Sylvain Cappell, for all of his generous and invaluable guidance. 

%% file: polyalg.tex

\section{Polynomial algebra}\label{S: poly alg}

Let $\Gamma=\Q[\Z]=\Q[t,t^{-1}]$ be the ring of Laurent polynomials with 
rational coefficients. In other words, the elements of $\Gamma$ are
polynomials $\sum_{i\in \Z} a_it^i$, such that each $a_i\in \Q$ and
$a_i=0$ for
all but a finite number of $i$. $\Gamma$ is a principal ideal domain
\cite[\S 1.6]{L66}. Unless otherwise specified, we will generally
not distinguish between elements of $\Gamma$ and their similarity classes
up to unit. In
this introductory section, we study some
basic facts, which will be used often, concerning torsion
$\Gamma$-modules and their associated
polynomials (the determinants of their square presentation matrices). In
analogy with homological algebra for modules, we refer to
this theory of the behavior of  the associated polynomials   as
\emph{polynomial
algebra}.

Let $\Lambda=\Z[\Z]=\Z[t,t^{-1}]$, the ring of Laurent polynomials
with integer coefficients. Then
$\Gamma=\Lambda\otimes_{\Z}\Q$. We call a polynomial in $\Lambda$
\emph{primitive} if its set of non-zero coefficients have no common
divisor
except for $\pm 1$. 
Any element of $\Gamma$ has an associate
in $\Gamma$ which is a primitive polynomial in $\Lambda$: Any
element $at^i\in \Gamma$ is a unit and, in particular then, any $a\in \Q$.
So given an element of $\Gamma$, we can first
clear denominators and then divide out any common divisors without
affecting similarity (associate) class in $\Gamma$. We will often choose
to represent an element of $\Gamma$ (technically, its associate class) by
such a primitive element of $\Lambda$.

\begin{proposition}\label{P:alt. poly.}
Suppose we have an exact sequence of finitely generated torsion
$\Gamma$-modules
\begin{equation}\label{E: exact1}
\begin{CD}
0 @>d_0>> M_1 @>d_1>> M_2 @>d_2>> \cdots @>d_{n-1}>> M_n @>d_n>> 0,
\end{CD}
\end{equation}
and suppose that $\Delta_i$ is the determinant of a square presentation
matrix of
$M_i$ (which we will refer to as the polynomial associated to the module). Then,  taking
$\Delta_{n+1}=1$ if $n$ is odd, the
alternating product $\prod_{i=1}^{\lceil n/2 \rceil}
\frac{\Delta_{2i-1}}{\Delta_{2i}}\in \Q(t)$ is equal to a unit of $\Gamma$,
and, in
particular, with a consistent choice of normalization within associate
classes for the elementary
divisors of the $M_i$ (in the language of \cite{HU}), this product is equal
to
$1$. 
\end{proposition}
\begin{proof}
It is well known (see, for example, \cite[p. 225]{HU}) that a finitely
generate torsion module over a principal ideal domain can be decomposed as the direct sum of cyclic torsion
summands of orders $p_j^{k_j}$, the $p_j$ not necessarily distinct primes in
the ground ring and the $k_j$ positive integers, also not necessarily
distinct. Furthermore, we know that
this decomposition is unique in the sense that the $p_j$ are determined up
to associate class, but the cyclic summands $\Gamma/(p_j^{k_j})$, being
independent of the choice of $p_j$ within the associate class, are uniquely
determined. Hence, in particular, each $M_i$ has a square presentation
matrix of the form
\begin{equation*}
\left( \begin{array}{cccc}
p_{i_1}^{k_{i_1}} & 0 & \cdots & 0\\
0 & p_{i_2}^{k_{i_2}} &      &0\\
\vdots & & \ddots &\vdots\\
0 & 0& \hdots &p_{i_{m_i}}^{k_{i_{m_i}}} 
\end{array}
\right),
\end{equation*}
and $\Delta_i=\prod_{j=1}^{m_i} p_{i_j}^{k_{i_j}}$.
Since we have a finite number of modules, each finitely generated, we have a
finite number of primes of $\Gamma$ occurring in the elementary divisors and 
in these matrices. We are free to choose these primes so that if two are in
the same associate class, then they are, in fact, the same element of
$\Gamma$, and we assign an order so that we may speak of the collection of
distinct primes $\{p_j\}_{j=1}^{m}$ which occur. Let $M_i(p_j)$ be the
summand
of $M_i$ which is the direct sum of cyclic modules of order a power of
$p_j$. This may be a trivial summand.  Then each $M_i$ decomposes as 

\[ M_i\cong M_i(p_1)\oplus M_i(p_2) \oplus \cdots \oplus M_i(p_m),
\]
and, if we set $\Delta_i(p_j)$ to be the determinant of the presentation matrix of $M_i(p_j)$,
then $\Delta_i=\prod_j \Delta_i(p_j)$ and $\Delta_i(p_j)=p_j^k$ where $k$ is the sum of
the powers of $p_j$ which occur in the elementary divisors of $M_i$.

\begin{lemma}\label{L: prime splitting}
Let $r$ and $s$ be powers of distinct (non-associate) prime elements of
$\Gamma$. Then the
only $\Gamma$-module morphism 
$f:\Gamma/(r)\to \Gamma/(s)$ is the $0$ map.
\end{lemma}
\begin{proof}
Suppose $f$ is such a map. Then, letting elements of $\Gamma$ stand for
their classes in $\Gamma/(s)$ and $\Gamma/(r)$ where appropriate, we have 
\[0=f(0)=f(r1)=rf(1)=ra\]
for some $a\in\Gamma/(s)$. But $ra=0$ implies that $ra=sb$ in $\Gamma$ for
some $b\in
\Gamma$. Since $s|sb$ but no prime divisor of $s$ divides $ r$, we must
have $s|a$ so that $a=sc$ for 
some $c\in \Gamma$. But then
$f(1)=a=sc=0$, for some $c$. This implies that $f$ is the $0$ map because
$f$ is completely determined by the image of the generator. 
\end{proof}
\begin{corollary}
With the notation above, the only $\Gamma$-module morphisms $f:M_i(p_k)\to
M_j(p_l)$, for $k\neq l$, are the $0$ maps.
\end{corollary}
\begin{proof}
This follows immediately from the lemma since the map on each summand must
be $0$. 
\end{proof}

\begin{corollary}\label{C: sequence splitting}
For any $p_j$,  the sequence 
\begin{equation}\label{E: exact2}
\begin{CD}
0 @>e_0>> M_1(p_j) @>e_1>> M_2(p_j) @>e_2>> \cdots @>e_{n-1}>> M_n(p_j) @>e_n>> 0 
\end{CD}
\end{equation}
is exact, where the maps $e_i$ are the restrictions of the the maps $d_i$ to the direct
summands $M_i(p_j)$.
\end{corollary}
\begin{proof}

First, we note that these maps are well-defined: Any element $a\in M_i$
can be represented as $\sum_l a_l$ where $a_l\in M_i(p_l)$, and so if
$a_j\in M_i(p_j)$, we can identify it via the inclusion of the summand
with $0+\cdots+0+a_j+0+\cdots+0\in M_i$, which we will also call $a_j$.  
Then $d_i(a_j)$ is represented by a sum $\sum_l b_l$, $b_l\in
M_{i+1}(p_l)$. If $r_l$ is the projection of $M_i$ to $M_i(p_l)$, then
$r_l d_i(a_j)=b_l$. But this gives a $\Gamma$-module morphism $M_i(p_j)\to
M_{i+1}(p_l)$ and hence $b_l=0$ if $l\neq j$ by the preceding corollary.
In other words, the image of the summand $M_i(p_j)$ under $d_i$ lies in
$M_{i+1}(p_j)$. Thus the maps of this  sequence are well-defined.

Next since $d_{i+1}d_i=0$, we also have $e_{i+1}e_i=0$.

It remains to show that Ker($e_{i}$)$\subset$Im($e_{i-1}$). Suppose that $e_{i}(a_j)=0$, $a_j\in M_i(p_j)$. Then we have also $d_i(a_j)=0$, since we have already observed that $d_i$ takes $a_j$ to $0$ in all of the other summands $M_{i-1}(p_l)$, $l\neq j$. But since \eqref{E: exact1} is exact, there is an element $c$ in $M_{i-1}$ such that $d_{i-1}(c)=a_j$, and we have $c=\sum_l c_l$, $c_l\in M_{i-1}(p_l)$ and $d_{i-1}(c)=\sum_l d_{i-1}(c_l)=a_j$. Since we know $d_{i-1}(c_l)\in M_i(p_l)$, we must have $d_{i-1}(c_j)=a_j$ and $d_{i-1}(c_l)=0$, $l\neq j$. But then $a_j= e_{i-1}(c_j)\in$Im($e_{i-1}$). 
\end{proof}

Note that this lemma together with its corollary 
allows us to write the exact sequence
\eqref{E: exact1} as the direct sum of exact sequences of the form \eqref{E: exact2}.

We will prove Proposition \ref{P:alt. poly.}
in the special
case that the exact sequence in its statement has the form of that in
equation \eqref{E: exact2}. The proposition  will then follow
for the general case by the formula
\[ \prod_{i=1}^{\lceil n/2 \rceil} \frac{\Delta_{2i-1}}{\Delta_{2i}}=
\prod_{i=1}^{\lceil n/2 \rceil} \frac{\prod_j \Delta_{2i-1}(p_j)}{ \prod_j \Delta_{2i}(p_j)}= 
\prod_j \prod_{i=1}^{\lceil n/2 \rceil}
\frac{\Delta_{2i-1}(p_j)}{\Delta_{2i}(p_j)}.\]

So it remains to prove that the exact sequence \eqref{E: exact2} implies
that\\  
$\prod_{i=1}^{\lceil n/2
\rceil}\frac{\Delta_{2i-1}(p_j)}{\Delta_{2i}(p_j)}$ is
a unit of $\Gamma$. In particular, with our choice of consistent $p_j$'s
within
the associated classes, this product will be $1$. For this, recall that we
have
already observed that $\Delta_i(p_j)=p_j^{k_i(p_j)}$, where $k_i(p_j)$ is
the
sum of the powers of $p_j$ which occur in the elementary divisors of
$M_i$.
Therefore $\prod_{i=1}^{\lceil n/2
\rceil}\frac{\Delta_{2i-1}(p_j)}{\Delta_{2i}(p_j)}=p_j^{k(p_j)}$ where
$k(p_j)=\sum_{i=1}^{\lceil n/2\rceil} (-1)^{i+1}k_i(p_j)$. We claim that
$k(p_j)=0$, which will complete the proof.

Of course, each $\Gamma$-module $M_i(p_j)$ has the underlying
structure of a rational vector space if we forget about the $t$ action,
and
similarly the exact sequence \eqref{E: exact2} can be regarded as an exact
sequence of vector spaces over $\Q$. Suppose $p_j=\sum_{l=a}^b c_l t^l$,
where
$a$ and $b$ are finite integers, $c_a\neq 0$ and $c_b\neq 0$ (we can
always
find such a representation of an element of $\Gamma$). Define $\Vert
p_j\Vert=b-a$. Then the dimension of $\Gamma/(p_j)$ as a rational vector
space
is $\Vert p_j \Vert$, and, more generally, the dimension of
$\Gamma/(p_j^k)$ as
a vector space is $\Vert p_j^k\Vert=k\Vert p_j\Vert$, for any non-negative
integer $k$. Therefore, the dimension of $M_i(p_j)$ as a rational vector
space
must be $k_i(p_j)\Vert p_j\Vert$. But since \eqref{E: exact2} is an exact
sequence of vector spaces, 
\begin{align*}
0&=\sum_i
(-1)^{i+1}\text{dim}(M_i(p_j))=\sum_i
(-1)^{i+1} k_i(p_j)\Vert p_j\Vert\\ 
&=\Vert p_j\Vert\sum_i (-1)^{i+1}
k_i(p_j)=\Vert p_j\Vert k(p_j).
\end{align*}
Since $\Vert p_j\Vert\neq 0$ (else $p_j$
would be a unit of $\Gamma$ and $\Gamma/(p_j^k)$ trivial), we must have
$k(p_j)=0$ as claimed. This completes the proof.

Note that had we not fixed the $p_j$ within their associate classes, the
product $\prod_{i=1}^{\lceil n/2 \rceil}
\frac{\Delta_{2i-1}}{\Delta_{2i}}$
would not necessarily be $1$, but it would still follow from minor
adjustments
to the above arguments that it would be a unit of $\Gamma$.

\end{proof}

\begin{corollary}\label{C: kern} 
With the notation and assumptions as above, each
$\Delta_i=\delta_{i}\delta_{i+1}$ where
$\delta_{i+1}|\Delta_{i+1}$ and $\delta_i|\Delta_{i-1}$. Furthermore, if we represent the $\Delta_i$ by the elements in their similarity classes in $\Gamma$  which are primitive in $\Lambda$, the $\delta_i$ will also be primitive in $\Lambda$.  
\end{corollary}
\begin{proof}
Let $Z_i\subset M_i$ denote the kernel of $d_i$. Then we have the short
exact
sequences
\[
\begin{CD}
0@>>> Z_i @>>> M_i @>>> Z_{i+1} @>>> 0
\end{CD}
\]

Let $\delta_i$ be the determinant of a square presentation matrix of
$Z_i$. Then, applying the above proposition for various choices of $i$, we
obtain $\Delta_i=\delta_i\delta_{i+1}$ up to associate classes, as well as
$\Delta_{i-1}=\delta_{i-1}\delta_i$ and
$\Delta_{i+1}=\delta_{i+1}\delta_{i+2}$.  This proves the first part of
the corollary. For the second, recall that we can always find an element
in the associate class of $\delta_i$ in $\Gamma$ which is primitive in
$\Lambda$, and this choice will be unique up to associate class in
$\Lambda$. Similarly for $\delta_{i+1}$. But the product of two primitive
elements of $\Lambda$ is again primitive in $\Lambda$ (the argument of
\cite[\S 3.10]{He} for $\Z[t]$ extends easily), so that, with this choice,
$\delta_i\delta_{i+1}$ is a primitive element of $\Lambda$ which is
equal to $\Delta_i$ up to associativity in $\Lambda$.  
\end{proof}

This corollary will be used often in what follows.

For convenience, we introduce the following notation. Suppose $\Delta_i\in \Gamma$. We will refer to an 
\emph{exact sequence of polynomials}, denoted by
\begin{equation*}
\begin{CD}
@>>> \Delta_{i-1} @>>>\Delta_i @>>> \Delta_{i+1} @>>>,
\end{CD}
\end{equation*}
to mean a sequence of polynomials such that each $\Delta_i\sim \delta_i \delta_{i+1}$, $\delta_i\in
\Gamma$. As we have seen, such a sequence arises  in the case of an exact sequence of
torsion
$\Gamma$-modules, $M_i$, and, in that case, the factorization of the polynomials is determined by the
maps of
the modules as in Corollary \ref{C:
kern}. In particluar, each $\delta_i$ is the polynomial of the module ker$(M_i\to M_{i+1})$. 
 
Observe that knowledge of two thirds of the terms of an exact sequence of polynomials (for
example, all $\Delta_{3i}$ and $\Delta_{3i+1}$, $i\in\Z$) and the common factors of those
terms (the $\delta_{3i+1}$), allows us to deduce the missing third of the sequence
($\Delta_{3i+2}= \delta_{3i+2}\delta_{3i+3}=
\frac{\Delta_{3i+1}}{\delta_{3i+1}}\cdot\frac{\Delta_{3i+3}}{\delta_{3i+4}}$).

Note also that for any bounded exact sequence of polynomials (or even a half-bounded 
sequence),
the collections $\{\Delta_i\}$ and
$\{\delta_i\}$ carry the same information. That is, suppose that one (or both) end(s) of the
polynomial sequence is an infinite number of $1$'s (by analogy to extending any bounded or 
half-bounded exact module
sequence to an infinite number of $0$ modules).  Clearly, the $\Delta_i$ can be reconstructed
from the
$\delta_i$ by $\Delta_i\sim \delta_i \delta_{i+1}$. On the other hand, if $\Delta_0$ is the
first nontrivial term in the polynomial sequence, then $\delta_0\sim 1$, $\delta_1\sim
\Delta_0$, and
$\delta_i\sim \Delta_{i-1}/\delta_{i-1}$ for all $i>1$. Similar considerations hold for a
sequence which is bounded on the other end. Therefore, we will often study properties of
the polynomials $\Delta_i$ in an exact sequence by studying the $\delta_i$ instead. We will
refer to the $\delta_i$ as the \emph{subpolynomials} of the sequence and to the process of
determining the subpolynomials from the polynomials as ``dividing in from the outside of
the sequence''.

%% file: alexa.fix.tex

\section{Sphere knots with point singularities and locally-flat disk
knots}\label{S: disk knots}

\subsection{Introduction}

Our goal in this section is to study the Alexander polynomials of a knot
with isolated
singularities. More specifically, let $\alpha: S^{n-2}\into S^n$, $n>3$ be
a PL-embedding such that for $x\in \alpha (S^{n-2})$, the link pair of $x$
in $(S^n, \alpha(S^{n-2}))$ is PL-homeomorphic to the standard unknotted sphere pair except at
finitely
many $x$, where it may be a knotted sphere pair. Henceforth, we will
dispense with $\alpha$ unless necessary and refer simply to the knot pair
$(S^n, K=\alpha(S^{n-2}))$ or the $n$-knot $K$. Just as in the classical
locally flat case, Alexander duality tells us that the homology of the knot
complement $S^n-K$ is that of a circle,  and this allows us to study the
infinite
cyclic cover of the knot complement and its homology regarded as a module
over $\Gamma=\Q [\pi_1(S^n-K)]=\Q [\Z ]=\Q [t,t^{-1}]$. We can then study
the Alexander invariants of these modules. 

We begin by seeing that the study of the homological
properties of the complements of sphere knots with isolated singularities
reduces to the study of the complements of locally-flat disk knots. 
This study of disk knots starts by emulating J. Levine's study of
Alexander invariants for the locally flat sphere knots \cite{L66}. In
Section \ref{S: nec cond disk}, we
introduce two sets of polynomial invariants, $\lambda_q^i$ and $\mu_q^i$,
corresponding to certain absolute and relative homology modules and show
that they satisfy certain duality and normalization conditions. From
these, we arrive at the corresponding definitions and properties for
the Alexander polynomials $\lambda_i$ and $\mu_i$ (see Section
\ref{S: cors and poly def}). 

In Section \ref{S: Realization}, we turn to the realization of  
locally-flat disk knots with given polynomial invariants which satisfy the
properties obtained in Section  \ref{S: nec cond disk}. We show that any
allowable set of $\lambda_i$ can be realized, first for a knotted $D^2$ in
$D^4$ (Section \ref{S: D2 in D4}) and then for arbitrary $D^{n-2}\subset
D^n$, $n>4$ (Section \ref{S: large n}). In Section \ref{S: real
disk knot}, we show that we can nearly completely characterize all three 
sets of Alexander polynomials which can occur for a locally flat disk knot
(the polynomials $\lambda_i$ and $\mu_i$, which we have already mentioned,
plus the Alexander polynomials of the boundary locally-flat sphere
knot). The barrier to a complete classification, at that point, is a
certain polynomial factor shared by $\lambda_q$ and $\mu_q$ for knotted
$D^{2q-1}\subset D^{2q+1}$, $q$ an even integer.

In Section \ref{S: middim disc}, we take up the study of this
middle-dimensional polynomial factor. We show that it is related to a
certain Hermitian self-pairing induced by the Blanchfield pairing on the
middle-dimension homology modules. We establish the realizability of such
pairings in disk knots and then study the relationship between the
Alexander polynomial factors and the
presentation matrices of the modules and their pairings. This allows us to
state necessary and sufficient conditions for this polynomial factor for
$n\neq 5$. 

Finally, in Section \ref{S: conclusion alexa}, we gather together the
results of Section \ref{S: disk knots}. Theorem \ref{T: alexa} states a
complete set of necessary and sufficient conditions for Alexander
polynomials for locally-flat disk knots $D^{n-2}\subset D^n$, $n\neq
5$. For $n=5$, the classification is nearly complete, but we obtain only
a partial characterization of the middle dimensional polynomial factor.

\subsection{The Knot Complement}
For technical simplicity, we will often study not the knot complement but rather a version of the the homotopy equivalent ``knot exterior''. For locally flat knots this is the exterior of an open tubular (PL-regular) neighborhood of the knot. Similarly, we can consider the exterior of a regular neighborhood of our singular knot. 

\begin{figure}[p]
\scalebox{1.1}{\includegraphics{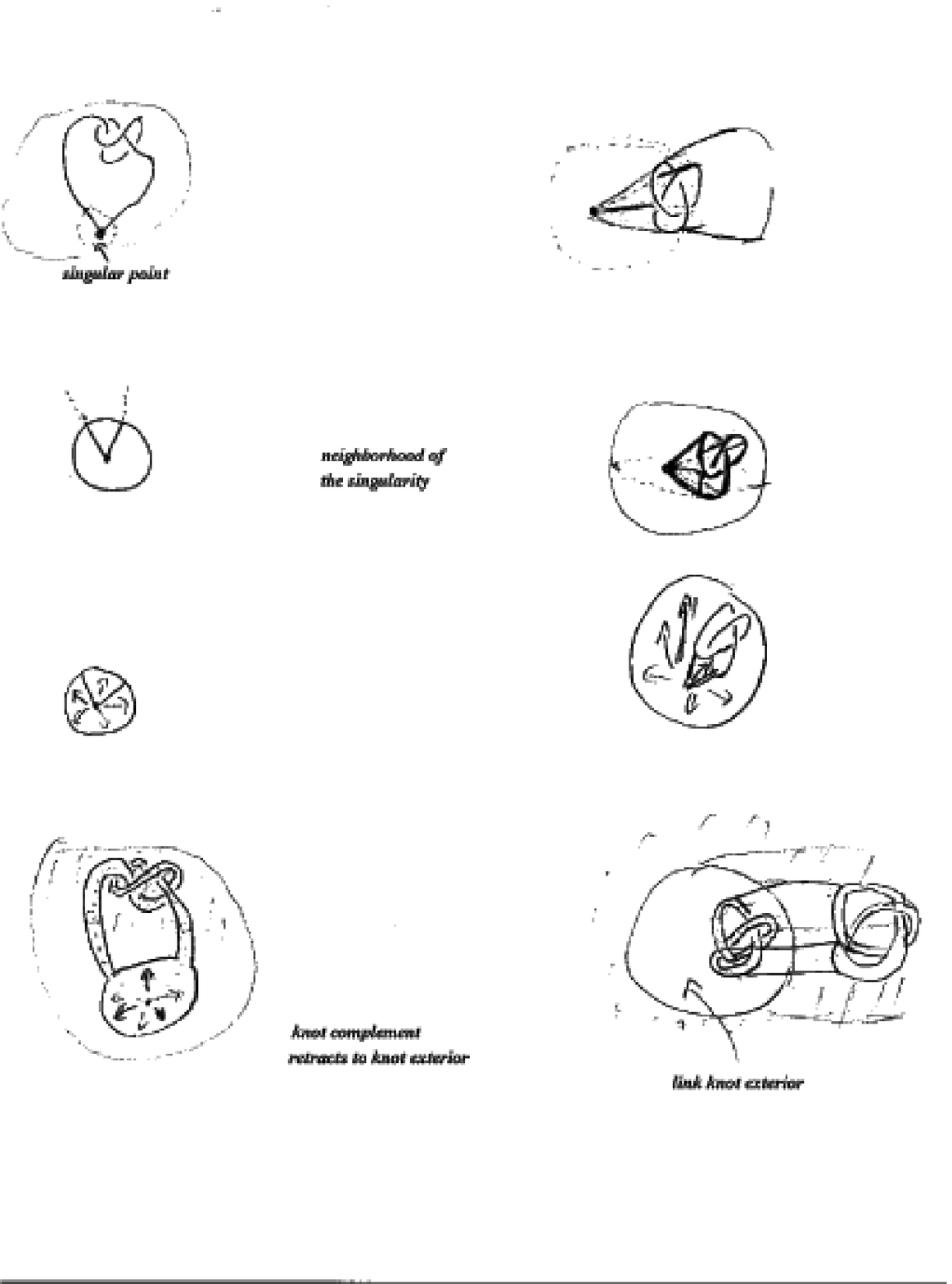}}
\caption{Diagrams of a knot with a point singularity and its
exterior. The pictures on the left are schematic; the corresponding
pictures on the right add a little more depth (the bottom right picture
knots up in dimensions too high to draw).}\label{F: point sing}
\end{figure}

First, assume that the knot, $K$, has only one singular point, $x$. Then the
neighborhood Star($x$) of $x$ in $S^n$ is a knotted ball pair
$(D^n,D^{n-2})$ which is (PL-homeomorphic to) the cone on Link($x$), which
is a knotted sphere pair $\bd (D^n,D^{n-2})=(S^{n-1},S^{n-3})=(S^{n-1},k)$,
where we let $k$ denote the locally flat (n-1)-knot of the sphere pair.
Since the cone point no longer remains when we consider only the knot
complement, we can retract what remains of the complement in Star($x$) out
to the boundary and see that our knot complement is homotopy equivalent to
the complement of a locally-flat knotted disk pair $(D^n, D^{n-2})$ where
this $D^n$ is
the complement of the open disk neighborhood of $x$ in $S^n$. This knotted disk pair in fact provides
a null-knot
cobordism of the slice knot $k$, and the study of the knot complement
reduces (up to homotopy equivalence) to the study of the cobordism
complement $D^n-D^{n-2}$, which we shall denote by $C$. If desired, we can
also retract this
complement to the complement of an open tubular (regular) neighborhood of
the locally-flatly embedded knotted disk in analogy with the usual notion of
knot exteriors. See Figure \ref{F: point sing}.

\begin{figure}[p]
\scalebox{1.1}{\includegraphics{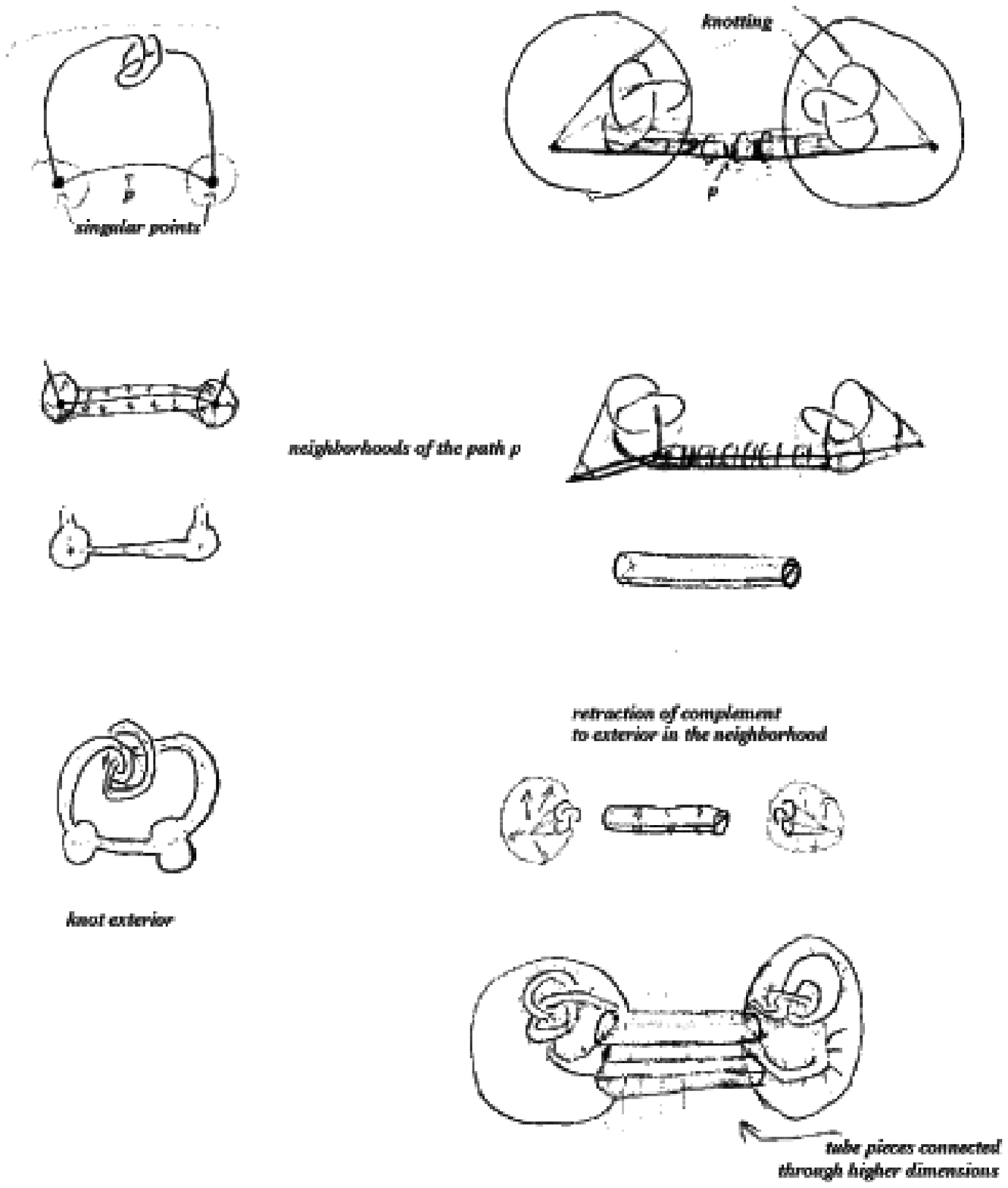}}
\caption{Diagrams of a knot with two point  singularities. Again, the
pictures on the left are schematic and the corresponding
pictures on the right add a little more depth. Note that the neighborhood 
of the path $p$ is a ball and that we retract from the complement of the
knot in this neighborhood to the boundary sphere knot complement of the
disk knot outside of the neighborhood. This boundary knot is the knots
some of the link knots around the points.} \label{F: mult sing}\end{figure}

If $K$ has multiple isolated singularities, $x_i$, the situation is slightly
more complicated but similar. Fox and Milnor's \cite{M66} analysis of the
case of a two-sphere with isolated singularities embedded in four-space
carries over to higher dimensions. In particular, we can choose a PL-arc,
$p$, embedded in $S^{n-2}$ which traverses each singular point (where here
we
confuse $S^{n-2}$ with K). Then a regular neighborhood, $N$, of $p$ is again
a disk pair $(D^n,D^{n-2})$ whose boundary is a knotted sphere pair
$(S^{n-1}, S^{n-3})=(S^{n-1},k)$, where the knot $k$ is the knot sum $\sum_i
k(x_i)$ of the knots of each sphere pair, Link($x_i$). As
in the case of a single singular point, the knot complement is homotopy
equivalent to the complement of the disk pair which is obtained from the
sphere pair by removing the open regular neighborhood of $p$. This can be
seen as follows: First retract
the star neighborhoods of the $x_i$ in the knot complement radially away
from the cone points, $x_i$, as in the last paragraph. The portion of
$S^n-K$ remaining in the interior of $N$ then consists of a disjoint set of
standard ball pairs $(D^n,D^{n-2})$ whose boundaries lie in $\bd N$ except
for two opposing sides (thinking of the balls as cubes) which lie in the
link pairs of $x_i$ and $x_{i+1}$ and can be identified as neighborhoods
there of points of $k(x_i)$ and $k(x_{i+1})$, respectively. But once we have
gone over to the knot complement and hence removed the (n-2)-balls, their
complements easily retract out to $\bd N-(\bd N \cap K)$. Once again, our
study is reduced to the complement of a knotted disk pair which forms the
null-cobordism of a slice knot. Henceforth, we refer to the knotted disk $L$
in $D^n$, $\bd L=k\in \bd D^n$. See Figure \ref{F: mult sing}.

By this discussion, our study of the homological properties of $S^n-K$
reduces to a study of the homological properties of $D^n-L$.

\subsection{Necessary conditions on the Alexander invariants}\label{S: nec
cond}\label{S: nec cond disk}

\subsubsection{Alexander invariants}

We now undertake a study of the Alexander polynomials of the
complements of locally-flat knotted disks following the pattern of Levine's
\cite{L66} study of the Alexander polynomials of locally flat sphere 
knots. In
particular, let $C$ be the disk-knot complement $D^n-L$ and let
$\tilde{C}$ be the infinite cyclic cover associated with the kernel of
the abelianization $\pi_1(C)\to\Z$. Let $t$ denote a generator for the
covering translation and $\Lambda$ the group ring $\Z[\Z]=\Z[t,t^{-1}]$. The
homology groups of $\tilde{C}$ are finitely generated $\Lambda$-modules
since $C$ has a finite polyhedron as a deformation retract, and
the rational homology groups $ H_*(\tilde{C};\Q) \cong
H_*(\tilde{C})\otimes_{\Z} \Q$ are, therefore, finitely generated modules
over the
principal ideal domain $\Gamma=\Q[t,t^{-1}]\cong \Lambda\otimes_{\Z}\Q$.
Therefore, letting $M_R(\lambda)$ denote the R-module of rank 1 with
generator of order $\lambda$, $H_q(\tilde{C};\Q)\cong \bigoplus_{i=1}^k
M_{\Gamma}(\lambda_q^i)$ (note that ``$i$'' is here an index, not a power). 
Furthermore, we can choose the $\lambda_q^i$ so that: 1)
The $\lambda_q^i$ are primitive in $\Lambda$ but are unique up to associate
class in $\Gamma$, and 2) $\lambda^{i+1}_q|\lambda^i_q$. For $0<q<n-1$,
the $\lambda_q^i$ are called the \emph{Alexander invariants} of the knot
complement. We
will also consider the relative homology modules $H_*(\tilde{C},
\tilde{X};\Q)$, where $X$ is the complement of $k$ in $S^{n-1}=\bd D^n$
and $\tilde{X}$ is its infinite cyclic covering. It will be clear from our
construction that $\tilde{X}$ and the cover of $X$ in $\tilde{C}$ are
equivalent. Then $H_*(\tilde{C}, \tilde{X};\Q)$ has the same properties
listed above for $H_*(\tilde{C};\Q)$ and its own Alexander invariants
$\{\mu_q^i\}$, $0<q<n-1$.

We will prove the following theorem concerning necessary conditions on these
polynomials:
\begin{theorem}\label{T: disk knot nec con}
Let $p=n-1-q$. With $\{\lambda_q^i\}$ and $\{\mu_q^i\}$ as above for a knotted disk pair $(D^n, D^{n-2})$, $n\geq3$, the following properties hold:
\begin{enumerate}
\item $\lambda^{i+1}_q|\lambda^i_q$ and $\mu^{i+1}_q|\mu^i_q$ in $\Lambda$;
\item $\lambda^{i}_q(1)=\pm 1$, $\mu^i_q(1)=\pm1$;
\item $\lambda^{i}_q(t)\sim \mu^i_{p}(t^{-1})$ in $\Lambda$, where $\sim$
denotes associativity of elements in $\Lambda$, i.e. $a\sim b$ implies $ a=
\pm t^k
b$ for some $k$.
\end{enumerate}
\end{theorem}  

The proof of the theorem is given over the following sections. 

\subsubsection{Construction of the covering}\label{S: covering 
construction}

We begin by finding $\Gamma$-module presentations for $H_*(\tilde{C};\Q)$ and $H_*(\tilde{C}, \tilde{X};\Q)$ by generalizing the usual technique of studying the Mayer-Vietoris sequence for the infinite cyclic cover obtained from cutting and pasting along a Seifert surface. 

\begin{proposition} \label{P:Seifert} Given a knotted disk $L\in D^n$, there
exists an $(n-1)$-di\-mens\-ion\-al connected bicollared submanifold $V\in
D^n$
such that $\bd V=L\cup F$, where F is a Seifert surface for $k$ in $\bd D^n$. 
\end{proposition}
\begin{proof}
Letting $T$ be a regular tubular neighborhood of $L$ in $D^n$, there is a map
$f: T-L\to S^1$ given by projection on the fibers. (A trivialization of the
disk bundle is provided by the restricting the trivialization of the disk
bundle constituting the tubular neighborhood of the locally-flat sphere knot
obtained by gluing our disk knot and its mirror image along the boundary
knots). As in the construction of the usual Seifert surfaces, this map can
be extended to the rest of $\bd D^n$ so that the inverse image there of a
regular value, $x$, of $S^1$ is a Seifert surface for $k$ after throwing
away extraneous components. In fact, the map can be easily modified so as to
avoid extraneous components (the extraneous components will be bicollared
close manifolds of $\bd D^n$ such that the fiber in the collar of each point
maps to an arc of $S^1$ containing $x$, and so on each fiber we can reverse
the map to run around $S^1$ the other way, avoiding $x$). Now we wish to
extend $f$ to the rest of $D^n-((T-L)\cup \bd D^n)$.

The obstructions to this extension lie in 
\begin{equation*}
H^{i+1}(C,
T\cup \bd D^n;\pi_i(S^1)) \cong H^{i+1}(C, \bd C;\pi_i(S^1))
\end{equation*} 
(see \cite[p. 54]{BR}). We know that
$\pi_i(S^1)=0$ for $i\neq 1$ and $\Z$ for $i=1$ so we need only calculate:
\[
 H^2(C,\bd C;\Z) \cong H_{n-2}(C) \cong H^{1}(L,\bd L=k) \cong
H^{1}(D^{n-2},\bd D^{n-2}) =0, \] 
since $n>3$. The first two isomorphisms are
due to Lefschetz duality and Alexander duality for a ball \cite[p. 426]{MK}.
Therefore, the obstruction is $0$, so the extension exists, and we can take
$V$
as the inverse image of a regular value in $S^1$ after, again, throwing away
extraneous components. If $x$ is no longer a regular value, we can instead
choose a new regular value, $y$, in an $\epsilon$-neighborhood of $x$ such
that $f^{-1}(y)\cap \bd D^n$ is isotopic (in $\bd D^n$) to $f^{-1}(x)$ and
hence gives ``the same'' Seifert surface for $k$. It is clear that $V$ has
the desired boundary. 
\end{proof}

We can now construct $\tilde{C}$ in the usual way by cutting along $V$ to
create a manifold $Y$ whose boundary is $(D^n-V)\cap S^{n-1}$ together with
two copies of $V$, $V_+$ and $V_-$, identified along $L$ and then by pasting
together a countably infinite number of disjoint copies $(Y_i,
V_+^i,V_-^i)$, $-\infty<i<\infty$, of $(Y-L, V_+-L, V_{-}-L)$ by identifying
$V^i_+ -L$ with $V^{i+1}_- -L$ for all $i$. Then $\tilde{X}$ is the
submanifold resulting from the restriction of this construction
to $\bd D^n\cap (Y_i, V_+^i,V_{-}^i)$. $\tilde{X}$ is thus the usual
infinite cyclic cover constructed for a classical knot complement as
claimed. Note that, just as in that case, $H_q(D^n-V)\cong H_q(Y)$ due to
homotopy equivalence. We also
denote $Y\cap \bd D^n$ by $Z$ and have $ H_q(\bd D^n-F)\cong H_q(Z)$.

The usual considerations (see, e.g., \cite{L66}) now allow us to set up the Mayer-Vietoris sequences for $\td{C}$ and $(\td{C},\td{X})$:
\begin{equation}\label{E:MV1}
\begin{diagram}
&\rTo& H_q(V;\Q)\otimes_{\Q}\Gamma&\rTo^{d_1}&
H_q(Y;\Q)\otimes_{\Q}\Gamma&\rTo^{e_1}&H_q(\td{C};\Q)&\rTo
\end{diagram}
\end{equation}
and
\begin{equation} \label{E:MV2}
\begin{diagram}
&\rTo & H_q(V,F;\Q)\otimes_{\Q}\Gamma&\rTo^{d_2}&
H_q(Y,Z;\Q)\otimes_{\Q}\Gamma &\rTo^{e_2}&H_q(\td{C},\td{X};\Q)&\rTo.
\end{diagram}
\end{equation}

We will see that $d_i$, $i=1,2$, is a monomorphism for $0\leq q<n-1$.
Hence $e_i$ is an epimorphism, $0<q<n-1$, and the $d_i$ provide presentation
matrices for the homology modules of the covers as $\Gamma$-modules. That
the $d_i$ are square matrices in this range follows from the following
proposition.

\begin{proposition}\label{P:square matrix}$\,$\newline
\begin{enumerate}
\item $H_q(V;\Q)\cong H_q(Y;\Q)$, $0\leq q<n-1$,
\item $H_q(V,F;\Q)\cong H_q(Y,Z; \Q)$, $0\leq q<n-1$.
\end{enumerate}
\end{proposition}

\begin{proof}
\begin{enumerate}

\item In the proof, we assume rational coefficients while omitting mention
for the sake of notational convenience. $H_q(Y)\cong H_q(D^n-V)\cong H^p(V,
F)$, $p+q=n-1$, by Alexander duality for the ball, while $H_q(V)=H^p(V, \bd
V)$ by Lefschetz duality. So we must show that $H^p(V, F)\cong H^p(V, \bd
V)$. Recall that $\bd V=F\cup_{k}L\cong F\cup_{k} D^{n-2}$. From the reduced
Mayer-Vietoris sequence, we get immediately that $H^k(F)\cong H^k(\bd V$)
for $k<n-3$, and the top of the sequence is

{\footnotesize
\begin{diagram}
\qquad& \lTo^0 & H^{n-2}(F) & \lTo & H^{n-2}(\bd V) & \lTo
&H^{n-3}(S^{n-3}) &
\lTo & H^{n-3}(F) & \lTo & H^{n-3}(\bd V) & \lTo^0\\
\qquad&& \dTo^{=} && \dTo^\cong && \dTo^\cong  && \dTo^{=} && \dTo^{=} &
\\
\qquad& \lTo & 0 & \lTo & \Q & \lTo & \Q & \lTo & H^{n-3}(F) & \lTo &
H^{n-3}(\bd
V) & \lTo^0\\
\end{diagram}
}

because $F$ is an $(n-2)$-manifold with boundary and $\bd V$ is a closed
$(n-2)$-manifold. Since
the map $\Q\to\Q$ must be an isomorphism, so must be the map $H^{n-3}(\bd V)\to H^{n-3}(F)$. Therefore, $H^k(F)\cong H^k(\bd V$) for $k<n-2$.  
\indent
Now turning to the long exact sequences of the pairs, the inclusion $(V,F)\hookrightarrow (V,\bd V)$ and naturality give a commutative diagram 
\begin{diagram}
H^{k}(F) & \lTo & H^{k}(V) & \lTo &H^{k}(V,F) & \lTo & H^{k-1}(F) & \lTo & H^{k-1}(V)\\
\dTo && \dTo && \dTo && \dTo && \dTo \\
H^{k}(\bd V) & \lTo & H^{k}(V) & \lTo &H^{k}(V,\bd V) & \lTo & H^{k-1}(\bd V) & \lTo &
H^{k-1}(V).\\
\end{diagram} 

By the five-lemma, $H^k(V,F)\cong H^k(V,\bd V)$ for $k<n-2$. For $k=n-2$, we
can use the facts that $H^{n-2}(F)=0$ and $H^{n-1}(V)=0$, since $F$ is an
$(n-2)$-manifold with boundary and $V$ is an $(n-1)$-manifold with boundary, and
that $H^{n-2}(\bd V) \to H^{n-1}(V,\bd V)$ is an
isomorphism $\Q\to\Q$ for
similar reasons. These allow us to extract the following commutative
diagram near the top of the sequence:
\begin{diagram}
0 & \lTo & H^{n-2}(V) & \lTo &H^{n-2}(V,F) & \lTo & H^{n-3}(F) & \lTo & H^{n-3}(V)\\
&& \dTo && \dTo && \dTo && \dTo\\
0 & \lTo & H^{n-2}(V) & \lTo &H^{n-2}(V,\bd V) & \lTo & H^{n-3}(\bd V) &
\lTo & H^{n-3}(V).\\
\end{diagram} 
Again using the five-lemma, we conclude that $H^k(V,F)\cong H^k(V,\bd V)$
for $k<n-1$.

\item For $q=0$, the statement is obvious as all of the spaces are connected. Otherwise,
from part (1), $H_q(V)\cong H_q(Y)$, $0\leq q<n-1$, and from the perfect linking pairings
$L'$ and $L''$ (see Section \ref{S: linking numbers}, below), these are dually
paired to
$H_{n-p-1}(Y,Z)$ and
$H_{n-p-1}(V,F)$, respectively, for $0<q<n-1$. From the induced perfect rational pairings
we get the statement of part (2).

\end{enumerate}
\end{proof}

Returning to the maps $d_i$; $i=1,2$; it follows from the construction of the covering and the action of the covering translation, $t$, that the maps can be written as 
\begin{align*}
d_i (\alpha \otimes 1)&=i_{-*}(\alpha)\otimes t -i_{+*}(\alpha)\otimes 1 \\
 &= t(i_{-*}(\alpha)\otimes 1) -i_{+*}(\alpha)\otimes 1, 
\end{align*}
where $i_{\pm}$ correspond to the identification maps of $(V,F)$ to $(V_{\pm},F_{\pm})$ and $\alpha\in H_q(V;\Q)$ or $H_q(V, F;\Q)$ according to whether $i=1$ or $2$.

\subsubsection{Linking Numbers}\label{S: linking numbers}

We now turn to the linking pairings on these homology groups. Let $p=n-1-q$. There are
perfect (modulo torsion) pairings
\begin{equation}\label{E:L'}
L':H_p(V, F)\otimes H_q(Y)\to\Z
\end{equation}
\begin{equation}\label{E:L''}
L'':H_p(Y, Z)\otimes H_q(V)\to\Z.
\end{equation}
These are the usual geometric linking pairings which are induced, via some
isomorphisms, by the classical Lefschetz dual intersection pairing $H_i(S^n,
Z)\otimes H_{n-i}(S^n-Z)\to \Z$, for $Z\subset S^n$. 
In particular, given the disk $D^n$ and a 
closed subpolyhedron  $B\subset D^n$ which meets $\bd D^n=S^{n-1}$ regularly,
let $A=B\cup \bd D^n=B\cup
S^{n-1}$. We also think of $D^n=D^n_+$ as the top hemisphere of
$S^n=D^n_+\cup_{S^{n-1}}D^n_-$.  Then a linking
pairing
\begin{equation}\label{E: link 1}
L':H_i(B, B\cap S^{n-1})\times H_{n-i-1}(D^n-B)\to\Z,
\end{equation}
$0<i<n-1$, can be defined by applying the following isomorphisms and then
applying the intersection pairing:
{\small
\begin{align*}
H_i(B, B\cap S^{n-1})&\cong H_i(B\cup S^{n-1}, S^{n-1}) &&
\text{by excision}\\
&=H_i(A, S^{n-1})&&\text{by definition of $A$}\\
&\cong H_i(A) && \text{by the long exact sequence of $(A, S^{n-1})$}\\
&\cong H_{i+1}(D^n,A) && \text{by the long exact sequence of $(D^n, A)$}\\
&\cong H_{i+1}(S^n, A\cup D^n_-) && \text{by excision},
\end{align*}}
and
{\small
\begin{align*}
H_{n-i-1}(D^n-B)&\cong H_{n-i-1}(D^n-A)& &\text{by the homotopy
equivalence}\\
&&&\text{of  $D^n-A$ and $D^n-B$}\\
&= H_{n-i-1}(S^n-A\cup D^n_-)& &\text{since $D^n-A=S^n-A\cup D^n_-$}. 
\end{align*}}
The linking pairing
\begin{equation*}
L'':H_i(D^n-B, (D^n-B)\cap S^{n-1})\times H_{n-1-i}(B)\to\Z
\end{equation*} 
can be obtained by considering an open regular neighborhood, $N$,
of
$B$ in $D^n$. Then $N$ deformation retracts to $B$ and
$(D^n-B, (D^n-B)\cap S^{n-1})$
deformation retracts to $(D^n-N,(D^n-N)\cap S^{n-1})$. So, if $U=D^n-N$, then 
$H_*(B)\cong H_*(D^n-U)$ and $ H_*(D^n-B, (D^n-B)\cap S^{n-1})\cong H_*(U,
U\cap S^{n-1})$. Then we can apply $L'$ with $U$ in place of $B$ in equation
\eqref{E: link 1}.

See \cite[Appendix]{GBF} for more details on the construction of these linking
pairings.

By taking tensor products, these pairings extend to perfect pairings from the rational
homology groups to $\Q$.
Let $\{\alpha^p_i\}$, $\{\beta^q_i\}$, $\{\gamma^p_i\}$, and $\{\delta^q_i\}$
represent dual bases for $H_p(V, F)$, $H_q(Y)$, $H_p(Y, Z)$, and
$H_q(V)$, all modulo torsion, so that
\begin{equation}\label{E:dualbasis}
L'(\alpha^p_i \otimes \beta^q_j)=L''(\gamma^p_i \otimes \delta^q_j)=\delta_{ij}.
\end{equation}
These collections also form bases then for the rational homology groups that result by tensoring with $\Q$, and the relations \eqref{E:dualbasis} hold under the induced perfect rational pairing.

Given $r\in H_p(V, F;\Q)$ and $s\in H_q(V;\Q)$, we also have the relation
\begin{equation}\label{E:commutivity}
L'(r \otimes i_{-*}(s))=L''(i_{+*}(r) \otimes s).
\end{equation}
This can be seen as follows: we can choose the inclusion maps $i_{\pm}:
(V,F)\hookrightarrow (Y,Z)$ as isotopies which push $V$ out along its collar
in one direction or the other. Then any chain representing $s$ gets pushed
into $Y$ under $i_{-}$ and the the linking form is the intersection of this
chain with a chain, $R$, representing the isomorphic image of $r$ in
$H_{p+1}(D^n,V\cup S^{n-1};\Q)$ (see \cite[Appendix]{GBF} ). The latter chain
can be taken as some chain in $D^n$ whose boundary, lying in $ V\cup
S^{n-1}$, is a chain representing $r$. Now, under
the isotopy which takes $V$ to $i_{+}(V)$ and $i_-(V)$ to $V$, the chain
representing $s$ gets pushed back into $V$ and $R$ gets pushed into a chain
in $D^n$ whose boundary, lying in $Y\cup S^{n-1}$, is $i_+$ of the chain
representing $r$. In particular, this latter chain
represents $i_{+*}(r)\in H_p(Y, Z)$. Thus this isotopy induces maps which
take $i_{-*}(s)$ to $s$ and $r$ to $i_{+*}(r)$, but since the geometric
relationship between the chains is unaffected by the isotopy, the
intersection
number is unaffected. The formula then follows
immediately using the definitions of $L'$ and $L''$ as geometric linking pairings
(again, see \cite[Appendix]{GBF} for more details).
Similarly, we get
\begin{equation}\label{E:commutivity2}
L'(r \otimes i_{+*}(s))=L''(i_{-*}(r) \otimes s).
\end{equation}

The final property of linking numbers which we will need is that given $r$ and $s$ as above
\begin{align}\label{E:link int}
L'(r\otimes i_{-*}(s))-L'(r\otimes i_{+*}(s))&=r\cap s\\
L''(i_{-*}(r)\otimes s)-L''(i_{+*}(r)\otimes s)&=r\cap s, \label{E:link int2}
\end{align}
where $r\cap s$ is the intersection pairing on $V$. The proof is analogous
to that in the usual case \cite[p. 542]{L66}.

\subsubsection{The proof of Theorem \ref{T: disk knot nec con}}

We can now complete the proof of the theorem: With the bases
$\{\alpha^p_i\}$,
$\{\beta^q_i\}$, $\{\gamma^p_i\}$, and $\{\delta^q_i\}$ as above,
$\{\alpha_i\otimes 1\}$, etc., give bases of $H_p(V,
F;\Q)\otimes_{\Q}\Gamma$, etc. Let
\begin{align*}
i_{+*}(\delta^q_j)&=\sum_i \lambda_{ij}^q \beta_i^q\\
i_{-*}(\delta^q_j)&=\sum_i \sigma_{ij}^q \beta_i^q\\
i_{+*}(\alpha^q_j)&=\sum_i \mu_{ij}^q \gamma_i^q\\
i_{-*}(\alpha^q_j)&=\sum_i \tau_{ij}^q \gamma_i^q.
\end{align*}
Note that the $\lambda$, $\sigma$, $\mu$, and $\tau$ will all be integers (by the chain map
interpretation of $i_{\pm}$ and the fact that the $\alpha$ and $\delta$ were initially chosen as generators of the torsion free parts of the appropriate integral homology groups).
Then
\begin{align*}
d_1(\delta^q_j\otimes 1)&=\sum_i(t\sigma_{ij}^q-\lambda_{ij}^q)(\beta_i^q\otimes 1)\\
d_2(\alpha^q_j\otimes 1)&=\sum_i(t\tau_{ij}^q-\mu_{ij}^q)(\gamma_i^q\otimes 1),
\end{align*}
and we obtain presentation matrices
\begin{align*}
P^q_1(t)&=(t\sigma_{ij}^q-\lambda_{ij}^q)\\
P^q_2(t)&=(t\tau_{ij}^q-\mu_{ij}^q)
\end{align*}
for $H_q(\td{C};\Q)$ and $ H_q(\td{C},\td{X};\Q)$. As we have already seen, these matrices will be square for $0<q<n-1$ by Proposition \ref{P:square matrix}.

Applying the perfect linking pairing gives
\begin{align*}
L'(\alpha_k^p\otimes i_{+*}(\delta^q_j))&=\sum_i \lambda_{ij}^q L'(\alpha_k^p\otimes \beta_i^q)=\lambda_{kj}^q\\
L'(\alpha_k^p\otimes i_{-*}(\delta^q_j))&=\sum_i \sigma_{ij}^q L'(\alpha_k^p\otimes \beta_i^q)=\sigma_{kj}^q\\
L''(i_{+*}(\alpha^q_j)\otimes \delta_k^p)& =\sum_i \mu_{ij}^q L''(\gamma_i^q\otimes \delta_k^p)=\mu_{kj}^q\\
L''(i_{-*}(\alpha^q_j)\otimes\delta_k^p )&=\sum_i \tau_{ij}^q L''(\gamma_i^q\otimes \delta_k^p)=\tau_{kj}^q,
\end{align*}
and so by \eqref{E:commutivity} and \eqref{E:commutivity2}, we have $\sigma^q_{jk}=\mu^p_{kj}$
and $\lambda^q_{jk}=\tau^p_{kj}$. This implies that $P^q_1(t)=- tP^p_2(t^{-1})'$, where $'$ indicates transpose. Further, 
\begin{align*}
P^q_1(1)&=(\sigma_{ij}^q-\lambda_{ij}^q)=(-\alpha^p_i\cap \delta^q_j)\\
P^q_2(1)&=(\tau_{ij}^q-\mu_{ij}^q)=(-\alpha^q_i \cap \delta^p_j)
\end{align*}
by \eqref{E:link int} and \eqref{E:link int2}. Since this is the (modulo
torsion) intersection pairing $\cap:H_p(V)\otimes H_q(V,F)\cong
H_p(V)\otimes H_q(V,\bd V)\to \Z$, which is non-singular modulo torsion,
these matrices are non-singular which implies that the maps $d_i$ of the
Mayer-Vietoris sequence \eqref{E:MV1} are injective as claimed above. In
fact, as the matrix of a perfect intersection pairing over $\Z$ of the free
summands of the relevant integral homology modules, the matrix is unimodular
with
determinant $\pm 1$. 
	
With only minor modifications, the conclusion of the theorem is now obtained
just as in \cite[\S 2.8]{L66} by looking at the $i^{th}$ order minors of
$P^q_1(t)$ and $P^p_2(t)$ and applying the properties of modules over
principal ideal domains. In particular, if $\Delta_i^q$ and $\bar
\Delta_i^p$ are the greatest
common divisors of the $i$th order minors of $P^q_1(t)$ and $P^p_2(t)$,
respectively, then they are elements of $\Lambda$, and $\lambda_q^i\sim
\Delta_i^q/\Delta_{i+1}^q$ and
$\mu_p^i\sim \bar\Delta_i^p/\bar\Delta_{i+1}^p$ in $\Lambda$. Furthermore,
by the properties of $P^q_1(t)$ and $P^p_2(t)$ proven above,
$\Delta_i^q(1)=\pm1$, $\bar
\Delta_i^p(1)=\pm1$, and $\Delta_i^q(t)\sim\bar
\Delta_i^p(t^{-1})$ in $\Lambda$. These imply that $\lambda_q^i(1)=\pm1$,
$\mu_p^i(1)=\pm1$, and  $\lambda_q^i(t)\sim \mu_p^i(t^{-1})$ in $\Lambda$.\hfill$\square$

\subsection{Some corollaries; Definition of Alexander
polynomials}\label{S: cors and poly def}
\begin{corollary}
With the notations as above, if the boundary slice knot $k$ is $j$-simple, meaning that
$S^{n-1}-k$ has the
homotopy of a circle for dimensions less than or equal to $j$, then for $0<i<j+1$ and $n-j-3<i<n-1$,
\begin{enumerate}
\item $\lambda_q^i(1)=\pm1$
\item $\lambda_q^i(t)\sim \lambda^i_{n-q-1}(t^{-1})$
\end{enumerate}  
\end{corollary}

\begin{proof} 
By \cite[p.14]{L65} the simplicity condition implies that we can modify the Seifert
surface, $F$, to be $j$-connected without changing $F$ near its boundary. We can then
use this
Seifert surface to redefine the extension of the map $f$ of Proposition \ref{P:Seifert} on $S^{n-1}$
so that it
yields this Seifert surface on the sphere. We then extend to the interior of $D^n-L$ as in
Proposition \ref{P:Seifert} to get a new $V$ with $\bd V=F\cup L$.

Since $\pi_i(F)=0$ for $0<i<j+1$, we get $H_i(F)=0$ in the same range, and so $H_i(V;\Q)\cong
H_i(V,F;\Q)$ by the long exact sequence of the pair. This isomorphism also holds for $n-j-3<i<n-1$.
In fact, by Lefschetz duality, $\td H_i(F;\Q)\cong H_{n-i-2}(F,\bd F;\Q)\cong
H_{n-i-2}(F,S^{n-3};\Q)$, since $F$ is $(n-2)$-dimensional and we are using field
coefficients. Thus
$H_i(F, S^{n-3};\Q)=0$ for $n-j-3<i<n-2$ and $H_i(F)=0$, $n-j-3<i<n-3$ by the long exact sequence of
the pair $(F,\bd F)$. $H_{n-2}(F)=0$ since $F$ is an $(n-2)$-manifold with boundary, and,
since the
top
of the sequence is (suppressing the coefficients)

{\small
\begin{diagram}
H_{n-2}(F)&\rTo& H_{n-2}(F,
S^{n-3}) &\rTo& H_{n-3}(S^{n-3})&\rTo& H_{n-3}(F)&\rTo&  H_{n-3}(F,
S^{n-3})\\
\dTo^{\cong}& &\dTo^{\cong} &  &\dTo^{\cong} &     &\dTo^{\cong} &
&\dTo^{\cong}       \\
0 &\rTo&  \Q  &\rTo&  \Q  &\rTo&  ? &\rTo&  0&,
\end{diagram}}

\noindent we must have that $H_{n-3}(F)=0$ also. Again by the long exact
sequence of the pair, we get $H_i(V;\Q)\cong H_i(V,F;\Q)$, $n-j-3<i<n-1$.

So, in this range, the perfect linking pairing $L'$ can be defined 
\begin{equation*}
L':H_p(V;\Q)\otimes H_q(Y;\Q)\to\Q
\end{equation*} 
using the isomorphisms $H_i(V)\cong H_i(V,F)$. From here, the proof follows as above and as in \cite{L66} using the dual bases of $H_p(V;\Q)$ and $H_q(Y;\Q)$ in the relevant dimensions.   
\end{proof}

\begin{corollary}\label{C: nec con}
If we define the \emph{Alexander polynomials} $\lambda_q(t)$ and $\mu_q(t)$
as the primitive polynomials in $\Lambda$ determined up to similarity class
by the determinants of the square presentation matrices of $H_i(\tilde C)$
and $H_i(\td C, \td X)$ as $\Gamma$-modules, then
\begin{enumerate}
\item $\lambda_q(1)=\pm1$ and $\mu_q(1)=\pm 1$ 
\item $\lambda_q(t)\sim\mu_{n-p-1}(t^{-1})$.
\end{enumerate}
\end{corollary}
\begin{proof}
From our earlier definitions, $\lambda^q=\prod \lambda^q_j$ and $\mu^q=\prod
\mu^q_j$. The corollary now follows immediately.
\end{proof}

\nopagebreak

%% file: alexb.tex

\nopagebreak

\begin{corollary}\label{C:division}
Let $\nu_i(t)$, $\lambda_i(t)$, and $\mu_i(t)$, $0<i<n-2$, be the
Alexander
polynomials corresponding to $H_i(\td X)$, $H_i(\td C)$, and $H_i(\td C,
\td X)$, respectively. 
Then 
\begin{equation}\label{E: alt prod}
\prod_{i>0}\frac{\mu_{2i-1}(t)\nu_{2i-1}(t)\lambda_{2i}(t)}{\mu_{2i}(t)\nu_{2i}(t)\lambda_{2i-1}(t)}=1,
\end{equation}
where, for this formula only, we define the polynomials to be $1$ for $i>n-2$.   
\end{corollary}
\begin{proof}

The Alexander polynomials are given by the determinants of the
presentation matrices of the terms of the exact sequence of $\Gamma$-modules
\begin{equation*}
\begin{CD}
@>>> \td H_i(\td X;\Q) @>>> \td H_i(\td C;\Q) @>>> \td H_i(\td C, \td X;\Q) @>>>.
\end{CD}
\end{equation*}
We know that each term is finitely generated as a $\Gamma$-module, so the
corollary follows immediately from Proposition \ref{P:alt. poly.}, the
triviality
of $\td H_0(\td X; \Q)$ (since $\td X$ is connected), and the triviality
of
$H_{n-2}(\td X; \Q)$ (by classical knot theory).
\end{proof}

\begin{corollary}
With the notation above, $\lambda_{n-2}(t)$ divides $\lambda_1(t^{-1})$.
\end{corollary}
\begin{proof}
From the proof of the last corollary and Corollary \ref{C: kern},
$\lambda_{n-2}(t)$
divides $\mu_{n-2}(t)$, but $\mu_{n-2}(t)\sim \lambda_1(t^{-1})$.
\end{proof} 

We can also use these methods to obtain the well-known fact:
\begin{corollary}\label{C: slice knots}
A classical slice 1-knot ($S^1\subset S^3$) has Alexander polynomial of
the form
$\nu_1(t)\sim p(t)p(t^{-1})$.
\end{corollary}
\begin{proof}
We take $n=4$ for our disk knot pair, so that the boundary slice knot will
be a knotted $S^1$ in $S^3$. Then the only non-trivial Alexander
polynomials are $\nu_1(t)$, $\lambda_1(t)$, $\lambda_2(t)$, $\mu_1(t)$, and $\mu_2(t)$. From
Corollary \ref{C:division}, 
\[\nu_1(t)\sim\frac{\mu_2(t)\lambda_1(t)}{\mu_1(t)\lambda_2(t)}\sim
\frac{\lambda_1(t^{-1})\lambda_1(t)}{\lambda_2(t^{-1})\lambda_2(t)}.\]
From here we can proceed more or less as in \cite{M66}: Let $d(t)$ be the
greatest common divisor of $\lambda_1(t)$ and $\lambda_2(t)$ so that
$\lambda_1(t)=d(t)a(t)$, $\lambda_2(t)=d(t)b(t)$, and $a(t)$ and $b(t)$ are relatively
prime. Then 
\[\nu_1(t)\sim\frac{d(t)d(t^{-1})a(t)a(t^{-1})}{d(t)d(t^{-1})b(t)b(t^{-1})}=
\frac{a(t)a(t^{-1})}{b(t)b(t^{-1})}.\]
Similarly, now let $c(t)$ be the greatest common divisor of $a(t^{-1})$ and
$b(t)$ so that $a(t^{-1})=p(t^{-1})c(t)$ and $b(t)=q(t)c(t)$. Then
\[\nu_1(t)=\frac{p(t)p(t^{-1})}{q(t)q(t^{-1})},\]
and the numerator and denominator are now relatively prime. But $\nu_1(t)$ is
actually a polynomial so $q(t)\sim 1$ and $\nu_1(t)\sim
p(t)p(t^{-1})$. 
\end{proof}

\subsection{Realization of given polynomials}\label{S: Realization}

In this section, we obtain results on the realization of knots with
prescribed Alexander polynomials. The construction of a knot $D^2\subset
D^4$ with a given polynomial is done by hand to get a feel for the
geometric concepts involved. This lays the foundation for realization
theorems in higher dimensions. 

Throughout this section, we continue to use $\lambda_q$, $\mu_q$, and
$\nu_q$ as defined in Corollary \ref{C:division}.

\subsubsection[Realizing $\lambda_i$ for $D^2\subset D^4$]{Realizing
{\mv $\lambda_i$} for {\mv $D^2\subset D^4$}}\label{S: D2 in
D4}

\begin{theorem}
Given any polynomial $p(t)\in\Lambda$ such that $p(1)=\pm1$,
there exists a knotted $D^2$ in $D^4$ with $\lambda_1(t)\sim p(t)$ and
$\lambda_2(t)\sim1$.  
\end{theorem}
\begin{proof}
For definiteness, let us normalize $p(t)$ so that
$p(t)=\sum_{i=0}^m a_i t^i$, $p(1)=1$, and $p(0)\neq0$. We
will construct a knotted disk with $H_2(\td C)=0$ and $H_1(\td C)\cong
\Lambda/p(t)$. 

We begin by embedding a 2-disk, $L$, in $S^1\times D^3$,
so that, in a neighborhood of a boundary point which is homeomorphic to the
half-space $\mathbb{R}^{4+}$, $L$ is embedded as a standard disk. In other
words, $\bd D^2=S^1$ is an unknotted circle within a neighborhood of a point
of $\bd (S^1\times D^3)=S^1\times S^2$, $\text{int}(L)$ lies in
$\text{int}(S^1\times D^3)$, and $(L, \bd L)$ is null-homotopic in ($S^1\times
D^3, \bd(S^1\times D^3))$. We also let $\bd L$ bound a disk $F$ in $S^1\times
S^2$ so that $F\cup L$ bounds a manifold $V$ homeomorphic to a standard $D^3$
such that int$(V)$ lies in int$(S^1\times D^3)$. Let $C_0=S^1\times D^3-L$,
and let $\td C_0$ be the infinite cyclic covering associated with the
kernel of the homomorphism $\pi_1(\td C_0)\to\Z$ defined by intersection
number with $V$. Forming the infinite cyclic cover by cutting along $V$, it
is clear that $H_2(\td C_0)=0$ and that $H_1(\td C_0)\cong \Lambda$, where we can take as
generator,
$\alpha$, the lift of a circle representing a generator of $\pi_1(S^1\times
D^3)$ and which does not intersect $F$.

We will prove the following lemma below: 
\begin{lemma}\label{L:embedding lemma} 
There exists an
embedding $f: S^1\hookrightarrow S^1\times S^2-\bd L$ which lifts to an embedding
$g:S^1\hookrightarrow \td C_0$ which represents the element $\lambda(t)\alpha\in
H_1(\td C_0)$. Furthermore, $f$ can be chosen isotopic to the standard
embedding which
takes $S^1$ to $S^1\times x_0$ for some $x_0\in S^2$. \end{lemma}

Now let $S=f(S^1) \in S^1\times S^2=\bd(S^1\times D^3)$. We will attach a
2-handle along $S$. In particular,
there is a neighborhood
$S\times D^2$ of $S$ in $\bd(S^1\times D^3)$ which we identify with half of
the boundary $\bd(I^2\times I^2)=(S^1\times D^2)\cup (D^2\times S^1)$. If $H$
denotes the handle, then $(S^1\times D^3)\cup_{S\times D^2} H\cong D^4$, and
we claim that $L$, which is now knotted in $D^4$, is the desired knotted
disk.

Still assuming the lemma, it remains to show only that we get the desired homology of the
cover. Notice that $S$ lifts to an infinite number of disjoint embeddings, $S_i$, in $\td C_0$
which correspond to $t^i\lambda(t)\alpha$, $-\infty<i<\infty$. If we attach an infinite number
of handles sewn along $S_i \times D^2$, we will obtain an infinite cyclic covering of $D^4-L$,
which we denote by $\td C$. That $\td C$ has the desired homology follows from the reduced
Mayer-Vietoris sequence:
\begin{equation*}
\begin{diagram}
&\rTo&\oplus_{i=-\infty}^{\infty}H_2(S^1) &\rTo&
(\oplus_{i=-\infty}^{\infty} H_2(D^4))\oplus H_2(\td C_0) &\rTo & H_2(\td
C) 
\end{diagram}
\end{equation*}
\begin{equation*}
\begin{diagram}
&\rTo& \oplus_{i=-\infty}^{\infty}H_1(S^1) &\rTo^{d}&
(\oplus_{i=-\infty}^{\infty} H_1(D^4))\oplus
H_1(\td C_0)& \rTo& H_1(\td C) 
\end{diagram}
\end{equation*}
\begin{equation*}
\begin{diagram}
&\rTo& \oplus_{i=-\infty}^{\infty}\td H_0(S^1) &\rTo &
(\oplus_{i=-\infty}^{\infty} \td
H_0(D^4))\oplus\td H_0(\td C_0)& \rTo & \td H_0(\td C)&\rTo.
\end{diagram}
\end{equation*}
The first two terms and the last are zero, as are $H_i(D^4)$, $i=1,2$, and
$\td
H_0(\td C_0)$. The map $\oplus_{i=-\infty}^{\infty}\td H_0(S^1) \to
\oplus_{i=-\infty}^{\infty} \td H_0(D^4)$ is an isomorphism
$\Lambda\to\Lambda$, and $\oplus_{i=-\infty}^{\infty}H_1(S^1)$ and $H_1(\td
C_0)$ are both isomorphic to $\Lambda$. Since we know that the generators
of $H_1(S_i)$ map
onto the the generators $t^i\lambda(t)\alpha$, the map $d$ must be an
injection. Hence, we can conclude from this information that $H_2(\td C)=0$
and $H_1(\td C)\cong \Lambda/(\lambda(t)\Lambda)$, which is the desired
result.  
\end{proof}

\begin{proof}[Proof of Lemma \ref{L:embedding lemma}] 
We wish to embed a circle $S$
into $S^1\times S^2$ so that it will be isotopic to a standard circle and so
that
that some lifting
will represent
$\lambda(t)\alpha$, where $\lambda(t)=\sum_{i=0}^m a_i t^i$ and $\alpha$ is some
generator of $H_1(\td C_0)$. It is possible, and simpler for visualization
purposes,
to embed the circle into the standard solid torus $S^1\times D^2\subset
S^1\times S^2$ obtained by
removing
a neighborhood of some $S^1\times x$, $x\in S^2$. We are also free to take $L$ in the
theorem so that $\bd L$ is a circle concentric to a standard meridian
inside this solid torus.
Then $F$ can be taken as the disk which fills in this circle. Note that
$\alpha$ can be
taken as a lift of a longitude, $\ell$,  which does not intersect $F$. 

We will construct $S$ primarily by running around the boundary $S^1\times S^1$ of the
solid torus with ever-increasing meridional angle. To be precise, we begin by choosing
an orientation for the longitude, $\ell= S^1\times 0$, which does intersect $F$, so that
its
lifts will be arcs
running from $\td x$ to $t\td x$, where $\td x$ is the lift of a point of the longitude.
Now, choose a point $x_0$ which lies in $S^1\times S^1$ on the meridian concentric to
$\bd L$ and $F$. We begin by running an arc around the torus $|a_0|$ times, choosing the
direction to agree with the that of $\ell$ if $a_0>0$ or to disagree if
$a_0<0$, 
while the meridional angle increases slightly to avoid self
intersection. Then run the arc into the interior through $F$ in the direction 
of $\ell$ and then back out
to the boundary of the torus. It is clear that this can be done in such a way that the
radial retraction of the arc to the torus will continue to be an embedding with increasing
meridional angle. Now, follow the same procedure for each of the $a_i$, doing nothing
but the final step of crossing $F$ if $a_i=0$. Clearly we can choose the rate of
increase of the meridional angle so that we never complete a full cycle meridionally.
Lastly, after wrapping around the torus the $a_m$th time, we run the arc back through $F$ $m$
times against the direction of $\ell$ (i.e. so that it links with $\bd F$ $m$
more times but the total linking number will be 0), still with increasing meridional angle, 
and then connect it back to the starting point along a meridian. 

To see that $S$ is isotopic to the standard longitude $S^1\time 0$, first observe
that our construction allows us to isotop $S$ out to the torus $S^1\times S^1$. In the torus,
the \emph{homotopy type} of $S$ is $(1,1)$ since $\lambda(1)=1$ and by the
method of construction. Now by \cite[p. 25]{Rolf}, $S$ is ambient isotopic in the torus
to the standard representation of the $(1,1)$ homotopy class, and this ambient isotopy
can be extended to a
neighborhood $S^1\times S^1\times [-1,1]$ of the torus in $S^1\times S^2$ (indeed, just
perform the isotopy, itself, on $S^1\times S^1\times [-1,0]$ and then its reverse on 
$S^1\times S^1\times [0,1]$). But the standard representation of $(1,1)$ is clearly
ambient isotopic to the standard $S^1\times 0\subset S^1\times D^2$
displacing it radially in each
meridional disk.   

It also apparent from the construction that $S$ will lift to the proper element of
$H_1(\td C_0)$. In fact, by considering the usual cut and paste construction of the
infinite cyclic covering, $\td C_0$ looks like an infinite number of $S^1\times
D^4$'s glued together, and, as remarked above, the generator of $H_1(\td C_0)$
corresponds to the generator of one of these ``solid tori'' and projects down to a
generator of the homology of $H_1(S^1\times D^4)$ which does not intersect $V$ (or
hence, $F$). So, by this construction, if we lift the starting point of $S$ to its
covering point in $\td C_0$ corresponding to the $t^0$ copy of $S^1\times D^4$, then
$S$ lifts to an arc which runs around this $S^1\times D^4$, parallel to $\alpha$,  
$a_0$ times in the correct
direction, then crosses into the $t^1$ copy of $S^1\times D^4$ and runs around it $a_1$
times parallel to $t\alpha$ and so on. After finishing its circuits in the $t^m$ copy
of $S^1\times D^4$, it
returns straight back to its starting point. Evidently, this lift represents the
homology class $\lambda(t)\alpha$ as desired.
\end{proof}

\begin{corollary} The conditions $\lambda_2(1)=\pm 1$, $\lambda_1(1)=\pm
1$, and $\lambda_2(t)|\lambda_1(t^{-1})$
completely characterize all of the Alexander polynomials, $\lambda_i$, of
a
disk knot $D^2\subset D^4$ and hence of a  singularly knotted 2-sphere in
$S^4$. \end{corollary} 
\begin{proof} 
We know that these conditions are necessary. To show that they are sufficient, let $\lambda_1(t)=\lambda_2(t^{-1})r(t)$. Then we must
also have $r(1)=\pm 1$. It follows from the preceding theorem that we can find a
knotted $D^2\subset D^4$ whose Alexander polynomials are $r(t)$ and $1$ in
dimensions 1 and
2, respectively. Taking the cone on the boundary sphere pair gives a singular knot with
the same Alexander polynomial. We can also find a locally-flat knot whose first and
second Alexander polynomials are $\lambda_2(t^{-1})$ and $\lambda_2(t)$, respectively
\cite[\S4]{L66}. Then the knot sum of these two knots has the desired Alexander
polynomials since Alexander polynomials multiply under knot sum.  
\end{proof}

Note that for a knot $D^2\subset D^4$ we have now completely classified
\emph{all} of the Alexander polynomials, since $\mu_1(t)\sim
\lambda_2(t^{-1})$, $\mu_2\sim \lambda_1(t^{-1})$, and $\nu_1$ is
completely determined by Corollary \ref{C:division}.

%% file: alexc.tex

\subsubsection[Realizing $\lambda_i$ for $D^{n-2}\subset D^n$, $n\geq
5$]{Realizing {\mv$\lambda_i$} for {\mv$D^{n-2}\subset D^n$, $n\geq
5$}}\label{S: large
n}

We next turn to realizing the Alexander polynomials, $\lambda_i(t)$, for
$n$-disk knots, $n\geq
5$. Our arguments are split into two propositions. The first provides
the realizability directly for the lower dimensional polynomials. The
second provides the realizability in higher dimensions by constructing
knots with the appropriate dual polynomials $\mu_i(t^{-1})$ in the lower
dimensions. The first result has been shown already by Sumners
in his thesis using similar methods (see \cite{Sum}, \cite{Sum2}). The
second result on the higher
dimensional $\lambda_i$ was also shown there but by different
methods. Consequences of our specific construction will be used again in
the proof of Theorem \ref{T:all polys} below.

\begin{proposition}\label{P:lower dims}
Given a polynomial $p(t)$ such that $p(1)=\pm 1$ and
integers $q$ and $n$ such that $1\leq q\leq \frac{n-2}{2}$ and $n\geq 5$, there exists a knotted $D^{n-2}\subset D^n$ such that $\lambda_q(t)\sim p(t)$ and
$\lambda_i(t)\sim 1$ for $0<i<n-1$, $i\neq q$, where $\lambda_i(t)$ is the Alexander
polynomial corresponding to $H_i(\td C)$, $\td C$ the knotted disk
complement.
\end{proposition}
\begin{proof}
We can normalize $p(t)$ so that $p(1)=1$. It suffices to construct a disk
knot
such that $H_q(\td C)\cong\Lambda/(p(t))$ and $H_i(\td C)=0$ for
$0<i<n-1$, $i\neq q$. The proof is a variation of that of Levine
\cite{L66} for locally-flat knots $S^{n-2}\subset S^n$.

We begin by embedding an $(n-2)$-disk, $L$, in $S^q\times D^{n-q}$,
so that, in a neighborhood of a boundary point which is homeomorphic to the
half-space $\mathbb{R}^{n+}$, $L$ is embedded as a standard disk. In other
words, $\bd D^{n-2}=S^{n-3}$ is an unknotted sphere within a neighborhood of a point in $\bd (S^q\times D^{n-q})=S^q\times S^{n-q-1}$, $\text{int}(L)$ lies in
$\text{int}(S^q\times D^{n-q})$, and $(L, \bd L)$ is null-homotopic in ($S^q\times D^{n-q}, \bd(S^q\times D^{n-q}))$. We also let $\bd L$ bound a disk $F$ in $S^q\times
S^{n-q-1}$ so that $F\cup L$ bounds a manifold $V$, homeomorphic to a standard $D^{n-1}$,
such that int$(V)$ lies in int$(S^q\times D^{n-q})$. Let $C_0=S^q\times D^{n-q}-L$,
and let $\td C_0$ be the infinite cyclic covering associated with the
kernel of the homomorphism $\pi_1(\td C_0)\to\Z$ defined by intersection
number with $V$. Similarly, we have the covering $\td X_0$ of $X_0=S^q\times S^{n-q-1}-\bd L$. Forming the infinite cyclic covers by cutting and pasting along $V$, it is clear that 
\begin{align*}
\td H_i(\td C_0)&\cong 
\begin{cases}
\Lambda, &i=q\\
0, &i\neq q
\end{cases}\\
\td H_i(\td X_0)&\cong 
\begin{cases}
\Lambda, & i=q,n-q-1\\
0, & i\neq q,n-q-1.
\end{cases}
\end{align*}
In fact $H_q(\td C_0)\cong H_q(\td X_0)$, and we can take the lift of a
sphere
representing a generator of $\pi_q(S^q\times S^{n-q-1})\cong \Z$ which
does not intersect $F$ as a
$\Lambda$-module generator, $\alpha$, of both modules. ($\td X_0$ is the connected sum along the
boundaries of an infinite number of copies of $S^q\times
S^{n-q-1}-\text{\{the open neighborhood of a point\}}$, and $\td C_0$ is
the boundary connected sum of an infinite number of copies of $S^q\times
D^{n-q}$.)

The Hurewicz map 
\begin{equation*}
h_q:\pi_q(\td X_0)\to H_q(\td X_0)\cong H_q(\td C_0)
\end{equation*} 
is an epimorphism. This follows immediately if $q=1$ by the abelianization
map. For $q>1$, we note that $\pi_1(\td X_0)=0$ using the Van Kampen
theorem, and then the Hurewicz theorem applies. Since $\pi_q(\td X_0)$ is
isomorphic to a subgroup of $\pi_q(X_0)$, we can represent any element of
$H_q(\td C_0)$ by the lift of an embedded sphere in $X_0$ to $\td X_0$.
The embedding is possible since $2q<n-1$. In particular, we choose an
embedded sphere, $S$, in $X_0$ whose lift represents $p(t)\alpha$ in
$H_q(\td C_0)$. We will attach a handle to $S^q\times D^{n-q}$ along $S$
to create a new manifold which will turn out to be $D^n$. The image
of
$L$ under the modification is the desired knotted disk.

In particular, the dimensions are sufficient for us to embed $S^q\times D^{n-q-1}$ as a tubular neighborhood of $S$ in $\bd(S^q\times D^{n-q})$, and we identify this neighborhood with the first term of the boundary 
\begin{equation*}
\bd(D^n)=\bd(D^{q+1}\times D^{n-q-1})=(S^q\times D^{n-q-1})\cup_{S^q\times S^{n-q-2}} (D^{q+1}\times S^{n-q-2})
\end{equation*}
to form the new manifold $\Delta \cong (S^q\times D^{n-q})\cup_{S\times D^{n-q-1}} D^n$. On the boundary, this gives the surgery which transforms $S^q\times S^{n-q-1}$ into 
\begin{equation*}
\bd \Delta=(S^q\times S^{n-q-1}-S\times D^{n-q-1})\cup_{S\times S^{n-q-2}} (D^{q+1}\times S^{n-q-2}).
\end{equation*}
We first show that $\Delta$ is in fact isomorphic to $D^n$. 

Since $p(1)=1$, $S$ represents the generator of $\pi_q(S^q\times D^{n-q})\cong \Z$ and hence
the generator of $H_q(S^q\times D^{n-q})\cong \Z$. The reduced Mayer-Vietoris sequence
immediately gives us that $\td H_i(\Delta)=0$ for all $i$. On the boundary, following Levine
\cite[p.547]{L66}, we can choose $S$ isotopic, in $S^q\times S^{n-q-1}$, to the standard embedded $S^q\times x_0$, $x_0\in S^{n-q-1}$, provided $q<\frac{n-2}{2}$, and the modified boundary is then diffeomorphic to $S^{n-1}$. In fact, we can extend the isotopy on the boundary radially into $S^q\times D^{n-q}$. Then the standard $(q+1)$-handle attachment to $S^q\times D^{n-q}$ along a tubular neighborhood of $S^q\times x_0$ yields the n-disk. 

For $q=\frac{n-2}{2}$, we have $n\geq 6$, and we will show that $\Delta$ is a disk through an application of the h-cobordism theorem \cite{RS}. First,
it must be that dim$(S^q\times S^{n-q-1})=n-1$ is odd, and since $n\geq 6$, $S^q\times S^{n-q-1}$ is simply-connected. It then follows from simply-connected surgery theory (see \cite[IV.2.13]{BW1}) that $\td H_i(\bd \Delta)=0$ for $i\leq \frac{n-2}{2}$ and then from Poincare duality that 
\begin{equation*}
\td H_i(\bd \Delta)\cong 
\begin{cases}
Z, &i=n-1\\ 
0, &i\neq n.
\end{cases}
\end{equation*}
Furthermore, $\bd \Delta$ is simply connected by the Van Kampen theorem since $D^{q+1}\times S^{n-q-2}$ and $S^q\times S^{n-q-1}-S\times D^{n-q-1}$ are both simply-connected, the latter because $S^q\times S^{n-q-1}$ is simply connected and $S\times D^{n-q-1}$ is homotopy equivalent to a subset of codimension $>2$.
We now wish to show that $\bd \Delta$ is homotopy equivalent to $\Delta-B_{\epsilon}(x)$,
where $ B_{\epsilon}(x)$ is a small open ball neighborhood of a point, $x$, in int$(\Delta)$. $\Delta$ is also simply connected by an easy application of the Van-Kampen theorem, and therefore so is $\Delta-B_{\epsilon}(x)$, which is homotopy equivalent to $\Delta-\{\text{a point}\}$.
Since $\td H_i(\Delta-B_{\epsilon}(x))=0$ except in dimension $n-1$ (by an easy long exact sequence argument for the pair $(\Delta, \Delta-B_{\epsilon}(x))$), the inclusion $\imath: \bd \Delta \to \Delta-B_{\epsilon}(x)$ induces an isomorphism of $H_i$, $i\neq n-1$. It also induces the isomorphism in dimension $n-1$ from the long exact sequence of the pair, since 
\begin{equation*}
H_n(\Delta-B_{\epsilon}(x), \bd \Delta)\cong H^0(\Delta-B_{\epsilon}(x),\bd \bar B_{\epsilon}(x))=0,
\end{equation*}
using Lefschetz duality, and 
\begin{align*} 
H_{n-1}(\Delta-B_{\epsilon}(x), \bd \Delta)&\cong
H^1(\Delta-B_{\epsilon}(x),\bd B_{\epsilon}(x))\\
& \cong
\text{Hom}(H_1(\Delta-B_{\epsilon}(x),\bd B_{\epsilon}(x)),\Z)=0,
\end{align*}
using Lefschetz duality, the universal coefficient theorem, the reduced long exact sequence of
the pair, and the simple-connectivity of each term of the pair. Therefore, $\imath: \bd \Delta
\to \Delta-B_{\epsilon}(x)$ is a homotopy equivalence by the Whitehead theorem, since it is a
homology equivalence of simply connected spaces. That $\imath: \bd \bar B_{\epsilon}(x)\cong
S^{n-1} \to \Delta-B_{\epsilon}(x)$ is a homotopy equivalence follows similarly, and the
h-cobordism theorem applies to tell us that $\Delta-B_{\epsilon}(x)\cong S^{n-1}\times I$.
Filling $B_{\epsilon}(x)$ back in gives us that $\Delta\cong D^n$ as claimed.

Finally, letting $C$ denote $\Delta-L$, we need to show that $\td C$ has the desired homology modules. But we can form $\td C$ by attaching an infinite number of handles to $\td C_0$, attached along the infinite number of lifts of $S$ which represent the homology elements $t^i p(t)\alpha$ obtained from $p(t)\alpha$ by the actions of the covering transformations. Then it is immediate from the Mayer-Vietoris sequence that
\begin{equation*}
\td H_i(\td C)\cong
\begin{cases}
\Lambda/(p(t)), &i=q\\
0, &i\neq q,
\end{cases}
\end{equation*}
which completes the proof of the proposition. 
\end{proof}

\begin{proposition}\label{P:upper dims}
Given a polynomial $p(t)$ such that $p(1)=\pm 1$ and
integers $q$ and $n$ such that $\frac{n-2}{2}<q<n-2$ and $n\geq 5$, there exists a knotted $D^{n-2}\subset D^n$  such that $\lambda_q(t)\sim p(t)$ and
$\lambda_i(t)\sim 1$ for $0<i<n-1$, $i\neq q$, where $\lambda_i(t)$ is the Alexander
polynomial corresponding to $H_i(\td C)$, $\td C$ the knotted disk
complement.
\end{proposition}
\begin{proof}

It suffices to construct a disk
knot such that $H_q(\td C)\cong\Lambda/(p(t))$ and
$H_i(\td C)=0$ for $0<i<n-1$, $i\neq q$. In fact, letting $X$ denote
$C\cap \bd D^n$ and $p=n-q-1$ (so that $1<p<\frac{n}{2}$), we
construct a disk knot that

\begin{equation*}
H_i(\td C, \td X)\cong 
\begin{cases}
\Lambda/(p(t^{-1})), &i=p\\
0, & 0<i<n-1, i\neq p,
\end{cases}
\end{equation*}
which will suffice since the Alexander polynomials corresponding to
$H_i(\td C)$ and $H_{n-i-1}(\td C,\td X)$ are related by $\lambda_i(t)\sim
\mu_{n-i-1}(t^{-1})$ according to Theorem \ref{T: disk knot nec con}. We
normalize $p(t)$ so that
$p(1)=1$.

We begin by embedding an $(n-2)$-disk, $L$, in $D^p\times S^{n-p}$, so that,
in a neighborhood of a boundary point which is homeomorphic to the
half-space $\mathbb{R}^{n+}$, $L$ is embedded as a standard disk. In other
words, $\bd D^{n-2}=S^{n-3}$ is an unknotted sphere within a neighborhood
of a point in $\bd (D^p\times S^{n-p})=S^{p-1}\times S^{n-p}$,
$\text{int}(L)$ lies in $\text{int}(D^p\times S^{n-p})$, and $(L, \bd L)$
is null-homotopic in ($D^p\times S^{n-p}, \bd(D^p\times S^{n-p}))$. We
also let $\bd L$ bound a disk $F$ in $S^{p-1}\times S^{n-p}$ so that
$F\cup L$ bounds a manifold, $V$, homeomorphic to a standard $D^{n-1}$,
such that int$(V)$ lies in int$(D^p\times S^{n-p})$. Let $C_0=D^p\times
S^{n-p}-L$, and let $\td C_0$ be the infinite cyclic covering associated
with the kernel of the homomorphism $\pi_1(\td C_0)\to\Z$ defined by
intersection number with $V$. Similarly, we have the covering $\td X_0$ of
$X_0=S^{p-1}\times S^{n-p}-\bd L$. Forming the infinite cyclic covers by
cutting and pasting along $V$, it is clear that
\begin{align*}
\td H_i(\td C_0)&\cong 
\begin{cases}
\Lambda, &i=n-p\\
0, &i\neq n-p
\end{cases}\\
\td H_i(\td X_0)&\cong 
\begin{cases}
\Lambda, & i=p-1,n-p\\
0, & i\neq p-1,n-p.
\end{cases}
\end{align*}
In fact, $H_{n-p}(\td C_0)\cong H_{n-p}(\td X_0)$, and we can take as
the generator, $\alpha$, of both modules the lift
of a sphere representing a generator
of $H_{n-p}(S^{p-1}\times S^{n-p})\cong \Z$ and which does not intersect $F$.
($\td X_0$ is the connected sum along the boundaries of a countably 
infinite number
of copies of $D^p\times S^{n-p}-\text{\{the open neighborhood of a
point\}}$, and $\td C_0$ is isomorphic to the boundary connected sum of an
infinite number of copies of $D^p\times S^{n-p}$.) By the long exact
sequence of the pair,
\begin{equation*}
H_i(\td C_0,\td X_0)\cong
\begin{cases}
\Lambda, &i=p\\
0, & i\neq p.
\end{cases}
\end{equation*}
Now for a lemma:

\begin{lemma}
The Hurewicz map $h_p:\pi_p(\td C_0, \td X_0)\to H_p(\td C_0, \td X_0)$ is an epimorphism.
\end{lemma}
\begin{proof}
We note first that $\td C_0$ is simply connected because it is the universal abelian cover of
$C_0$ whose fundamental group is $\pi_1(C_0)=\Z$. This last statement is true because we can
decompose $C_0$ into $N-L$, where $N$ is the contractible neighborhood of the boundary in
which we have embedded $L$, and $D^p\times S^{n-p}-N$. The latter is homotopy equivalent to
$D^p\times S^{n-p}$, and $\pi_1(D^p\times S^{n-p})=0$ due to the range of $p$. $N-L$ is
homotopy equivalent to the complement of the trivial sphere pair $(S^n, S^{n-2})$, and so
$\pi_1(N-L)=\Z$. Since $n\geq 5$, an easy application of the Van-Kampen theorem proves
the claim.

Then, using our knowledge of the homology of $\td C_0$ and the Hurewicz theorem \cite[\S
VII.10]{BR}, $\pi_i(\td C_0)=0$, $i<n-p$. This implies by the long exact homotopy sequence
that $\pi_i(\td C_0, \td X_0)\cong \pi_{i-1}(\td X_0)$, $1<i<n-p$. Furthermore, $H_i(\td C_0,
\td X_0)\cong H_{i-1}(\td X_0)$, $1<i<n-p$, for the same homological reasons. Now, as in the
proof of Proposition \ref{P:lower dims}, $h_{p-1}:\pi_{p-1}(\td X_0)\to H_{p-1}(\td X_0)$ is
an epimorphism, and since $p<n-p$, we have the following commutative diagram as a piece of the
``homotopy-homology ladder'':

\begin{equation*}
\begin{CD}
\pi_p(\td C_0,\td X_0) @>\cong >> \pi_{p-1}(\td X_0)\\
@V h_p VV @V h_{p-1} V \text{onto} V\\
H_p(\td C_0,\td X_0) @>\cong >> H_{p-1}(\td X_0).
\end{CD}
\end{equation*}
The truth of the lemma is now apparent. 
\end{proof}

Since $\pi_p(\td C_0, \td X_0)\cong \pi_p(C_0,X_0)$, there is therefore a
map $(D^p,S^{p-1})\to (C_0,X_0)$ whose lift represents the element
$p(t^{-1})\beta$, where $\beta$ is some generator of $H_p(\td C_0, \td
X_0)$ as a $\Lambda$-module. Let $(D,S)$ represent the image of the
disk-sphere pair. Since $p<\frac{n}{2}$, we can choose $(D,S)$ to be an
embedded disk in $C_0$ whose boundary is an embedded sphere in $X_0$ and
such that $\text{int}(D)\subset\text{int}(C_0)$. Note that, chasing the
exact sequences around, this boundary must lift to the element
$p(t^{-1})\alpha\in H_{p-1}(\td X_0)$, for some generator, $\alpha$, of
$H_{p-1}(\td X_0)$ as a $\Lambda$-module. Let $R$ be an open tubular
(regular) neighborhood of $(D,S)$. We will show that $(D^p\times
S^{n-p}-R, L)$, denoted by $(\Delta, L)$, is our desired knotted disk
pair.

We begin by showing that $D^p\times S^{n-p}-R$ is the $n$-disk. As in the proof of Proposition
\ref{P:lower dims}, the fact that $p(1)=1$ implies that $S$ represents a generator of $\pi_{p-1}(S^{p-1}\times S^{n-p})\cong \Z$, and hence, using the long exact homotopy sequence, $(D,S)$ represents a generator of $\pi_p(D^p\times S^{n-p},S^{p-1}\times S^{n-p})\cong H_p(D^p\times S^{n-p},S^{p-1}\times S^{n-p})\cong\Z$. Hence $(D,S)$ is homotopic in $(D^p\times S^{n-p},S^{p-1}\times S^{n-p})$ to $D^p\times x_0$ for some $x_0\in S^{n-p}$. If $2(p+1)\leq n$, then this can be taken as an ambient isotopy, and then clearly $ D^p\times S^{n-p}-R\cong D^p\times D^{n-p}\cong D^n$.

If $p=\frac{n}{2}-\frac{1}{2}$, then $2(p-1)+2=n-1$, so there is still an
ambient isotopy of the boundary, $S^{p-1}\times S^{n-p}$, which takes $S$
to $S^{p-1}\times x_0$ for some $x_0\in S^{n-p}$, and this isotopy can be
extended to an ambient isotopy of all of $D^p\times S^{n-p}$, \cite{RS}.
When we form $\Delta=D^p\times S^{n-p}-R$, the new boundary will therefore
be \[(S^{p-1}\times S^{n-p}- S^{p-1}\times x_0\times
D^{n-p})\cup_{S^{p-1}\times x_0\times S^{n-p-1}} (D\times S^{n-p-1})\cong
S^{n-1},\] since this is the standard torus decomposition of $S^{n-1}$. We
will next show that $\Delta$ is contractible. Then, since the manifold
$\Delta$ will be a homotopy $n$-disk bounded by an $n-1$ sphere, $\Delta\cup
\bar c(\bd \Delta)$, (where $\bar c(\bd \Delta)$ indicates the closed cone
on the
boundary), will be a homotopy $n$-sphere. But $n\geq 5$, so $\Delta\cup
\bar c(\bd \Delta)$ is in fact a true sphere by the Poincare conjecture
and
$\Delta$ will be a true $n$-disk. It remains to show that $\Delta$ is
contractible. $\Delta$ is simply-connected because, for $p$ as given,
$\pi_1(D^p\times S^{n-p})=0$ and dim$(D)=p<n-2$. Together, these imply
by a general
position argument that $\pi_1(D^p\times S^{n-p}-D)\cong \pi_1(D^p\times
S^{n-p}-R)=0$, as well. To compute the homology of $\Delta$, we observe
that $H_i(\Delta)=0$, $i\geq n$, since $\Delta$ is an $n$-manifold with
boundary; $H_n(\Delta,\bd \Delta)\cong H_{n-1}(\bd \Delta)\cong \Z$, since
$\bd \Delta\cong S^{n-1}$ and these are generated by the orientation
classes; and $H_i(\Delta)\cong H_i(\Delta, \bd \Delta)$, $0<i<n$, by the
long exact sequence of the pair. By excision and homotopy equivalence,
\begin{align*}
H_i(\Delta, \bd \Delta)&\cong H_i(D^p\times S^{n-p}, (S^{p-1}\times S^{n-p})\cup
\bar{R})\\
&\cong H_i(D^p\times S^{n-p}, (S^{p-1}\times S^{n-p})\cup D).
\end{align*}
By the Mayer-Vietoris sequence, and using the fact that $S$ is a generator of $\pi_{p-1}(S^{p-1}\times S^{n-p})$ and hence of $H_{p-1}(S^{p-1}\times S^{n-p})$, 
\begin{equation*}
\td H_i((S^{p-1}\times S^{n-p})\cup D)\cong
\begin{cases}
\Z, &i=n-p, n-1\\
0,  &i\neq n-p,n-1.
\end{cases} 
\end{equation*}
But the generators of $H_i((S^{p-1}\times S^{n-p})\cup D)$ in dimensions
$n-p$ and $n$, respectively, are the generators of $H_i(S^{p-1}\times
S^{n-p})$ in the same dimensions (note that $D$ has no simplices of
dimension $>p$). The former is also a generator of $H_{n-p}(D^p\times
S^{n-p})$, and the latter is the boundary of the orientation class of
$H_n(D^p\times S^{n-p}, S^{p-1}\times S^{n-p})$. Therefore, using these
isomorphisms, the long exact sequence yields that $H_{n-p}(D^p\times
S^{n-p}, (S^{p-1}\times S^{n-p})\cup D)=0$, $0<i<n$, which, by our
calculations, shows that $\td H_i(\Delta)=0$, $i>0$. Therefore, by the
Whitehead Theorem, $\Delta$ is contractible, and we have finished proving
that $\Delta\cong D^n$.

For the last step in the proof of the proposition, we begin by fixing some
notation. Let $C$ denote $\Delta-L$. We can lift $(D,S)\subset (C_0,X_0)$
to an infinite number of copies $(D_i,S_i)$, $-\infty<i<\infty$,
corresponding to the translates of a lift of $(D,S)$ under the covering
translations, and similarly we lift the neighborhood $R$ to an infinite
number of $R_i$. Let $\td D$, $\td S$, and $\td R$ denote the disjoint
unions $\amalg_i D_i$, $\amalg_i S_i$, and $\amalg_i R_i$, respectively.
Then $\td C_0-\td R$ covers $C_0-R\cong C$. Furthermore, let $X$ denote
the manifold $(X_0-S\times D^{n-p})\cup_{S\times S^{n-p-1}} (D\times
S^{n-p-1})$ which results as the new complement of $\bd L$ in $\bd
\Delta$. Then the cover $\td X$ corresponds to $(\td X_0-\td R)\cup_{\td
S_\times S^{n-p-1}} (\td D\times S^{n-p-1})$. We show that the homology of
$(\td C, \td X)$ is as desired.

By excision, $H_i(\td C,\td X)\cong H_i(\td C_0,\td X_0\cup \td R)$. Let us denote $\td X_0\cup \td R$ by $\td Y$. Since 
\begin{align*}
\td H_i(\td X_0)\cong 
\begin{cases}
\Lambda, & i=p-1,n-p\\
0, & i\neq p-1,n-p
\end{cases}
\end{align*}
and since the lifts $S_i$ represent $t^i\lambda(t^{-1})\alpha\in H_p(\td X_0)$, an easy Mayer-Vietoris sequence argument gives 
\begin{align*}
\td H_i(\td Y)\cong 
\begin{cases}
\Lambda/p(t^{-1}), & i=p-1\\
\Lambda, & i=n-p\\
0,	& i\neq p-1,n-p.
\end{cases}
\end{align*}
Then using 
\begin{align*}
\td H_i(\td C_0)&\cong 
\begin{cases}
\Lambda, &i=n-p\\
0, &i\neq n-p,
\end{cases}
\end{align*}
the long exact sequence of the pair $(\td C_0,\td Y)$ gives $H_p(\td
C_0,\td Y)= \Lambda/p(t^{-1})$. The only other part of the sequence that
bears checking is where we have 

{\footnotesize
\begin{diagram}
H_{n-p+1}(\td C_0)&\rTo& H_{n-p+1}(\td C_0, \td Y) &\rTo& H_{n-p}(\td Y) &\rTo^i& H_{n-p}(\td C_0) &\rTo & H_{n-p}(\td C_0, \td Y) &\rTo &  H_{n-p-1}(\td Y)\\
\dTo && \dTo && \dTo && \dTo && \dTo && \dTo\\
0 &\rTo & ? &\rTo & \Lambda &\rTo & \Lambda &\rTo & ? &\rTo & 0.
\end{diagram}}

\noindent But the isomorphism $H_{n-p}(\td Y)\cong \Lambda$ comes from the isomorphism
$H_{n-p}(\td
X_0)\overset{\cong}{\longrightarrow} H_{n-p}(\td Y)$ induced by inclusion in the Mayer-Vietoris sequence we used, and we already know that $H_{n-p}(\td X_0)\overset{\cong}{\longrightarrow} H_{n-p}(\td C_0)$ is induced by inclusion (see above). Therefore, the map $i$ is an isomorphism and 
\begin{equation*}
H_i(\td C, \td X)\cong 
\begin{cases}
\Lambda/(p(t^{-1})), &i=p\\
0, & 0<i<n-1, i\neq p
\end{cases}
\end{equation*}
as claimed.
\end{proof}

Putting the results of these propositions together yields the following
classification of polynomials which can be realized as the
Alexander polynomials, $\lambda_i$, of a disk knot.
 
\begin{theorem}\label{T:suff}
Given polynomials $p_i(t)\in \Lambda$, $0<i<n-1$, $n\geq 4$, such that
$p_i(1)=\pm 1$ for each $i$ and $p_{n-2}(t)|p_1(t^{-1})$,
there
exists a knotted embedding $S^{n-2} \hookrightarrow S^n$ with at most
isolated point singularities such that the Alexander polynomials, $\lambda_i(t)$, of the 
knot are the given polynomials. 
\end{theorem}
\begin{proof}
The case $n=4$ has already been show. Suppose $n\geq5$. The necessity of the
conditions on the $\lambda_i(t)$ has been shown above in Section \ref{S: nec cond}.
For the sufficiency, let
$p_1(t)=p_{n-1}(t^{-1})r(t)$. By \cite{L66},
there is a locally-flat knot $S^{n-2}\subset S^n$ whose first
and $(n-2)$nd Alexander polynomials are $p_{n-1}(t^{-1})$ and
$p_{n-1}(t)$, respectively, and whose other Alexander polynomials
are all $1$. By Propositions \ref{P:lower dims} and
\ref{P:upper dims}, above, we can form $n-3$ separate knotted disk pairs
such that the first pair has first Alexander polynomial $r(t)$,
the $i$th pair has $i$th Alexander polynomial $\lambda_i(t)$, $1<i<n-2$,
and all the rest of the Alexander polynomials are $1$.
Then, taking the cone on the boundary of each knotted disk pair gives a
knotted sphere pair with point singularity,  $S^{n-2}\subset S^n$, and 
with the same Alexander
polynomials. Taking the knot sum of all of these knots (with the
connections being made in the neighborhoods of locally-flat points of the
embeddings) gives the desired knot because Alexander polynomials multiply
under knot sum.
\end{proof}

As a corollary and sample application, we can re-prove the following known
result concerning the Alexander polynomials of locally-flat slice
sphere knots.

\begin{corollary}\label{C: slice}
For any $n\geq3$ and collection of polynomials $p_i(t)\in\Lambda$,
$0<i\leq\lfloor
\frac{n-1}{2}\rfloor$, such that $p_i(1)=\pm 1$ and, if $n$ is odd,
$p_{ \frac{n-1}{2}}(t)\sim r(t)r(t^{-1})$ for some $r(t)$, there
is a locally flat slice knot $S^{n-2}\subset S^n$ whose \emph{i}th Alexander
polynomials, $\lambda_i(t)$, $0<i\leq \lfloor \frac{n-1}{2}\rfloor$, are
the
$p_i(t)$. These conditions on $p_{\frac{n-1}{2}}(t)$ are also
necessary. (Note that this also determines the Alexander polynomials for
$\lfloor \frac{n-1}{2}\rfloor<i<n-1$, as well, since
$\lambda_i(t)\sim \lambda_{n-i-1}(t^{-1})$ for locally flat knots.)
\end{corollary}
\begin{proof}
The necessity that $p_i(1)=\pm 1$ is proven in \cite{L66}.

We construct $\lfloor \frac{n-1}{2}\rfloor$ distinct locally flat slice
knots such that the \emph{i}th Alexander polynomial of the \emph{i}th knot
is
$p_i(t)$ and the rest of the Alexander polynomials (for $i\leq \lfloor
\frac{n-1}{2}\rfloor$) are 1.  Then our desired knot is the knot sum of
these, since Alexander polynomials multiply under knot sum and the knot
sum
of slice knots is slice.

Consider the long exact sequence of the pair $(\td C, \td X)$ for the
complement of a knotted disk pair $D^{n-1}\subset D^{n+1}$ . By Theorem
\ref{T:suff}, there is such a knotted disk pair whose Alexander polynomial
corresponding to $H_i(\td C)$ is $p_i(t)$ and such that $H_p(\td C)=0$
for all other $p$, $0<p<n$. This implies by Corollary \ref{C: nec con}
 that
$H_{n-i}(\td C, \td X)$ has Alexander polynomial $p_i(t^{-1})$ and all
other $H_p(\td C, \td X)=0$. For $i< \frac{n-1}{2}$, we obtain immediately
from the long exact sequence of the pair that the boundary knot with
complement $X$ has the desired homology. In fact, since $n-i>i+1$ in this
case, the exact sequence implies that $H_i(\td X)\cong H_i(\td C)$ and
$H_p(\td X)=0$, $0<p\leq\lfloor \frac{n-1}{2}\rfloor$, $p\neq i$. So the
boundary knot is the desired slice knot with
$\lambda_i(t)\sim p_i(t)$, $\lambda_{n-i-1}(t)\sim p(t^{-1})$, and no
other
non-trivial Alexander polynomials.

For $i=\frac{n-1}{2}$, the necessity that
$\lambda_{\frac{n-1}{2}}(t)\sim r(t)r(t^{-1})$ follows just as in Corollary \ref{C: slice
knots} for the
case of the classical slice knots, where $n=3$,
by using the product formula
\eqref{E: alt prod} which relates the Alexander polynomials corresponding to the homology
modules of
$\td X$, $\td C$, and $(\td C, \td X)$. Note that for slice knots, this condition along with
$p_{\frac{n-1}{2}}(1)=\pm 1$ implies Levine's necessary condition (d)
for Alexander polynomials \cite{L66}. 

For the construction, we consider the
knotted disk pair whose $\frac{n-1}{2}$-th Alexander polynomial
corresponding
to $H_{\frac{n-1}{2}}(\td C)$ is $r(t)$ and whose other Alexander
polynomials, corresponding to the other $H_p(\td C)$, are all
trivial. Such a
disk pair exists by the theorem and the fact that we must have $r(1)=\pm
1$ in order to have $p_{\frac{n-1}{2}}(1)=\pm1$. Then the Alexander
polynomials corresponding to the $H_p(\td C, \td X)$ are $r(t^{-1})$ for
$p=\frac{n+1}{2}$ and trivial otherwise. It then follows from the long exact
sequence of the pair $(\td C, \td X)$ that $H_p(\td X)=0$ for
$0<p<n-1$, $p\neq \frac{n-1}{2}$. Lastly, from the product formula of
Corollary \ref{C:division}, which relates
the three sets of Alexander polynomials, it must be that the
Alexander polynomial $\lambda_{\frac{n-1}{2}}(t)\sim r(t)r(t^{-1})$.

For the cases $i<\frac{n-1}{2}$, we can also note the existence of the given slice knots by observing that our procedure for
creating disk pairs with given Alexander polynomials in this range restricts on the boundary to Levine's procedure
\cite{L66} for creating knotted sphere pairs with the same prescribed Alexander polynomials.
\end{proof}

\subsubsection{Realization of all Alexander polynomials}\label{S: real
disk knot}

So far, we have stated all of our realizability conditions for disk
knots in terms of the Alexander polynomials $\lambda_i$ which correspond
to the $\Gamma$ modules $H_i(\td C;\Q)$. We now turn to a characterization
which simultaneously involves all of the Alexander polynomials we have
discussed: $\lambda_i$, $\mu_i$, and $\nu_i$, which correspond,
respectively, to $H_i(\td C;\Q)$, $H_i(\td C,\td X;\Q)$, and $H_i(\td
X;\Q)$. It will prove more natural, however, to consider the corresponding
\emph{subpolynomials} (see Section \ref{S: poly alg}. In fact, the long
exact reduced homology sequence of the pair $(\td C,\td X)$ yields an
exact
polynomial sequence
\begin{equation*}
1\to \lambda_{n-2}\to\mu_{n-2}\to \nu_{n-3}\to\cdots\to
\nu_1\to\lambda_1\to \mu_1\to 1 
\end{equation*}
with all polynomials in primitive form.
By the discussion in Section \ref{S: poly alg}, this gives rise to a
sequence of primitive polynomials of the form
\begin{equation*}
1\to c_{n-2}\to c_{n-2}a_{n-3}\to a_{n-3}b_{n-3}\to\cdots\to
a_1 b_1 \to b_1c_1\to c_1 \to 1.
\end{equation*}
As noted there, knowledge of the $a_i$, $b_i$, and $c_i$ is equivalent to
knowledge of the $\lambda_i$, $\mu_i$, and $\nu_i$. While we have been
referring to $\lambda_i$, $\mu_i$, and $\nu_i$ as the Alexander polynomials
of the disk, we will refer to $a_i$, $b_i$, and $c_i$ as the \emph{Alexander
subpolynomials}.

With this notation, we can observe the following lemma:
\begin{lemma}\label{L: subduality}
For a locally-flat disk knot $D^{n-2}\subset D^n$, $c_i(t)\sim
c_{n-i-1}(t^{-1})$ for $0<i<n-1$, $a_j(t)\sim b_{n-j-2}(t^{-1})$ for
$0<j<n-2$, and each of these
polynomials evaluates to $\pm 1$ at $t=1$.
\end{lemma} 
\begin{proof}
The last statement follows from the fact that each of the $a_i$, $b_i$,
and $c_i$ is a primitive polynomial in $\Lambda$ which divides another
primitive polynomial which evaluates to $\pm  1$ at $t=1$.

The other results follow by induction from the outside of the sequence to
the inside using Theorem \ref{T: disk knot nec con} and Corollary
\ref{C: nec con} (which states that $\lambda_i(t)\sim
\mu_{n-i-1}(t^{-1})$) together with Levine's \cite{L66} necessary
conditions for the Alexander polynomials of the locally-flat boundary
sphere knot (which states that $\nu_{i}(t)\sim
\nu_{(n-1)-i-1}(t^{-1})$). First,
we note that $c_{n-2}(t)\sim\lambda_{n-2}(t)\sim
\mu_1(t^{-1})\sim c_1(t^{-1})$. Next, $\mu_{n-2}(t)\sim
\lambda_{1}(t^{-1})$, but $\mu_{n-2}(t)\sim c_{n-2}(t)a_{n-3}(t)$ and
$\lambda_{1}(t^{-1})\sim  b_1(t^{-1})c_1(t^{-1})$. It follows that
$a_{n-3}(t)\sim b_1(t^{-1})$. The lemma is established by continuing this
procedure to the middle of the exact polynomial sequence.
\end{proof}

In these terms, we can now completely classify the Alexander polynomials
of disk knots with the exception of some extra middle-dimensional
information which we will study in the next section. For now, we impose
one unnecessary condition for the purpose of collecting the results that
derive from our work in this section.

\begin{theorem}\label{T:all polys}
Let $\nu_i(t)$, $\lambda_i(t)$, and $\mu_i(t)$ denote the Alexander
polynomials of a knotted $D^{n-2}\subset D^n$ corresponding to
$H_i(\td X)$, $H_i(\td C)$, and $H_i(\td C, \td X)$, respectively, and suppose $n\geq
4$. Recall that we can assume these to be primitive in $\Lambda$. Let $q_i(t)$, $0<i\leq\lfloor
\frac{n-2}{2}\rfloor$;
$r_i(t)$, $0<i\leq\lfloor \frac{n-2}{2}\rfloor$; and $p_i(t)$, $0<i\leq\lfloor
\frac{n-1}{2}\rfloor$, be polynomials in
$\Lambda$ satisfying the following properties:
\begin{enumerate}
	\item \begin{enumerate}
			\item $q_i(1)=\pm 1$
			\item $r_i(1)=\pm 1$
			\item $p_i(1)=\pm1$

	and hence, in particular, they must each also be primitive in $\Lambda$.
		\end{enumerate}
	\item There exist polynomials $a_i(t)$, $b_i(t)$, and $c_i(t)$,
primitive in $\Lambda$, such that
		\begin{enumerate} 
			
		\item $q_i(t)\sim a_i(t)b_i(t)$, $0<i<\lfloor
\frac{n-2}{2}\rfloor$
				
		\item \label{cond 3} $q_{\frac{n-2}{2}}(t)\sim
b_{\frac{n-2}{2}}(t)b_{\frac{n-2}{2}}(t^{-1})$, if $n$ is even

		\item $r_i(t)\sim b_i(t)c_i(t)$, $0<i\leq\lfloor
\frac{n-2}{2}\rfloor$
			
		\item $p_i(t)\sim c_i(t)a_{i-1}(t)$,  $0<i\leq\lfloor \frac{n-1}{2}\rfloor$
					(taking  $a_0(t)=1$)
		\item \label{cond 1} $c_{\frac{n-1}{2}}(t)\sim c_{\frac{n-1}{2}}(t^{-1})$,  if $n$ is odd
		\item \label{cond 2}$c_{\frac{n-1}{2}}(-1)=\pm \text{an odd
						square}$ if $n=2k+1$, $k$ 
even, and $c_{\frac{n-1}{2}}(t)$ is in normal form (defined below).
			
		\end{enumerate} 
		\end{enumerate} 
Then, there exists a knotted $D^{n-2}\subset D^n$ such that $\nu_i(t)\sim
q_i(t)$, $\lambda_i(t)\sim r_i(t)$, and $\mu_i(t)\sim p_i(t)$ in the
relevant
ranges. Note that this determines all of the Alexander polynomials using
$\nu_i(t)\sim \nu_{n-i-2}(t^{-1})$ and $\lambda_i(t)\sim \mu_{n-i-1}(t^{-1})$.
Furthermore, these conditions are necessary except for condition \ref{cond 2}. \end{theorem}

\begin{proof}

Most of the necessity has already been shown either above or in \cite{L66}.
Condition \ref{cond 1} and \ref{cond 3} are necessary by Lemma
\ref{L: subduality} since $n-\frac{n-1}{2}-1=\frac{n-1}{2}$ is an integer
for $n$ odd, and, if $n$ is even, $\frac{n-2}{2}$ is an integer,
$n-\frac{n-2}{2}-2=\frac{n-2}{2}$, and $\nu_{\frac{n-2}{2}}\sim
a_{\frac{n-2}{2}}(t)b_{\frac{n-2}{2}}(t)$. 
Condition
\ref{cond 2} is not necessary and will be weakened in the following section. Following Levine
\cite{L66}, the normal form for $c_{\frac{n-1}{2}}(t)$ is the member,
$c(t)$,  
of its similarity class in $\Lambda$ which satisfies $c(t)=c(t^{-1})$ and
$c(1)=1$. Due to
the other conditions on  $c_{\frac{n-1}{2}}(t)$, it is always possible to
find such a similar polynomial with these properties (see \cite[\S
1.5]{L66}).

To show that there is such a knotted disk, we will find a series of knots
whose only nontrivial Alexander polynomials are the $a_i(t)$, $b_i(t)$, or
$c_i(t)$ in the correct dimesnions and then take a knot sum. In this case, by a
knot sum of two knots we mean the following: Suppose we have two knotted disk
pairs $D^{n-2}_1\subset D^n_1$ and $D^{n-2}_2\subset D^n_2$. We can first
take the connected sum of the $D^n_i$ to form a new disk, $D^n$, in which the
two knotted $(n-2)$-disks are embedded disjointly. We can then connect the
two
knotted disks by an unknotted tube, $T\cong D^1\times D^{n-2}$, which connects a
neighborhood of a point of $\bd D_1^{n-2}$ in $D_1^{n-2}$ with a neighborhood
of a point of $\bd D_2^{n-2}$ in $D_2^{n-2}$ with the reverse orientation.
The knot sum is the closure of \[(D_1^{n-2}\cup D_2^{n-2}\cup T)
-(\text{int}(T)\cup(T\cap (D_1^{n-2}\cup D_2^{n-2}\cup \bd D^n)).\] On the
boundary this is the usual knot sum of the knots given by $\bd D_1^{n-2}$ and
$\bd D_2^{n-2}$ in $\bd D^n$. From the usual Mayer-Vietoris considerations,
the Alexander polynomials multiply under this knot sum. (Alternatively, we
could define knot sum by coning on the boundary of the knotted disk pairs to
create possibly non-locally flat sphere pairs, taking the usual knot sum
with
the connections in neighborhoods of locally flat points, and following our
original procedure for turning such a singular knot back into a knotted
disk pair.)

For the remainder of the proof, the term ``lower dimensional Alexander
polynomial'' will refer to the Alexander polynomials in the lower
dimensions listed above. We make this definition  in order to avoid
repetition of conditions that
arise from
the duality in the upper dimensions.

The knotted disk pairs with $\nu_j(t)\sim \lambda_j(t)\sim b_j(t)$ for a single
$j$, $0<j<\lfloor \frac{n-2}{2}\rfloor$, and all other lower dimensional
Alexander polynomials trivial is constructed in the proof of  Corollary \ref{C: slice}, as is the case where
$\nu_{\frac{n-2}{2}}(t)\sim b_{\frac{n-2}{2}}(t)b_{\frac{n-2}{2}}(t^{-1})$,
$\lambda_{\frac{n-2}{2}}(t)\sim b_{\frac{n-2}{2}}(t)$ and all other lower
dimensional Alexander polynomials are trivial. 

To construct the knotted disk pair with $\lambda_j(t)\sim \mu_j(t)\sim
c_j(t)$ for a single $j$, $0<j\leq\lfloor \frac{n-1}{2}\rfloor$, and all
other lower dimensional Alexander polynomials trivial, we have stipulated
the sufficient conditions \cite{L66} to construct a locally flat knot
$S^{n-2}\subset S^n$ with $c_j(t)$ as its only non-trivial lower
dimensional Alexander polynomial in the usual sense. Then, we can take the
trivial slicing of this knot by excising a ball neighborhood of a point of
the knot. Then $\lambda_j(t)\sim c_j(t)$; all of the $\nu_i(t)\sim
1$; $\mu_j(t)\sim c_j(t)$, by
the long exact sequence; and all other lower
dimensional Alexander polynomials are trivial.

To construct the knotted disk pair with $\mu_{j+1}(t)\sim \nu_j(t)\sim a_j(t)$ for
a single $j$, $0<j< \lfloor \frac{n-2}{2}\rfloor$, and all other lower
dimensional Alexander polynomials trivial, we recall our procedure from
Proposition \ref{P:upper dims}. As in that proof, we can construct a knotted
disk pair with $\mu_{j+1}(t)\sim a_j(t)$, $0<j< \lfloor
\frac{n-2}{2}\rfloor$,
and all other $\mu_i(t)$ and all lower dimension $\lambda_i(t)$ trivial. But the
construction restricted to the boundary was a surgery which we can
check to be equivalent to Levine's method \cite{L66} for creating a knotted
$S^{n-1}\subset S^{n-3}$ whose only lower dimensional Alexander polynomial is
$\nu_j(t)\sim a_j(t)$. Therefore, this knotted disk pair is the desired one. 

By taking the knots sums of these constructions as indicated above, we
obtain our desired knot. 
 \end{proof}

%% file: alexd.fix.tex

\subsection{The middle dimension polynomial}\label{S: middim disc}

We now turn to the case of realizing Alexander polynomials in the middle
dimension of a $(2q+1)$-disk knot, $q$ even. In particular, for $q>2$, we give a
characterization of the polynomials $c(t)\in \Lambda$ such that there exists a  locally-flat knotted
disk pair $D^{2q-1}\subset D^{2q+1}$ such
that $c(t)$ is the Alexander polynomial factor shared by $H_q(\td C)$ and
$H_q(\td C, \td X)$. Equivalently, $c(t)$ is the Alexander polynomial
associated to the modules \emph{ker}$(\bd_*)$ and \emph{cok}$(i_*)$ in the
long exact sequence of the pair
\begin{equation*}
\begin{CD}
@>>> H_k(\td X) @>i_*>> H_k(\td C) @>>> H_k (\td C, \td X)@>\bd_*>>
H_{k-1}(\td X)@>>>. 
\end{CD}
\end{equation*} 
We will show, in particular, that for any realizable $c(t)$ there exists such a knot with all other
Alexander polynomials (and Alexander subpolynomials) equal to $1$ so
that $c(t)$ will be the only non-trivial Alexander polynomial of $H_*(\td C)$ and
$H_*(\td C, \td X)$. We can then use the usual procedure of  taking connected sums
of disk knots
to combine this with other Alexander polynomials. 

We will in fact show something more. We will realize entire
$\Lambda$-modules and intersection pairings. First, we need a few
definitions. Following Levine \cite{L77}, we say that a $\Lambda$ module,
$A$,
is of \emph{type K} if it is finitely generated and multiplication by $t-1$ induces an automorphism of
$A$. It is a standard fact, see e.g. \cite{L77}, that the Alexander
modules of locally-flat sphere knots must be of type K. The standard proof
following Milnor
\cite{M68} extends easily to
disk knots. We provide it here to add the few words
relevant for the cases we will consider. 

\begin{lemma}
Let $D^{n-2}\subset D^{n}$ be a locally flat disk knot. Then the
$\Lambda$-modules  $H_i(\td X)$, $H_i(\td C)$, and $H_i (\td C,
\td X)$, $i>0$, are all of type K.
\end{lemma}  
\begin{proof}
That the modules are finitely generate follows from the usual argument
stemming from the fact that there is a one-to-one correspondence between
generators of the chain complexes of the knot exteriors (which are finite
complexes) and  the generators of the chain complexes of the infinite cyclic covers as
$\Lambda$-modules. Specifically, we choose one lift of each simplex.

 Now, let $W$ stand for $C$, $X$, or the pair $(C,X)$. Then
we have an exact sequence 
\begin{equation*}
\begin{CD}
0 @>>> C_i(\td W) @>t-1>> C_i(\td W) @>>> C_i(W) @>>>0,
\end{CD}
\end{equation*}
which generates the long exact homology sequence
\begin{equation}\label{E: Milnor sequence}
\begin{CD}
@>>> H_i(\td W) @>t-1>> H_i(\td W) @>>> H_i(W) @>\bd_*>>.
\end{CD}
\end{equation}
But by Alexander duality (Alexander duality for a ball), $X$ and $C$ are  
homology circles, and it is easy to see that $(C,X)$ is a homology ball.
Therefore, it is immediate for $i\geq2$ that $t-1$ is an automorphism of
the homology groups of the covers of $C$ and $X$, and in all dimensions $i>0$ for
$W=(C,X)$. 

For the remaining cases, we note that the long exact sequence
must terminate as
\begin{equation*}
\begin{CD}
@>0>> H_1(W)  @>\cong>> H_0(\td W) @>0>> H_0(\td W) @>\cong >>
H_0(W) @>>>0.
\end{CD}
\end{equation*}
The rightmost isomorphism is induced by the projection of a point that
generates
$H_0(\td W)$ to a point that generates $H_0(W)$. To see the other
isomorphism, observe that in the diagram chase that defines the boundary
map of the long exact sequence of homology, the generator of $H_1(W)$, a
meridian of the boundary sphere knot (disk knot), gets lifted to a
1-chain in the cover whose boundary is $(t-1)x$ for some point $x$ in $\td
W$ representing an element of $C_0(\td W)$. This is the image of $x$
under the map $t-1$. Thus the image of a generator of $H_1(W)$ goes to a
generator of $H_0(\td W)$ under the boundary map of the long exact
sequence. The sequence now shows that $t-1$ is also a homology
automorphism of the cover for $i=1$. 
\end{proof}

\subsubsection{The Blanchfield pairing}
 
We will also need the \emph{Blanchfield pairing} on the infinite cyclic cover. We summarize its construction and properties following \cite{L77}. More details can be found in the references cited there. (Note: for notational convenience we introduce the symbol $\bl{\,}{}$ to represent the general Blanchfield pairing and reserve $\blm{\,}{}$ for the induced middle dimensional self-pairing.)

For $M$ a compact $m$-dimensional PL-manifold with boundary which admits a regular
cover with group of covering transformations $\pi$, one first defines an
intersection pairing on the chain groups of the covers, $C_q(\td M,\bd \td
M)\times
C_{m-q}(\td M^1)\to \Z[\pi]$, $\alpha\times \beta\to \alpha\cdot\beta$, where the
chain groups are thought of as (left) $\Z[\pi]$ modules and $M^1$ represents the
dual complex to the triangulation of $M$. This pairing is bilinear over $\Z$ and
satisfies
\begin{enumerate}
\item $(g\alpha)\cdot\beta=g(\alpha\cdot\beta)$, for $g\in\pi$, 
\item $\alpha\cdot\beta=(-1)^{q(m-q)}\overline{\beta\cdot\alpha}$, where the bar denotes the
antiautomorphism of $\Z[\pi]$ induced by $\bar g=g^{-1}$ for $g\in \pi$,
\item $(\bd \alpha)\cdot \beta =(-1)^q \alpha\cdot (\bd \beta)$.
\end{enumerate}
This induces a pairing on the appropriate homology groups. 

Now assume $\pi=\Z$ and that $\alpha\in H_q(\td M, \bd \td M)$ and $\beta\in
H_{m-q-1}(\td X)$ are $\Lambda$-torsion elements represented by chains $z\in
C_q(\td M, \bd \td M)$ and $w\in C_{m-q-1}(\td M^1)$. Then $z=\frac{1}{\lambda}\bd
c$ for some $\lambda\in \Lambda=\Z[\Z]$ and $c\in C_{q+1}(\td M, \bd \td M)$.
Define $\bl{\alpha}{\beta}= \frac{1}{\lambda} c\cdot w\, \text {mod }
\Lambda$. This
induces a well-defined pairing $H_q(\td M, \bd \td M) \times H_{m-q-1}(\td M)\to
Q(\Lambda)/\Lambda$, where $Q(\Lambda)$ is the field of fractions of $\Lambda$.
The pairing $\bl{\,}{}$ is conjugate linear, meaning that it is additive in each
variable and $\bl{\lambda
\alpha}{\beta}=\lambda\bl{\alpha}{\beta}=\bl{\alpha}{\bar\lambda \beta}$. The
conjugation on elements of $\Lambda$ is induced by $\bar t=\frac{1}{t}$.
Furthermore, if $m=2q+1$, $\bl{\,}{}$ induces a self-pairing $\blm{\,}{}:H_q(\td
M)\times H_q (\td M) \to Q(\Lambda)/\Lambda$ by $\blm{\alpha}{\beta}=
\bl{j_*\alpha}{\beta}$, where $j_*: H_q(\td M)\to H_q(\td M,\bd\td M)$ is the map
of the long exact sequence. This pairing is $(-1)^{q+1}$-Hermitian, meaning that
$\blm{\alpha}{\beta}=(-1)^{q+1} \overline{\blm{\beta}{\alpha}}$.

We observe that, in the case of a disk knot, the arguments of \cite[\S 5]{L77}
carry over to show that $\bl{\,}{}$ is a non-singular pairing on the $Z$-torsion
free parts of $H_q(\td C)$ and $H_{n+1-q}(\td C, \bd \td C)$.  If $n=2q+1$, the
induced pairing $\blm{\,}{}$ on $H_q(\td M)$ further induces a nondegenerate
(though possibly singular) conjugate linear $(-1)^{q+1}$-Hermitian pairing on the
$\Z$-torsion free part of \emph{coim}$(j_*)$ (though we will keep the same
notation $\blm{\,}{}$): To see that this is well-defined, we observe that if
$\alpha+\beta, \gamma\in H_q(\td M)$, $\beta\in$\emph{ker}$(j_*)$, then
$\blm{\alpha
+\beta}{\gamma}=\bl{j_*(\alpha+\beta)}{\gamma}=\bl{j_*\alpha}{\gamma}=
\blm{\alpha}{\gamma}$. By the Hermitian property, similar considerations hold for
the second argument so that $\blm{\,}{}$ only depends on $H_q(\td
M)/$\emph{ker}$(j_*)$. For the non-degeneracy, note that the non-singularity of
$\bl{\,}{}$ implies that for every non-zero, non-$\Z$ torsion element
$j_*\alpha\in H_q(\td M, \bd \td M)$, there is a non-$\Z$ torsion $\gamma\in
H_q(\td M)$ such that $\blm{\alpha}{\gamma}=\bl{j_*\alpha}{\gamma}\neq 0$. But
from the well-definedness argument above, if $\gamma\in$\emph{ker}$(j_*)$, then
$\blm{\alpha}{\gamma}=0$. So $\gamma$ has non-zero image when projected into
\emph{coim}$(j_*)$. This establishes the non-degeneracy since such a $\gamma$
exists for all such $j_* \alpha$.

In the above arguments, we can replace $H_i(\td C, \bd\td C)$ with $H_i(\td C, \td
X)$, $i<n-2$, as follows: First observe that $\bd \td C=\td X \cup_{S^{n-3}\times
\R} D^{n-2}\times \R$ so that $H_i(\bd \td C, \td X)=H_i(D^{n-2}, S^{n-3})$, by
excision and homotopy equivalence. Therefore, by the long exact sequence of the
pair, the map induced by inclusion, $j_*: H_i(\td X)\to H_i(\bd\td C)$, is an
isomorphism for $i<n-3$ and onto for $i=n-3$. Using long exact sequences and the
five-lemma, this implies that $H_i(\td C, \td X)\cong H_i(\td C, \bd \td C)$,
$i<n-2$.

Summarizing part of this discussion gives:

\begin{proposition}\label{P: Blanchfield pairing} 
Let $D^{n-2}\subset D^{n}$ be a disk knot, $n=2q+1$, $k>0$. Let $f(A)$ denotes the
$\Lambda$-module $A$ modulo its $\Z$-torsion.
Then $H_q(\td C)$
and
$H_q(\td C, \td X)$ are $\Lambda$-modules of type K, and the non-singular pairing
$\bl{\,}{}: f(H_q(\td C, \td X))\times f(H_q(\td C))\to Q(\Lambda)/\Lambda$
induces a nondegenerate conjugate linear $(-1)^{q+1}$-Hermitian pairing
$f($\emph{coim}$(j_*))\times f($\emph{coim}$(j_*))\to Q(\Lambda)/\Lambda$.
\end{proposition}

\subsubsection{Realization of middle dimensional pairings}
We will establish a converse to Proposition \ref{P: Blanchfield pairing}:

\begin{proposition}\label{P: pairing}
Let $A$ be a $\Z$-torsion free $\Lambda$ module of type K with a non-degenerate
conjugate linear $(-1)^{q+1}$-Hermitian pairing $\blm{\,}{}:A\times A \to
Q(\Lambda)/\Lambda$. Then there exists a disk knot $D^{n-2}\subset D^{n}$,
$n=2q+1$, $q>2$, such that:
\begin{enumerate}
\item $H_q(\td C)=A$,
\item $H_i(\td C)=0$, $0<i<n-1$, $i\neq q$,
\item \label{I:3}$H_i(\td X)=0$, $0<i<n-2$, $i\neq q-1$,
\item \label{I:4}$H_{q-1}(\td X)=0$ is a $\Z$-torsion module,
\item \label{I:5}$H_i(\td C, \td X)=0$, $0<i<n-1$, $i\neq q$,

\item the pairing on $H_q(\td C)$ is given by $\blm{\,}{}$. (Note that $H_q(\td
X)=0$ implies that $H_q(\td C)\cong A\cong$\emph{coim}$(j_*)$ in the long exact
sequence). 
\end{enumerate}
\end{proposition}
\begin{proof}

By \cite[Proposition 12.5]{L77}, given such an $A$ and $\blm{\,}{}$, there exists a smooth compact
$(2q+1)$-dimensional manifold, $C$, such that $\pi_1(C)=\Z$, $H_q(\td C)=A$, $H_i(\td C)=0$ for
$i\neq 0$,$q$, and the given pairing $\blm{\,}{}$ corresponds to the pairing on $H_q(\td C)$. The
proof consists of first being able to write the defining matrices for the presentation of $A$ and
the pairing with respect to the basis of presentation in certain forms, which follows from
\cite[Proposition 12.3]{L77} and the remarks before \cite[Proposition 12.5]{L77} because $A$ is of type K; and then (\cite[Lemma 12.2]{L77}) constructing $C$ using the matrix information to attach appropriate $q$-handles to 
\begin{equation*}
C_0=(\#_{i=1}^m S^q\times D^{q+1})\#(S_1\times D^{2q}),
\end{equation*}
where the presentation matrix has size $m\times m$ and $\#_{i=1}^m S^q\times D^{q+1} $ denotes the connected sum of $m$ copies of $ S^q\times D^{q+1} $. This $C$ will be our disk knot complement.

We observe that $C$ is a homology circle: $H_1(C)=\pi_1(C)=\Z$ as above, and the
triviality in the remaining dimensions, $i>1$, follows from Milnor's exact
sequence \eqref{E: Milnor sequence} and $A$ being of type K. As Levine notes in
Proposition 12.6 of \cite{L77}, we also have $\Z=\pi_1(C)\cong \pi_1(C-K)\cong
\pi_1(\bd C) $, where $K$ is the $(q+1)$-dimensional subcomplex formed from the
cores of the handles added onto $C_0$: $C-K$ deformation retracts to $\bd C$, and
the claim follows from general position since $q>2$. Thus, we can add a $2$-handle
onto $C$ along a generator of $\pi_1(\bd C)$ to obtain a manifold which is
contractible (using the Hurewicz and Whitehead theorems) with simply-connected
boundary, hence a disk by \cite{Sm62}. If $D^2\times D^n$ is the attached handle,
then our disk knot is $0\times D^n$, the ``cocore'' of the handle. Clearly then C
is the knot's exterior with modules and pairings as claimed.

It remains to show that properties (\ref{I:3}), (\ref{I:4}), and (\ref{I:5}) hold.
Again from the proof of \cite[Proposition 12.6]{L77}, $H_i(\bd \td C)=0$ for
$i\neq 0$, $q-1$, $q$, $2q-1$, and $H_i(\td C, \bd\td C)=0$, $i\neq q-1$, $q$,
$2q-1$. The argument uses the Hurewicz theorem, a version of Poincare duality for
coverings (\cite{M68} and \cite{L77}), and a universal coefficient short exact
sequence for torsion $\Lambda$-modules. As noted above, $j_*: H_i(\td X)\to
H_i(\bd\td C)$ is an isomorphism for $i<n-3$, so $H_i(\td X)=0$, $i<q-1$, and
therefore $H_i(\td X)=0$, $q<i<n-2$, by the duality of sphere knot modules
\cite{L77}. Similarly, $H_i(\td C, \td X)=0$, $i<q$ or $q<i<n-1$, using the
long exact sequence of the pair $(\td C, \td X)$.

At this point we have all of the Alexander modules $0$ except for $H_q(\td C)$,
$H_q(\td C, \td X)$, $H_{q-1}(\td X)$, and $H_{q}(\td X)$. But $H_{q}(\td X)$ must
be $0$ because the non-degeneracy of the pairing on $C$ implies that the map $j_*:
H_q(\td C)\to H_q(\td C, \td X)$ of the long exact sequence must be injective. It
now follows from Levine's duality properties for the Alexander modules of locally-flat sphere
knots (see \cite{L77}) that $H_{q-1}(\td X)$ is a $\Z$-torsion module.
\end{proof}

\subsubsection[Matrix representations of the middle dimension module and
its pairing]{Matrix representations of the middle dimension module and
its pairing; Characterization of the middle dimensional polynomial in
these terms}

It is also useful to study these middle-dimensional Alexander modules using
presentation matrices. We first examine the form that these matrices take. From
the proof of Corollary \ref{C: kern}, we know that $c_q(t)$ is the determinant of
the presentation matrix of the kernel of the map
\begin{equation*}
\bd_*: H_q(\td C, \td X; \Q)\to H_{q-1}(\td X;\Q)
\end{equation*}
in the long exact sequence of the pair. Let us denote this kernel module by $H$. Equivalently, it is the determinant of the presentation matrix of the isomorphic coimage of the map
\begin{equation*}
p_*: H_q(\td C; \Q)\to H_{q}(\td C, \td X;\Q).
\end{equation*}
We will refer to this module as $\bar H$. 

To obtain a presentation matrix for $H$ (or $\td H$), recall the Mayer-Vietoris
sequences \eqref{E:MV1} and \eqref{E:MV2} used to obtain the presentation matrices
for the Alexander polynomials.
The long exact sequences of the \emph{rational} homology of the pairs $(V, F)$ and
$(Y,Z)$ must split at each term as exact sequences of vector spaces; in other 
words, each is isomorphic to an exact sequence of vector spaces of the form
$\to A\oplus B\to B\oplus C\to C\oplus D \to$. This splitting and
exactness is preserved under the tensor product with the free module $\Gamma$ over
$\Q$. Hence we obtain the following diagram which commutes owing to the obvious
commutativity at the chain level induced by the maps in the Mayer-Vietoris
sequence
and by naturality of the homology functor. The $0$ terms arise by truncation,
using
our knowledge that the Mayer-Vietoris sequences break into short exact sequences.
\begin{equation}\label{E:grid}
\begin{diagram}
& & \dTo && \dTo && \dTo\\
0 & \rTo & H_q(F;\Q)\otimes_{\Q}\Gamma & \rTo & H_q(Z;\Q)\otimes_{\Q}\Gamma & \rTo &  H_q(\td X;\Q) & \rTo & 0 \\
& & \dTo && \dTo && \dTo \\
0 & \rTo &  H_q(V;\Q)\otimes_{\Q}\Gamma & \rTo^{d_1} &  H_q(Y;\Q)\otimes_{\Q}\Gamma & \rTo^{e_1} &  H_q(\td C;\Q) & \rTo & 0\\
 & & \dTo^{r} && \dTo^{s} && \dTo \\
0 & \rTo &  H_q(V, F;\Q)\otimes_{\Q}\Gamma & \rTo^{d_2}& H_q(Y,Z;\Q)\otimes_{\Q}\Gamma & \rTo^{e_2} &  H_q(\td C,\td X;\Q) & \rTo & 0\\
& & \dTo^{\bd_*} && \dTo^{\bd_*} && \dTo^{\bd_*} 
\end{diagram}
\end{equation}
Let $E$ and $G$ denote, respectively, the kernels of the boundary maps $\bd_*$ in
$H_q(V, F;\Q)\otimes_{\Q}\Gamma$ and  $H_q(Y,Z;\Q)\otimes_{\Q}\Gamma$. Let $J$,
$K$, and $L$ be the respective cokernels of the boundary maps of which $E$, $G$, 
and $H$ are the kernels. Then, by the snake lemma, we obtain an exact sequence
\begin{equation*}
\begin{CD}
0 @>>> E @>d>> G @>>> H @>>> J @>>> K @>>> L @>>>0.
\end{CD}
\end{equation*}
But note that by the splitting of the two leftmost (non-zero) vertical sequences
in the diagram \eqref{E:grid}, $J$ and $K$ are direct summands of
$H_{q-1}(F;\Q)\otimes \Gamma$ and $H_{q-1}(Z;\Q)\otimes \Gamma$, respectively.
Hence the injectivity of the map $H_{q-1}(F;\Q)\otimes \Gamma \to
H_{q-1}(Z;\Q)\otimes \Gamma$ in the Mayer-Vietoris sequence implies that the
induced map $J\to K$ must also be injective. Therefore, we get an exact sequence
\begin{equation*}
\begin{CD}
0 @>>> E @>d>> G @>>> H @>>>0.
\end{CD}
\end{equation*} 
This sequence gives a presentation for $H$. In
fact, $E$ and $G$ are certainly free $\Gamma$-modules (each being a rational
vector space tensored with $\Gamma$ over $\Q$), and the matrix representing
$d$ gives a presentation matrix for $H$. Note that the matrix for $d$ is a  submatrix (which we can arrange to be the upper left submatrix)
of the matrix representing $d_2$. The
generators of $E$ and $G$ are the elements $\{e_i\otimes 1\}$ and
$\{g_i\otimes
1\}$,
where $\{e_i\}$ and $\{g_i\}$ are the generators of the direct summands of
$H_k(V,
F;\Q)$ and $H_q(Y,Z;\Q)$ which are the images of $H_q(V;\Q)$ and $H_q(Y;\Q)$
under the projection maps of the exact sequences of the pairs. Furthermore,
$d$ must be represented by a square matrix: If it had more columns than
rows, then there would be more generators than relations in $H$ which is
impossible since we know that $H$ is a $\Gamma$-torsion module; and if it
had
more rows than columns, then since the elements in the summand $E$ map only
into the summand $G$ and $d_2$ is square, $d_2$ would be forced to have
determinant $0$, which is also impossible as we saw in the proof of Theorem
\ref{T: disk knot nec con}. Hence the matrix of $d$ gives a
square presentation of $H$, which we can take to be the upper left $m\times
m$ submatrix of $d_2$, by changing bases if necessary.
Similar considerations give the isomorphic presentation of the coimages
\begin{equation*}
\begin{CD}
0 @>>> \bar E @>\bar d>> \bar G @>>> \bar H @>>>0.
\end{CD}
\end{equation*}

From the termwise splitting of the leftmost column of \eqref{E:grid} before
tensoring with $\Gamma$, there exist vectors space summands $\td E$ and
$\Td{\Bar E}$ in $H_q(V,F;\Q)$ and $H_q(V;\Q)$, respectively, such that
$E=\td E\otimes\Gamma$ and $\bar E=\Td{\Bar E}\otimes \Gamma$.
Furthermore, $r$ can be written as $\td r\otimes \text{id}$, where $\td r:
H_q(V;\Q)\to H_q(V,F;\Q)$ is the map of the long exact sequence induced by
inclusion (and induces the isomorphism of the summands $\td E\cong
\Td{\Bar E}$). We can make similar conclusions about $G$ in the second
column of \eqref{E:grid} and carry over all of the bar and tilde
notations. Identifying quotient vector spaces with summands, for
convenience, we obtain the diagram:

\begin{equation*}
\begin{CD}
\Td{\Bar E} \subset H_q(V; \Q) @>>> H_q( Y;\Q)\supset \Td{\Bar G}\\
@V\td r VV @V\td s VV\\
\td E \subset H_q(V,F; \Q)@>>> H_q(Y,Z; \Q) \supset \td G.
\end{CD}
\end{equation*}

We will now choose suitable bases for $\td E$, $\Td{\Bar E}$, $\td G$, and $\Td {\Bar
G}$. Consider
now the \emph{integral} homology groups and long exact
sequence maps $\td r_{\Z}: H_q(V)\to H_q(V,F)$ and $\td s_{\Z}: H_q(Y)\to
H_q(Y,Z)$. As abelian groups, each of these is the direct sum of its free part and
its torsion part, and we can choose bases so that maps between the free summands
are represented by diagonal matrices ordered so that all of the zero diagonal
entries are moved to the bottom right \cite[\S 11]{MK}. Clearly then when we
tensor with $\Q$, we get the maps in the above diagram with the vector space
summands $\td E$, $\Td {\Bar E}$, $\td G$, and $\Td {\Bar G}$ being represented by
the $\Q$ spans of the first $m$ basis elements of the groups, i.e. we can now
choose bases $\{\alpha_i\}$, $\{\beta'_i\}$, $\{\gamma_i\}$, $\{\delta'_i\}$, of the
free parts of $H_q(V,F)$, $H_q(Y)$, $H_q(Y, Z)$, $H_q(V)$, such that, upon
tensoring with $\Q$, the first $m$ elements of each basis will span $\td E$,
$\Td{\Bar G}$, $\td G$, and $\Td {\Bar E}$, respectively, and the maps $\td
r_{\Z}\otimes \Q$ and $\td s_{\Z}\otimes \Q$ induce the appropriate vector space
isomorphisms.  Furthermore, $\{\alpha_i\otimes 1\}_{i=1}^m$, $\{\beta'_i\otimes
1\}_{i=1}^m $, $\{\gamma_i\otimes 1\}_{i=1}^m $, $\{\delta'_i\otimes 1\}_{i=1}^m $
now span $E$, $\bar G$, $G$, and $\bar E$.

We claim also that with these choices $\td E$ and $\Td{\Bar G}$ are dual
with respect to the linking pairing $L'$ (see Section \ref{S: nec
cond})
and $\Td{\Bar
E}$ and $\td G$ are dual with respect
to $L''$, which will allow us to perform changes of bases of $\Td{\Bar G}$ (to
$\{\beta_i\}$) and
$\Td{\Bar E}$ (to $\{\delta_i\}$)  such that 
\begin{equation*} 
L'(\alpha_i \otimes \beta_j)=L''(\gamma_i \otimes \delta_j)=\delta_{ij},
\end{equation*}
$1\leq i \leq m$. The changes of bases can be taken to be integrally unimodular
(see below). 

We proceed
by  first proving that the duals to the $\{\gamma_i\}_{j=1}^m$ under $L''$ span
$\Td{\Bar{E}}\subset H_q(V;\Q)$. To see this, we first observe that, up to sign,
$L''([v],\td s([y]))=L'(\td r([v]),[y])$ for $[v]\in H_q(V;\Q)$ and $[y]\in
H_q(Y;\Q)$. This follows by considering the definition of the linking pairings. If
$v$ and $y$ are chains representing $[v]$ and $[y]$, then they also represent $\td
r[v]$ and $\td s[y]$ (as relative chains modulo the chain complexes $C_q(F)$ and
$C_q(V)$). Then $L''([v], \td s([y]))$ is the intersection number of $y$ with a
chain in $D^n$ whose boundary is $v$, while $L'(\td r([v]),[y])$ is the
intersection number of $v$ with a chain in $D^n$ whose boundary is $y$. By the
properties of intersection numbers, these agree up to sign. Now suppose that $v$
is an element of $H_q(V; \Q)$ which lies in the summand \emph{ker}$(\td r)$ and
that $\{\td s^{-1}\gamma_i\}_{i=1}^m$ are elements of $H_q(Y)$ which map onto the
$\gamma_i$. Then $0=L'(\td r(v),\td s^{-1}\gamma_i)=L''(v,\gamma_i)$. Therefore,
\emph{ker}$(\td r)$ is orthogonal to $\td G$ under $L''$. Thus, the dual subspace
to $\td G$, spanned by $\{\delta_i\}_{i=1}^m$, must lie outside \emph{ker}$(\td
r)$ and project onto an $m$-dimensional subspace of \emph{coim}$(\td r)=\Td{\Bar
E}$. But dim$(\Td {\Bar E})=$dim$(\td E)=m$ by isomorphism and dim$(\td
E)=$dim$(\td G)$ because the map $d$ was a square presentation. This proves that
$\Td{\Bar E}$ and $\td G$ are dual. 

It also follows from the discussion of the last paragraph that we must have
$\delta_i\in \text{ker}(\td r)$ for $i>m$: Suppose not. Without loss of generality,
suppose $\delta_{m+1}\notin \text{ker}(\td r)$. Then, in the rational vector space
$\Td{\Bar{E}}\otimes \Q$, there will be (at least ) $m+1$ linearly independent vectors,
$\{\delta_i\}_{i=1}^m$, which do not lie in the kernel. But since the kernel  has
dimension $n-m$ (rationally), the span of $\{\delta_i\}_{i=1}^m$ must intersect the
kernel. Therefore
there is a vector $v\in  \text{ker}(\td r)\otimes \Q$ such that $v=\sum_{i=1}^{m+1}n_i
\delta_i$, $n_i\in \Q$. Furthermore, there must be some $n_j$, $j\leq m$, such that
$n_j\neq 0$ (else $v=n_{m+1}\delta_{m+1}\notin \text{ker}(\td r)\otimes \Q$). Then
$L''(v,\alpha_j)=n_j\neq 0$, contrary to the results of the last paragraph. Therefore,
$\delta_i\in \text{ker}(\td r)$ for $i>m$. Now, since each $\delta'_i$ is an integral
linear combination of the $\{\delta_i\}$ (since each is a basis for
$H_q(V)$), the same must be true under the projection to $\Td{\Bar{E}}$, i.e. the
projection of each $\delta'_i$ is an integral linear combination of the projections of
the $\{\delta_i\}$. But since $\delta_i\in \text{ker}(\td r)$ for $i>m$, each projected  
$\delta'_i$ is a linear combination of the projections of
$\{\delta_i\}_{i=1}^m$. Since the projected $\{\delta'_i\}_{i=1}^m$ form a basis for
$\Td {\Bar E}$, it
follows that the projections $\{\delta_i\}_{i=1}^m$ also form a basis for $\Td {\Bar
E}$. In particular, we see that $\Td {\Bar
E}$ is integrally dual to $\td G$ (and hence also rationally when tensored with
$\Q$).
In what follows, we shall
also refer to the projections of the $\{\delta_i\}_{i=1}^m$ into $\Td {\Bar
E}$ as $\{\delta_i\}_{i=1}^m$.

Similar considerations apply for the other case to show that $\Td {\Bar 
G}$ with basis $\{\beta_i\}_{i=1}^m$ is dual to $\td E$.

Next, we can apply our previous notations, procedures, and results (see Section \ref
{S: nec
cond}) to
these modules to obtain the formulae: 
\begin{align*}
i_{+*}(\delta_j)&=\sum_i \lambda_{ij} \beta_i\\
i_{-*}(\delta_j)&=\sum_i \sigma_{ij} \beta_i\\
i_{+*}(\alpha_j)&=\sum_i \mu_{ij} \gamma_i\\
i_{-*}(\alpha_j)&=\sum_i \tau_{ij} \gamma_i\\
L'(\alpha_k\otimes i_{+*}(\delta_j))&=\sum_i \lambda_{ij} L'(\alpha_k\otimes \beta_i)=\lambda_{kj}\\
L'(\alpha_k\otimes i_{-*}(\delta_j))&=\sum_i \sigma_{ij} L'(\alpha_k\otimes \beta_i)=\sigma_{kj}\\
L''(i_{+*}(\alpha_j)\otimes \delta_k)& =\sum_i \mu_{ij} L''(\gamma_i\otimes \delta_k)=\mu_{kj}\\
L''(i_{-*}(\alpha_j)\otimes\delta_k )&=\sum_i \tau_{ij} L''(\gamma_i\otimes \delta_k)=\tau_{kj},
\end{align*}
where all of the indices run only to $m$ and everything is of dimension $q$. We get presentation matrices 
\begin{align*}
P_1(t)&=(t\sigma_{ij}-\lambda_{ij})\\
P_2(t)&=(t\tau_{ij}-\mu_{ij})
\end{align*}
for $\bar H$ and $H$, and we know that $\sigma_{jk}=\mu_{kj}$ and $\lambda_{jk}=\tau_{kj}$. 

We are further furnished with one more relation between the matrices $\mu$ and
$\tau$. Let $R=(R_{ij})$ be the matrix representation of $\td r|\Td{\Bar E}$.  
Let $v_i$ be a chain representing $\delta_i\in \Td{\Bar E}$, $1\leq i \leq m$, and
observe that the same chain (modulo chains in F) represents $\td r(\delta_i)\in
\td E$. Thus, using chains interchangeably with their appropriate homology
classes,
\begin{multline*}
L''(i_{+*}(\td r\delta_j)\otimes \delta_i)= L''(i_{+}(v_j)\otimes v_i) \\
=(-1)^{q+1}L''(i_{-}(v_i)\otimes v_j)= (-1)^{q+1}L''(i_{-*}(\td r\delta_i)\otimes \delta_j),
\end{multline*}
where the middle equality comes from the usual geometry of the isotopies obtained
by ``pushing along the bicollar'' (see Section \ref{S: linking numbers}), and 
the sign change is the usual sign change in the commutativity formula for a linking pairing induced by an intersection pairing (see \cite[Appendix]{GBF}). But
\begin{align*}
L''(i_{+*}(\td r\delta_j)\otimes \delta_i)
&= L''(i_{+*}(\sum_{k=1}^m R_{kj}\alpha_k)\otimes \delta_i)\\
&=\sum_{k=1}^m R_{kj} L''(i_{+*}(\alpha_k)\otimes \delta_i)\\
&=\sum_{k=1}^m R_{kj} \mu_{ik}.
\end{align*}
Similarly, we get that
\begin{align*}
L''(i_{-}(\td r\delta_i)\otimes \delta_j)
&=\sum_{k=1}^m R_{ki} \tau_{jk}.
\end{align*}
This yields the matrix equations 
\begin{equation*}
\mu\cdot R= (-1)^{q+1} (\tau \cdot R)' = (-1)^{q+1} R'\cdot \tau',
\end{equation*}
and we can conclude the following:

\begin{proposition}
The $\Gamma$-module $H$ has a presentation matrix of the form $\tau t - (-1)^{q+1}
R'\tau'R^{-1}$, where $R$ is
the matrix of the map $\Td{\Bar E}\to \td E$ induced by $\td r: H_q(V)\to H_q(V, F)$. $\bar H$ has presentation matrix $(-1)^{q+1}(R^{-1})'\tau R t-\tau'$. 
\end{proposition}

\begin{remark}
Both of these presentation matrices have the same determinant, up to sign, as expected.
\end{remark}

In this situation, we can say the following about the matrix of the pairing $\blm{\,}{}: \bar H \times \bar H \to Q(\Lambda)/\Gamma=Q(\Gamma)/\Gamma$:

\begin{proposition}\label{P: pairingmatrix}
In the above situation, taking $\{B_i\}_{i=1}^m$ as the generators of $\bar H$,
where $B_i$ is the image of the $\{\beta_i\otimes 1\}\in \Td{\Bar
E}\otimes\Gamma=\td E$ in $\bar H$, a matrix representative of the pairing
$\blm{\,}{}: \bar H \times \bar H \to Q(\Lambda)/\Gamma$ is
given by $\frac{t-1}{(R^{-1})'\tau -(-1)^{q+1}t\tau'R^{-1}}$.
\end{proposition}

\begin{proof}
The proof follows closely that of \cite[Proposition 14.3]{L77}.  We choose
particular lifts of $V$ and $Y$ which adjoin (i.e. any path from $t^{-1}Y$ to $Y$
must cross $V$, identifying $t$ as the covering translation) and identify
$\delta_i\in V$ with $\delta_i\otimes 1$, which we will call $\td \delta_i$ for
convenience. Set $p(t)\td \delta_i=\delta_i\otimes p(t)$ for $p(t)\in \Gamma$.
We treat the other bases similarly.

Since $ i_{+*}(\td \delta_i)=\sum_i \lambda_{ji} \td\beta_j$
and $i_{-*}(\td\delta_i)=\sum_j \sigma_{ji} \td\beta_j$ are induced by homotopies, there are chains $c_i$ and $c'_i$ such that 
\begin{align*}
\bd c_i&=\td\delta_i-\sum_j \lambda_{ji} \td\beta_j\\
\bd t c'_i&=\td\delta_i-\sum_j \sigma_{ji} t\td\beta_j.
\end{align*}
Thus 
\begin{align*}
\bd (tc'_i- c_i)&= -\sum_i t\sigma_{ji} \td\beta_j + \sum_i \lambda_{ji} \td\beta_j\\
&= \sum_i ( \lambda_{ji} - t\sigma_{ji})\td \beta_j.
\end{align*}

As usual, let $\lambda$ and $\sigma$ denote the matrices $(\lambda_{ji})$ and $(\sigma_{ji})$. Let $\Delta(t)=$det$( \lambda - t\sigma)$ and $M(t)=\Delta(t) ( \lambda' - t\sigma')^{-1}$, i.e. the matrix of cofactors of $( \lambda' -t \sigma')$.
Thus 
\begin{equation}\label{E: cancel}
\delta_{jk} \Delta(t) =\sum_i M_{ki}(t) ( \lambda_{ji} - t\sigma_{ji}),
\end{equation}
so that 
\begin{align*}
\Delta(t)\td \beta_k
&= \sum_j \delta_{jk} \Delta(t) \td\beta_j\\
&= \sum_{i,j} M_{ki}(t) ( \lambda_{ji} - t\sigma_{ji})\td \beta_j\\
&= \sum_i M_{ki}(t)\bd (tc'_i- c_i)\\
&= \bd(\sum_i M_{ki}(t) (tc'_i- c_i)).\\
\end{align*}

Now, as outlined above, to compute $\blm{B_k}{B_l}$, we choose representative chains for the $B_i$ (denoting both the chains and classes by the same symbol for simplicity) and find a chain $c$ such that $\bd c=p(t) B_k$ for some $p(t)\in \Lambda$. Then $\blm{B_k}{B_l}= \frac{1}{p(t)}c\cdot B_l$ mod $\Lambda$. 

Based upon the above computations, we can take $p(t)=\Delta(t)$ and $c(t)=\sum_i
M_{ki}(t) (tc'_i-c_i)$ from which
\begin{equation*}
\blm{B_k}{B_l}= \frac{\sum_i M_{ki}(t) (t(c'_i\cdot B_l) - (c_i\cdot B_l))   }{\Delta(t)}.
\end{equation*}
Since the $c_i$, $c'_i$, and $B_i$ all lie in the same lift of $Y$, the
intersection numbers in this formula are the ordinary intersection numbers in $Y\subset S^n$ and are thus the same as the usual linking numbers of the chains $\bd c_i$ and $\bd c'_i$ with $B_l$. Since a chain representing $\td \delta_i$ represents $\sum_j R_{ji}\td\alpha_j$ in $\td E$, we get
\begin{align*}
c_i\cdot B_l&=L''(\sum_j R_{ji}\td\alpha_j, \td\beta_l)- \sum_j \lambda_{ji} \ell(\td\beta_j,\td\beta_l)\\
&= R_{li}- \sum_j \lambda_{ji} \ell(\td\beta_j,\td\beta_l)\\
c'_i\cdot B_l&=L''(\sum_j R_{ji}\td\alpha_j, \td\beta_l)- \sum_j \sigma_{ji} \ell(\td\beta_j,\td\beta_l)\\
&= R_{li}- \sum_j \sigma_{ji} \ell(\td\beta_j,\td\beta_l),
\end{align*}
where $\ell(\td\beta_j,\td\beta_l)$ is the linking number in $S^n$ of chains representing $\beta_j$ and $\beta_l$. Thus

\begin{align*}
\blm{B_k}{B_l}&=
\sum_i \frac{ M_{ki}(t)}{\Delta(t)}[(t-1)R_{li} + \sum_j (\lambda_{ji}-t\sigma_{ji})\ell(\td\beta_j,\td\beta_l)]\\
&= \sum_i\frac{ M_{ki}(t) (t-1)R_{li}}{\Delta(t)} +\sum_{ij}\frac{ M_{ki}(t) (\lambda_{ji}-t\sigma_{ji} )\ell(\td\beta_j,\td\beta_l)}{\Delta(t)}\\
&= \sum_i\frac{ M_{ki}(t) (t-1)R_{li}}{\Delta(t)} +\sum_j \delta_{jk}\ell(\td\beta_j,\td\beta_l),
\end{align*}
where we have used equation \eqref{E: cancel} to simplify in the last step.
Since $\ell(\td\beta_j,\td\beta_l)$ is an integer, 
\begin{equation*}
\blm{B_k}{B_l}=\sum_i\frac{ M_{ki}(t) (t-1)R_{li}}{\Delta(t)} \text{  mod }\Lambda.
\end{equation*}
Thus the matrix of the pairing is given by 
\begin{align*}
\frac{t-1}{\Delta(t)}M(t)R' &=(t-1)( \lambda' - t\sigma')^{-1}R'\\
&= (t-1)(\tau-(-1)^{q+1}tR'\tau'R^{-1})^{-1}R'\\
&= \frac{t-1}{(R^{-1})'\tau -(-1)^{q+1}t\tau'R^{-1}}
\end{align*}
\end{proof}

Conversely, suppose we are given integer matrices $\tau$ and $R$ such that $R$
has non-zero determinant, $(R^{-1})'\tau R$ is an integer matrix, and
det$[M(1)]=\pm 1$, where $M$ is the matrix $M(t)= (-1)^{q+1}(R^{-1})'\tau Rt
-\tau'$. Let $A$ be the $\Lambda$-module whose presentation matrix is $M(t)$,
i.e. $A=\Lambda/M\Lambda$. Then $N(t)=\frac{t-1}{(R^{-1})'\tau
-(-1)^{q+1}\tau'tR^{-1}}= \frac{1-t}{(R^{-1})'M(t)'}$ determines a nondegenerate
$(-1)^{q+1}$-Hermitian form $\blm{\,}{}:A\times A\to Q(\Lambda)/\Lambda$ by
$\blm{a_1}{a_2}=a_1' N(t)\bar a_2$. (For a more general discussion of the
construction of which this is a minor modification, see \cite[\S 1]{T73}.) A
simple calculation shows that $N(t)$ is $(-1)^{q+1}$-Hermitian. The pairing is
well-defined because if $a_1=0$ in $A$, then $a_1\in M(t)\Lambda$ so that it can
be represented as $M(t)a_0$. Then $\blm{a_1}{a_2}=a_1' N(t) a_2= (M(t)a_0)'
N(t) \bar a_2= (1-t)a_0'M(t)' (M(t)')^{-1}R' \bar a_2= (1-t) a_0' R' \bar a_2 \in
\Lambda$. For the non-degeneracy, the work of Blanchfield \cite[pp. 350-1]{B57}
implies that $N_0(t)=[(R^{-1})'M(t)']^{-1}=[M(t)']^{-1}R'$ is a non-singular
$\Gamma$-module pairing $B\times B\to Q(\Gamma)/\Gamma$, where
$B=\Gamma/[N_0(t)']^{-1}\Gamma$, provided this is a $\Gamma$-torsion module. But
since $R$ is rationally unimodular,
$\Gamma/[N_0(t)']^{-1}\Gamma=\Gamma/M(t)R^{-1}\Gamma=\Gamma/M(t)\Gamma=A\otimes
\Q$. Hence, $B$ is $\Gamma$-torsion because $A$ is $\Lambda$-torsion. Thus
$N_0(t)$ can have no rows or columns composed completely of elements of $\Gamma$,
hence of $\Lambda$. This together with the fact that $(t-1)$ is an isomorphism on
$A$ (which is clearly of type $K$) shows that the pairing $N(t)$ is
non-degenerate.

Given any module and pairing as defined in the last paragraph, it is
realizable as the middle-dimensional module and pairing of a disk knot
$D^{2q-1}\subset D^{2q+1}$, $q>2$, by Proposition \ref{P: pairing}. Thus, we have proven:

\begin{theorem}\label{T:middim}
A polynomial $c(t)\in\Lambda$ can be realized as the Alexander subpolynomial
factor shared by $H_q(\td C)$ and $H_q(\td C, \td X)$ for the
locally-flat knotted disk pair $D^{2q-1}\subset D^{2q+1}$, $q>2$, if and only if 
$c(t)=$det$[M(t)]$, where $M(t)= (-1)^{q+1}(R^{-1})'\tau R t-\tau'$
for integer matrices $\tau$ and $R$, such that $R$ has non-zero determinant,
$(R^{-1})'\tau R$ is an integer matrix, and \emph{det}$[M(1)]=\pm 1$.
\end{theorem}

\begin{remark}
If the boundary knot is trivial, then we will have $R=I$, and we expect our formulae to look like those in \cite{L77} for the middle-dimensional duality of a sphere knot. That these formulae do not agree identically is due to two differences in conventions: The first is that we have chosen to use Levine's original convention of \cite{L66} for which map to label $i_{-}$ and which to label $i_+$ (these choices are reversed in \cite{L77}). The second is that while we have employed presentation matrices acting on the left, so that the matrix $A$ corresponds to the module $\Lambda^k/A\Lambda^k$, in \cite{L77} Levine allows his presentation matrices to act on the right so that $A$ corresponds to $\Lambda^k/\Lambda^kA$. Thus our presentation matrices are  transposed compared to those in \cite{L77}. 
\end{remark}

\subsubsection{Characterization of the middle dimension polynomial in
terms of pairings}

An alternative way of formulating  Theorem \ref{T:middim} is  the
following:
\begin{theorem}\label{T:middim2}
A primitive polynomial $c(t)\in\Lambda$ can be realized as the Alexander
polynomial factor shared by $H_q(\td C)$ and $H_q(\td C, \td X)$ for the locally-flat knotted disk
pair $D^{2q-1}\subset D^{2q+1}$, $q>2$, if and only if
$c(1)=\pm1$ and there exist an integer $\rho$ and a non-negative integer $\omega$
such that
$\frac{(t-1)^{\omega}\rho}{\pm c(t)}$ is the discriminant of a
$(-1)^{q+1}$-Hermitian form on a $\Lambda$-module of type $K$.
\end{theorem}
\begin{proof}

If $c(t)$ is the Alexander subpolynomial in primitive form, we know that $c(1)=\pm
1$, $c(t^{-1})\sim c(t)=\emph{det}\left[(-1)^{q+1}(R^{-1})'\tau R t-\tau'\right]$, and
$\frac{t-1}{(R^{-1})'\tau -(-1)^{q+1}\tau'tR^{-1}}$ is the matrix of a form of
the given type on a $\Lambda$-module of type $K$.  Letting $|\tau|$ stand for 
the number of rows (or columns) of the square matrix $\tau$, the discriminant
of the form is 
\begin{align} \label{E: discr}
\text{det}\left[\frac{t-1}{(R^{-1})'\tau -(-1)^{q+1}\tau'tR^{-1}}\right]
&=\frac{(t-1)^{|\tau|} \text{det}(R)}{ \pm c(t)}. 
\end{align}
Setting $\rho=\text{det}(R)$ and $\omega=|\tau|$ proves the claim in
this direction.

Conversely, suppose that we are given a primitive polynomial $p(t)\in \Lambda$
such that $p(1)=\pm 1$ and there exist $\rho$ and $\omega$ such that
$\frac{(t-1)^{\omega}\rho}{p(t)}$ is the discriminant, $D$, of a
$(-1)^{q+1}$-Hermitian form on a $\Lambda$-module of type $K$. Then by
Propositions \ref{P: pairing} and \ref{P: pairingmatrix}, the module and pairing
can be realized as an appropriate middle-dimensional knot pairing such that the
module has a presentation matrix of the form $(-1)^{q+1}(R^{-1})'\tau R t-\tau'$
and the pairing has a matrix of the form $\frac{t-1}{(R^{-1})'\tau
-(-1)^{q+1}\tau'tR^{-1}}$. The associated Alexander polynomial is then
$c(t)=$\emph{det}$((-1)^{q+1}(R^{-1})'\tau R t-\tau')$, while the discriminant is
$D=$
\emph{det}$[\frac{t-1}{(R^{-1})'\tau
-(-1)^{q+1}\tau'tR^{-1}}]$. Thus we have
\begin{align*}
c(t)&\sim \frac{(t-1)^{|\tau|}\text{det}(R)}{D}&&\text{as in the
last paragraph}&\\
   &\sim\frac{p(t) (t-1)^{|\tau|} \text{det}(R) } {( t-1)^{\omega}\rho
} &&\text{by assumption}\\ 
&\sim p(t) (t-1)^{|\tau|-\omega}\left(\frac{\text{det}(R)}
{\rho}\right). \end{align*}
But since we know that both $c(1)$ and $p(1)$ are equal to $\pm 1$, we must have
$\omega=|\tau|$ and $\rho=\text{det}(R)$, so that $c(t)\sim p(t)$ and $p(t)$ is
an Alexander polynomial of the desired type. \end{proof}

For the case where $q$ is odd, we already know from Sections \ref{S: nec
cond} and \ref{S: Realization}
that these polynomials must
be completely classified as those such that $c(1)=\pm 1$ and $c(t)\sim
c(t^{-1})$. I do not know of such a similarly straightforward 
classification for the case where
$q$ is even, although we will show that the previously imposed condition
that $|c(-1)|$ be a square is not necessary. In fact, we will show that
any quadratic polynomial, $c(t)\in \Lambda$, satisfying
\begin{enumerate}\item $c(1)=\pm1$ \item $c(t)\sim c(t^{-1})$
\end{enumerate} can be realized. It is easy to show that any such
polynomial has the form $at^2+(\pm 1-2a)t+a$. Now, we can just take
\begin{equation*}
R=\left(
\begin{matrix}
\pm 1+4a & 0\\
0 & 1
\end{matrix}
\right)
\hskip5mm
\tau=\left(
\begin{matrix}
a & 0\\
1 & 1
\end{matrix}
\right).
\end{equation*}
Then \begin{align*}
c(t)&=\text{det}[(R^{-1})'\tau R t-(-1)^{q+1}\tau']=
\text{det}\left( 
\begin{matrix}
at+a  &  1\\
(\pm 1+4a)t+1  & t+1
\end{matrix}
\right)\\
&=at^2 + (\pm1 -2a)t  + a. 
\end{align*}

Note that $c(-1)=4a\pm 1$, so that, by choosing $a$ suitably, we can realize any odd number as $c(-1)$. Observe that $c(-1)$ must be odd for any $c(t)$ satisfying conditions (1) and (2) above (see $\cite{L66}$).

For $q>2$, we can now replace condition \ref{cond 2} of Theorem \ref{T:all polys} with the
necessity statement of Theorem \ref{T:middim} or Theorem \ref{T:middim2}. The constructibility
follows by taking an appropriate connected sum with the knots constructed
in Proposition \ref{P: pairing}.

For $n=2q+1$, $q=2$, the methods employed above break down. The
difficulties in this case are clearly related to the difficulties of
classifying the $Z$-torsion part of the dimension-one Alexander module of a
locally-flat knot $S^2\subset S^4$ (see \cite{L77}).

\subsection{Conclusion}\label{S: conclusion alexa}

We summarize our results on the Alexander polynomials of locally-flat disk
knots, or equivalently, sphere knots with point singularities.

\begin{theorem}\label{T: alexa}
For $n\neq 5$ and $0<i<n-1$,   
$0<j<n-2$,
the following conditions are necessary and sufficient for $\lambda_i$,
$\mu_i$, and
$\nu_j$, to be the polynomials associated to the $\Gamma$-modules
$H_i(\td
C;\Q)$, $H_i(\td C,\td X;\Q)$, and $H_j(\td X;\Q)$ of a locally flat
disk knot $D^{n-2}\subset D^n$ or a
knot $S^{n-2}\subset S^n$ with point singularities  (see
Section \ref{S: nec cond disk} for the definitions of $C$ and $X$):  
There exist polynomials $a_i(t)$, $b_i(t)$, and $c_i(t)$, primitive in
$\Lambda$, such that
\begin{enumerate}
\item\begin{enumerate}
\item $\nu_i\sim a_i b_i$
\item $\lambda_i\sim b_ic_i$
\item $\mu_i\sim c_i a_{i-1}$
\end{enumerate}
\item \label{I: dual}\begin{enumerate}
\item $c_i(t)\sim c_{n-i-1}(t^{-1})$
\item $a_i(t)\sim b_{n-i-2}(t^{-1})$
\end{enumerate}
\item \label{I: disc} $a_i(1)=\pm1, b_i(1)=\pm1 , c_i(1)=\pm 1$, $a_0(t)=1$. 
\item \label{I: mddisk} If $n=2q+1$ and $q$ is even, then there exist an integer $\rho$ and a non-negative
integer $\omega$ such that $\frac{(t-1)^{\omega}\rho}{\pm c_q(t)}$ is the discriminant of a
$(-1)^{q+1}$-Hermitian form on a $\Lambda$-module of type $K$ (or equivalently, $c_q(t)=$det$[M(t)]$, where
$M(t)=
(-1)^{q+1}(R^{-1})'\tau R t-\tau'$ for integer matrices $\tau$ and $R$ such that $R$ has non-zero determinant
and $(R^{-1})'\tau R$ is an integer matrix). See Section \ref{S: middim disc} for definitions and more details.
\end{enumerate}

For a locally-flat disk knot $D^{3}\subset D^5$ or a knot $S^{3}\subset S^5$ with
point singularities, these conditions are all necessary. Furthermore, we can
construct any knot which satisfies both these
conditions and the added, perhaps unnecessary, condition that $|c_2(-1)|$ be
an odd square.
\end{theorem}
\begin{proof}
This is simply a conglomeration of the results of this section. Note that the
duality statements of \eqref{I: dual} follow from the duality results of Section
\ref{S: nec cond disk} and some simple polynomial algebra (see Lemma
\ref{L: subduality} in Section \ref{S: real    
disk knot}).
\end{proof}

\begin{remark}
For a locally-flat $D^1\subset D^3$, the boundary modules are all trivial in
dimensions greater than $0$. In fact the only nontrivial Alexander modules will be
$H_1(\td C;\Q)\cong H_1(\td C,\td X;\Q)$, and the only non-trivial polynomial
$c_1\sim
\lambda_1\sim \mu_1$ is completely classified by $c_1(t)\sim c_1(t^{-1})$ and
$c_1(1)=\pm 1$. Noting that the complement of a locally-flat 1-disk knot is
the same 
as that of the $S^1$ knot obtained by coning on the boundary (such a cone
remains locally-flat at all points), this follows by Levine's conditions
\cite{L66}. These conditions are equivalent
to the conditions stated above, taking $n=3$, although we have not proved here that
any such knot can be constructed. (The necessity could follow from our proof
for higher dimensional knots as the assumption $n>3$ was imposed only to focus our 
attention on knots which could have point singularities.) 
 \end{remark}

%% file: sstrata.fix.tex

\section{Knots with more general singularities}\label{S: sing}

\subsection{Introduction}
We now study the Alexander polynomials of non-locally-flat knots with
singularities more general than the point singularities of the last
section. To
be specific,
let $\alpha:S^{n-2}\hookrightarrow S^n$, $n\geq 3$, be a PL-embedding which is
locally-flat except on a singular set $\Sigma_{n-k}\subset S^{n-2}\subset
S^n$. Note that if we view $S^n$ as a PL-stratified space with singular locus
$S^{n-2}$, then $\Sigma$ will be a subpolyhedron of dimension less than $n-3$  
(see
Section 5 of \cite{GBF} or \cite{GBF2} for a more detailed discussion of
knots as stratified spaces). 

By analogy with Section \ref{S: disk knots}, we can study the homology
modules of the
infinite cyclic covers of the complement of the knot in the exterior of a regular neighborhood
of the singularity  and of the complement of the
intersection of the knot with the boundary of this regular neighborhood (see Section \ref{geom 
prel}). We can
also study the relative homology of the pair. These will all be
torsion $\Gamma$-modules, and thus we again obtain
three sets of polynomials to study: $\nu_i$, $\lambda_i$, and $\mu_i$. 

In Section \ref{S: duality*}, we show that $\lambda_i$ and $\mu_i$ satisfy a modified
version of the duality and normalization results for disk knots (Theorem \ref{T: disk
knot nec con}). In fact, the results are the same except for the appearance of a power of
$t-1$ as a factor in each $\mu_i$. In Section \ref{S: boundary knot*}, we show that the
$\nu_i$ also satisfy  self-duality and normalization conditions which generalize
Levine's conditions for the boundary locally-flat sphere knots of Section 
\ref{S: disk
knots}. Again there are extra $t-1$ factors to account for.

In Section \ref{S: subpolys*}, we discuss the factorization of the Alexander polynomials
into subpolynomials and rephrase our results in that context. In Section \ref{S: high
dim *}, we show that the Alexander polynomials are all trivial for sufficiently large
dimension index $i$.

\subsection{Necessary conditions on the Alexander
invariants}\label{S: sstrata}

\subsubsection{Geometric preliminaries}\label{geom prel}

Letting $K=\alpha(S^{n-2})$, $S^n-K$ is a homology
circle by Alexander duality, and again, just as in Section \ref{S: disk
knots}, we can
study the rational homology of its infinite cyclic cover viewed as a
module over $\Gamma=\Q[\Z]=\Q[t, t^{-1}]$. The Alexander polynomials are
the determinants of the presentation matrices of these modules. We begin
with some geometric preliminaries and notations.

As usual, we use the homotopy equivalent knot complement or knot exterior
as it suits our needs (the knot exterior being the complement in $S^n$ of
an open regular neighborhood of $K$). As already stipulated, $K$ is
locally flat away from $\Sigma$ so that each point in $K-\Sigma$ has a
distinguished neighborhood homeomorphic to $D^{n-2}\times c(S^1)$ (where
$c(X)$ denotes the open cone on $X$). Let $D$ denote the
manifold $S^n-N(\Sigma)$, where
$N(X)$ is the open regular neighborhood of $X$, and let $\bd D=S=\bd
\overline{N(\Sigma)}$.
Note that the boundary of the knot exterior
is the union of two pieces: a circle bundle in $D$ over $S^{n-2}-N(\Sigma)$  and the exterior in $S$ of an open neighborhood of $S\cap K$ in $S$. These pieces are joined along their intersection, a circle bundle in $S$ over $S\cap K$ which is the boundary of the closed regular neighborhood of $S\cap K$ in $S$.

We now construct a version of the Seifert surface in this context.

\begin{lemma}
There is a retract $R:S^n-K\to S^1$, where $S^1$ is a given PL-meridian of the
knot $K$.
\end{lemma}
\begin{proof}
Let $i:S^1\to S^n-K$ be the inclusion of the meridian $S^1$. Since $S^n-K$ is a
homology circle and the meridians generate its first homology group, $i_*$ is a
homology isomorphism in all dimensions. Hence, $\td H_*(S^n-K, S^1)\cong
\td H^*(S^n-K,S^1)=0$. But by Eilenberg-MacLane theory, since $S^1$ is a $K(\Z,1)$
and
since $(S^n-K, S^1)$ as a simplicial pair can also be considered a CW pair, this
implies
that the identity map $S^1\to S^1$ can be extended to a map $R: S^n-K\to S^1$ (see
\cite[8.1.12]{Sp}).
\end{proof}

\begin{proposition}\label{Seifert*}
With $D$ and $S$ as above, there is a bicollared  $(n-1)$-manifold $V\subset D$,
such that $\bd V=(K\cap D) \cup F$, where $F$ is a bicollared $(n-2)$-manifold in $S$ with $\bd
F=K\cap S$.
\end{proposition}  
\begin{proof}
Consider the regular neighborhood $\overline{N(K\cap D)}$ in $D$ which is
a 2-disk bundle over $K\cap D$. Let $\gamma$ be a fiber of the boundary
circle bundle $\bd \overline{N(K\cap D)}$ in $D$. Then $\gamma$, with the
proper choice of orientation, generates $H_1(S^n-N(K))=\Z$. In fact,
$\gamma$ certainly has linking number $\pm 1$ with $K$, corresponding to
the intersection point of $K$ with the obvious $2$-disk in $S^n$ which
bounds $\gamma$ and makes up a fiber of the regular neighborhood of $K$.  
Since the linking pairing $H_1(S^{n}-K)\otimes H_{n-2}(K)\cong \Z\otimes\Z
\to \Z$ is perfect (see \cite[Appendix]{GBF}) any element of $H_1(S^n-K)$ which
maps to $1$ under the pairing must be a generator.

Now, by the lemma, there is a retract $R:S^n-K$ to $\gamma$ which, by
restricting to $\bd \overline{N(K\cap D)}$ in $D$, provides a homotopy
trivialization of this circle bundle and hence of the disk bundle $N(K\cap
D)$ in $D$ by extending in the obvious way to the interior of the bundle.  
This homotopy trivialization is homotopic to an actual trivialization,
i.e. a projection $\bd \overline{N(K\cap D)}\cong (K\cap D)\times S^1 \to
S^1$, and by the homotopy extension principle, we can obtain a map $r$,
homotopic to $R$, such that $r|\bd \overline{N(K\cap D)}$ is the
projection to $S^1$. Consider now $r$ restricted to $\bd( S^n-N(K))$. We
wish to obtain our Seifert surface, $V$, by taking the transverse inverse
image of a generic point under a PL-approximation of $r$, but first we
must take care to avoid getting excess boundary components.

We can first take a PL-approximation to $r|\bd(S^n-N(K))$ which remains
the projection on $(K\cap D)\times S^1$. Now, we take the transverse
inverse image in $\bd S^n-N(K)$ of a sufficiently generic point, say $y$,
of $S^1$, which gives us a bi-collared $(n-2)$-submanifold. One component
of this submanifold consists of the union of $(K\cap D)\times y\subset
(K\cap D)\times S^1$ and a manifold $F\subset S$, with the union taken
along their common boundary $(K\cap S) \times y\subset (K \cap S) \times
S^1$. This can be seen by considering $F$ to be a component of the
transverse inverse image of the restriction to $S-S\cap N(K)$ of the
PL-approximation to $r|\bd( S^n-N(K))$. Unfortunately, there may be excess
closed components of the inverse image in $S$, but these can be removed by
replacing the approximation to $r'$ with the map to $S^1$ determined by
the connected bicollared submanifold consisting of the main component
discussed above (in particular, the map which takes the submanifold to the
point $y\in S^1$, the hemispheres of the bicollar to the two halves of the
circle, and the rest of $\bd (S^n-N(K))$ to the point antipodal to $y$).
Since $S^n-N(K)\sim_{h.e.} S^n-K$ is a homology circle, $H^2(S^n-N(K), \bd
(S^n-N(K)))\cong H_{n-2}(S^n-N(K))=0$, (recall $n\geq 4$), and therefore
there is no obstruction to extending this new map to a map $r'':
S^n-N(K)\to S^1$. Now we take the transverse inverse image of a 
PL-approximation to $r$ at $y$ or another sufficiently close point and
discard excess components to obtain a bicollared submanifold in $S^n-N(K)$
which will have the desired properties once we extend it trivially to the
interior of the disk bundle $\overline{N(K\cap D)}$ in $D$.

\end{proof}

We now establish some notation. We have already denoted $V\cap S$ by $F$. Let $Y=D-V$, $Z=Y\cap S$,
$W=V\cup \overline{N(\Sigma)}$, and $\Omega=Y \cup \overline{N(\Sigma)}$. Note that both $W$ and
$\Omega$ contain $\Sigma$. We observe that $D-(K\cap D)$ is homotopy equivalent to $S^n-K$, and so
we
can consider the homology of either to study the Alexander polynomials.

We begin our study of the Alexander invariants with the following observations
and definitions: Let $C$ be the knot complement $ D-(K\cap D)\sim_{h.e.} S^n-K$, and let $\tilde{C}$
be the infinite cyclic cover associated with the
kernel of the abelianization $\pi_1(C)\to\Z$. Letting $t$ denote a generator of
the covering translation, the homology groups of $\tilde{C}$ are finitely
generated $\Lambda$-modules ($\Lambda=\Z[\Z]=\Z[t, t^{-1}]$) since $C$ has a
finite polyhedron as a deformation retract, and the rational
homology groups $ H_*(\tilde{C};\Q) \cong H_*(\tilde{C})\otimes_{\Z} \Q$ are
finitely generated modules over the principal ideal domain
$\Gamma=\Q[t,t^{-1}]\cong \Lambda\otimes_{\Z}\Q$. Therefore, letting
$M_R(\lambda)$ denote the R-module of rank 1 with generator of order $\lambda$,
$H_q(\tilde{C};\Q)\cong \bigoplus_{i=1}^k M_{\Gamma}(\lambda_{q_i})$, where we
can choose the $\lambda_{q_i}$ so that: 1) The $\lambda_{q_i}$ are primitive in
$\Lambda$ but are unique up to associate class in $\Gamma$ and 2)
$\lambda_{q_{i+1}}|\lambda_{q_i}$. For $0<q<n-1$, these are called the
\emph{Alexander invariants} of the knot complements. The polynomial
$\lambda_q=\prod_{i=1}^k \lambda_{q_i}$, which is also primitive in $\Lambda$,
is the Alexander polynomial of the knot complement.  We will also consider the
relative homology modules $H_*(\tilde{C}, \tilde{X};\Q)$, where $X$ is the
complement in $S$ of $K\cap S$ (the ``link complement'','' as $S$ is the link of
$\Sigma$) and $\tilde{X}$ is its infinite cyclic covering. It will be clear
from our construction that $\tilde{X}$ and the cover of $X$ in $\tilde{C}$ are
equivalent. Then $H_q(\tilde{C}, \tilde{X};\Q)$ has the same properties listed
above for $H_q(\tilde{C};\Q)$ and its own Alexander invariants $\{\mu_{q_i}\}$,
$0<q<n-1$, and relative Alexander polynomial $\mu_q=\prod_i \mu_{q_i}$.

\subsubsection{Duality and normalization theorem}\label{S: duality*}

We will prove the following theorem analogous to that already established in
Section \ref{S: disk knots} for
the case of a point singularity:

\begin{theorem}\label{T:duality*}
Let $p+q=n-1$ with knot and notation as above. The following properties hold:
\begin{enumerate}
\item $\lambda_p(1)=\pm 1$,
\item $\mu_q(t)\sim \lambda_p(t^{-1})(t-1)^{\td B_{q-1}}$, where $\sim$ denotes
associativity of elements in $\Lambda$ (i.e. $a\sim b$ if and only if $a=\pm
t^k b$ for some
$k$) and $\td B_i$ is the $i$th reduced Betti number of $\Sigma$ (i.e. the Betti number of the reduced homology).
\end{enumerate} 
\end{theorem}

The proof will occupy the next several pages. We begin the proof by finding
$\Gamma$-module presentations for
$H_*(\tilde{C};\Q)$ and $H_*(\tilde{C}, \tilde{X};\Q)$ by studying the
Mayer-Vietoris sequences for the infinite cyclic cover obtained by cutting and
pasting along the Seifert surface $V\subset D$.

We construct $\tilde{C}$ as in Section \ref{S: disk knots} by
first cutting $D$ open along $V$ to
create a manifold, $Y'$, which is homotopy equivalent to $Y$ and whose boundary
is $Z$ plus
two copies of $V$, $V_+$ and $V_-$, identified along $K\cap D$, and by then
pasting
together a countably infinite number of disjoint copies $(Y^i, V_+^i,V_-^i)$,
$-\infty<i<\infty$, of $( Y'-K, V_+-K, V_{-}-K)$ by identifying $V^i_+ -K$ with
$V^{i+1}_- -K$ for all $i$. Then $\tilde{X}$ is the sub-manifold resulting from
looking at the restriction of this construction to $S\cap (Y^i,
V_+^i,V_{-}^i)$. $\tilde{X}$ is thus an infinite cyclic cover of $X$ as
claimed. We note once again that $H_i(D-V)\cong H_i(Y)$ and that $
H_i(S-F)\cong H_i(Z)$.

The usual considerations (see \cite{L66} and Section \ref{S: covering construction}) now allow us to set up the Mayer-Vietoris sequences for $\td{C}$ and $(\td{C},\td{X})$:
\begin{equation}\label{E:MV1*}
\begin{CD}
\to H_q(V;\Q)\otimes_{\Q}\Gamma@>d^1_q>>
@>>>H_q(Y;\Q)\otimes_{\Q}\Gamma@>e^1_q>>H_q(\td{C};\Q)@>>>
\end{CD}
\end{equation}
and
\begin{equation} \label{E:MV2*}
\begin{CD}
\to H_q(V,F;\Q)\otimes_{\Q}\Gamma@>d^2_q>>
H_q(Y,Z;\Q)\otimes_{\Q}\Gamma@>e^2_q>>H_q(\td{C},\td{X};\Q)\to.\\
\end{CD}
\end{equation}

We will see later that $d^i_q$, $i=1,2$, is a monomorphism for $0\leq q<n-1$.
Hence $e^i_q$ is an epimorphism, $0<q<n-1$, and the $d^i_q$ provide
presentation matrices for the homology modules of the covers. In fact, the
surjectivity of the $e^i_1$ and the equivalent injectivity of the $d^i_0$
follows from standard connectedness considerations or by replacing homology
with reduced homology, so we need only show injectivity of $d^i_q$, $0<q<n-1$.
That the $d^1_q$ are square matrices in this range follows from:

\begin{proposition}\label{P:d1 square*}
$H_i(Y;\Q)\cong H_i(V;\Q)$, $0<i<n-1$.
\end{proposition}
\begin{proof}
In the proof (and often from here out) we suppress the rational coefficients for simplicity of notation.
 
On the one hand, for $0<i<n-1$:
{\footnotesize
\begin{align*}
H_i(Y)& \cong  H_i(D-V) & &\text{by definition of Y}\\
&\cong H_i(D-(V\cup S)) & &\text{by homotopy equivalence}\\
&\cong H_i(S^n-W) && \text{from the definition of $W$}\\
& \cong H_{n-i-1}(W) & &\text{by Alexander duality.}
\end{align*}}

On the other hand, for $0<i<n-1$:
{\footnotesize
\begin{align*}
H_i(V)&\cong H_{n-i-1}(V,\bd V) & &\text{Poincare Duality and the universal coefficient theorem}\\
&\cong H_{n-i-1}(W, \bd V\cup \overline{N(\Sigma)}) & &\text{by excision}\\ 
&\cong H_{n-i-1}(W, K\cup \overline{N(\Sigma)}) &&\text{by the definitions of the spaces}\\
&\cong H_{n-i-1}(W, S^{n-2}) & &\text{by homotopy equivalence
($\overline{N(\Sigma)}$ collapses to $\Sigma$)}\\
&\cong H_{n-i-1}(W) & &\text{for $i\neq 1$ by the reduced long exact sequence
of the
pair.}
\end{align*}}
For $i=1$, we examine the top of the long exact sequence of the pair
$(W,S^{n-2})$:

\begin{equation*}
\begin{CD}
&@>>>&H_{n-1}(W)&@>>>& H_{n-1}(W,S^{n-2})&@>>>& H_{n-2}(S^{n-2})\\
&@>>>& 
H_{n-2}(W)&@>>>&
H_{n-2}(W,S^{n-2})&@>>>&H_{n-3}(S^{n-2})&@>>>.
\end{CD}
\end{equation*}
Of course $H_{n-3}(S^{n-2})=0$ and $H_{n-2}(S^{n-2}) \cong \Q$. We claim that
$H_{n-1}(W)\cong 0$ and $H_{n-1}(W,S^{n-2})\cong \Q$. This will suffice because
any injection $\Q$ to $\Q$ must be an isomorphism. 

$H_{n-1}(W,S^{n-2}) \cong \Q$ because, as above, $H_{n-1}(W,S^{n-2})\cong
H_{n-1}(V,\bd V)$, which is isomorphic to $\Q$ since $V$ is a connected
$(n-1)$-manifold with boundary. To see that $H_{n-1}(W)\cong 0$, consider the
Mayer-Vietoris sequence of $W\cong V\cup_{F} \overline{N(\Sigma)}$. We know
that: $H_{n-1}(V)= 0$, since $V$ is an $(n-1)$-manifold with
boundary; $H_{n-1}(\overline{N(\Sigma)})= 0$, since $\overline{N(\Sigma)}$ is
homotopy equivalent to a complex ($\Sigma$) of dimension less than $n-1$; and
$H_{n-2}(F)=0$, since $F$ is an $(n-2)$-manifold with boundary. Therefore,
$H_{n-1}(W)= 0$ by the Mayer-Vietoris sequence. 

Therefore, $H_i(V)\cong H_{n-i-1}(W)\cong H_i(Y)$ for all $i$, $0<i<n-1$. 
\end{proof}

We will show later that it is also true that $H_i(Y,Z;\Q)\cong H_i(V,F;\Q)$, $0<i<n-1$, and so the maps $d^2_q$ also give square presentation matrices. 

It follows from the construction of the covering and the action of the covering translation, $t$, that the maps can be written as 
\begin{align*}
d^i_q (\alpha \otimes 1)&=i_{-*}(\alpha)\otimes t -i_{+*}(\alpha)\otimes 1 \\
 &= t(i_{-*}(\alpha)\otimes 1) -i_{+*}(\alpha)\otimes 1, 
\end{align*}
where $\alpha\in H_q(V;\Q)$ or $H_q(V, F;\Q)$ according to whether $i=1$ or $2$, and $i_{\pm}$ correspond to the identification maps of $(V,F)$ to $(V_{\pm},F_{\pm})$ obtained by pushing chains out along the bicollar.

To identify these maps (and their matrices) more specifically, we will turn from the context of $D$, which was useful for the the geometric construction, back to the context of $S^n$, which will be more useful for the following algebraic constructions. In particular, we will make use of the facts that 
\begin{enumerate}
\item $H_i(Y)\cong H_i(S^n-W)$ since the two spaces are homotopy equivalent,
\item $H_i(V,F)\cong H_i(W, \overline{N(\Sigma)})\cong H_i(W, \Sigma)$ by excision and homotopy equivalence,
\item $H_i(V)\cong H_i(S^n-\Omega)$ by homotopy equivalence, and
\item $H_i(Y,Z)\cong H_i(\Omega, \overline{N(\Sigma)}) \cong H_i(\Omega,\Sigma)$ by excision and homotopy equivalence.  
\end{enumerate}
We also need to define maps $j_{\pm}:W\to \Omega$ which extend the
maps $i_{\pm}$ which push $V$ out along its collar isotopically. Let
$N'(\Sigma)$ be another
regular neighborhood of $\Sigma$ such that $\overline{N'(\Sigma)}$ lies in
the interior of $N(\Sigma)$.
Then the closure of $N(\Sigma)-N'(\Sigma)$ is a collar of $\bd
\overline{N(\Sigma)}$ using
the ``generalized annulus property'' (see \cite[Proposition 1.5]{Sto}). Define
$j_{\pm}$ to be $i_{\pm}$ on $V$ and the identity on $N'(\Sigma)$. Extend it to
$N(\Sigma)-N'(\Sigma)$ as the homotopy induced on $\bd \overline{N(\Sigma)}$ by
$i_{\pm}$.
It is easily seen that with the canonical identifications of homology groups
above, $i_{\pm*}: H_i(V)\to H_i(Y)$ corresponds to $j_{\pm*}:H_i(S^n-\Omega)\to
H_i(S^n-W)$ and $i_{\pm*}:H_i(V,F)\to H_i(Y,Z)$ corresponds to $j_{\pm*}:
H_i(W,\Sigma)\to H_i(\Omega,\Sigma)$. This follows by making the correct
identifications at the chain level. Of course we also get maps $j_{\pm *}:
H_i(W)\to H_i(\Omega)$. Therefore, to study the matrices $d^1$ and $d^2$ we can
use
\begin{align*}
d^i (\alpha \otimes 1)&=j_{-*}(\alpha)\otimes t -j_{+*}(\alpha)\otimes 1 \\
 &= t(j_{-*}(\alpha)\otimes 1) -j_{+*}(\alpha)\otimes 1, 
\end{align*}
where $\alpha\in H_q(S^n-\Omega;\Q)$ or $H_q(W, \Sigma;\Q)$ according to whether $i=1$ or $2$.

We will make use of the rational perfect linking pairings (suppressing the ``$\Q$'' in the homology notation)
\begin{align*}
L': &H_p(W)\otimes H_q(S^n-W)\to \Q \\ 
L'': &H_p(\Omega)\otimes H_q (S^n-\Omega) \to \Q, 
\end{align*}
$p+q=n-1$ and $0<p<n-1$, which derive from the perfect intersection pairings
\begin{align*}
\cap: &H_{p+1}(S^n, W)\otimes H_q(S^n-W)\to \Q \\ 
\cap: &H_{p+1}(S^n,\Omega)\otimes H_q (S^n-\Omega) \to \Q 
\end{align*}
and the isomorphisms $H_{p+1}(S^n,W)\cong H_p(W)$ and $H_{p+1}(S^n,\Omega)\cong
H_p(\Omega)$, $0<p<n-1$, obtained from the long exact sequences of the pairs.
Note, to be technically precise, there is the issue that $\Omega$ is an open
set and not a closed subcomplex of $S^n$, but we can get around this by
replacing $Y=D-V$ with the homotopy equivalent $D-N(V)$ (the neighborhood taken
in $D$) and then forming $\Omega$ from this $Y$ and $\overline{N(\Sigma)}$ as
above. This new $\Omega$ will be homotopy equivalent to the old one but have
the benefit of being a closed subcomplex. We will usually avoid the distinction
since the two versions are equivalent for homological purposes.  Recall also
that these pairings are induced (after tensoring with the rationals) from
perfect pairings on the integral homology groups modulo their torsion
subgroups.

Given $r\in H_p(W;\Q)$ and $s\in H_q(S^n-\Omega;\Q)$, we have
\begin{equation}\label{E:commutivity*}
L'(r \otimes j_{-*}(s))= L''(j_{+*}(r) \otimes s).
\end{equation}
This can be seen as follows: Any chain representing $s$ (which lies in
$W-\overline{N(\Sigma)}$  by the definition of $\Omega$) gets pushed into
$S^n-W$ under
$j_{-}$ and the the linking form is the intersection of this chain with a
chain, $R$, representing the isomorphic image of $r$ in $H_{p+1}(S^n,W;\Q)$
(see \cite[Appendix]{GBF}). The latter chain can be taken as some chain $R$ in
$S^n$ whose boundary is a chain representing $r$. Now, under the isotopy of
$S^n$ which takes $W$ to $j_{+}(W)$ and $j_-(W)$ to $W$, the chain representing
$s$ gets pushed back into $W$ and $R$ gets pushed into a chain in $S^n$ whose
boundary is $j_+$ of the chain representing $r$ (which represents $j_{+*}(r)\in
H_p(\Omega)$). Thus this isotopy induces maps which take $j_{-*}(s)$ to $s$ and
$r$ to $j_{+*}(r)$, but since the geometric relationship between the chains is
unaffected by isotopy, the intersection number is unaffected. The formula then
follows immediately using the definitions of $L'$ and
$L''$ in the appendix. Similarly, we get
\begin{equation}\label{E:commutivity2*}
L'(r \otimes j_{+*}(s))=L''(j_{-*}(r) \otimes s).
\end{equation}

We will need one other property of linking numbers. Given $r$ and $s$ as above
\begin{align}\label{E:link int*}
L'(r\otimes j_{-*}(s))-L'(r\otimes j_{+*}(s))&=r\cap s\\
L''(j_{-*}(r)\otimes s)-L''(j_{+*}(r)\otimes s)&=r\cap s, \label{E:link int2*}
\end{align}
where $r\cap s$ is the intersection pairing of $r$ and $s$ on $W$. The
geometric proof is analogous to that in the usual case \cite[p. 542]{L66}.
Recall that $r\in H_p(W;\Q)\cong H_p(V,\bd V;\Q)$ (see the proof of Proposition
\ref{P:d1 square*}) and $s\in H_q(S^n-\Omega;\Q)\cong H_q(V;\Q)$. We claim that
this intersection pairing is equivalent to the perfect intersection pairing
$\cap: H_p(V,\bd V; \Q)\otimes H_q(V;\Q) \to \Q$ (using $V\sim_{h.e.} V-\bd
V$). In either case the intersection pairing is given by the sum of signed
point intersections (assuming general position) of chains in the manifold
$W-N(\Sigma)$. If $\bar s$ is a chain representing $s$, then, since it lies in
$S^n-\Omega\subset V$, it also represents
the corresponding class in $H_q(V;\Q)$. Meanwhile, by tracing back what happens
at the chain level in the equations of the second half of the proof of
Proposition \ref{P:d1 square*}, any chain $\bar r$ representing $r$ also
represents its image under the isomorphism $H_p(W;\Q)\cong H_p(V,\bd V;\Q)$.
But none of this affects the geometric intersection, and the choice of chain is
irrelevant since the intersection pairing is well-defined up to homology.
Therefore, the pairings correspond under the isomorphisms.

We now show the following:
\begin{proposition}\label{square2*}
$H_m(W, \Sigma; \Q)\cong H_m(\Omega, \Sigma; \Q)$, $0<m<n-1$, and hence
$H_m(V,F;\Q)\cong H_m(Y,Z;\Q)$ and $d^2_m$ is a square matrix in the same range. 
\end{proposition}
\begin{proof}
Again we suppress the ``$\Q$'' in the proof for notational convenience. 

We begin with the claim that 
\begin{align*}
H_m(W,\Sigma)&\cong H_m(W)\oplus \td
H_{m-1}(\Sigma)\\
\intertext{and}
H_m(\Omega,\Sigma)&\cong H_m(\Omega)\oplus \td
H_{m-1}(\Sigma)
\end{align*}
for $0<m<n-1$. This will follow from the fact that, for
$0<m<n-1$, the inclusion map $i_*$ of each of the long
exact reduced homology sequences of the pairs $(W,\Sigma)$ and $(\Omega,
\Sigma)$ is the $0$
map. For $i_*: H_m(\Sigma)\to H_m(W)$, $0<m<n-2$, this follows because the
inclusion map can be factored $\Sigma\hra K=S^{n-2}\hra W$ since $\Sigma
\subset K \subset W$. Then $i_*$ factors through $H_m(S^{n-2})$ which is $0$ in
the appropriate range. For $m=n-2$, the equation still holds from the long
exact sequence since $\Sigma$ has dimension $n-k$, $k\geq 4$, so that
$H_{n-2}(\Sigma)\cong 0$. The idea is the same for $i_*: H_m(\Sigma)\to
H_m(\Omega)$ except that a little more care must be taken to identify the
$S^{n-2}$ that the inclusion factors through. This can be done by employing one
of the maps $j_{\pm}$ to $K$. Since the $j_{\pm}$ are ends of isotopies they
take the
$(n-2)$-sphere $K$ to another $(n-2)$-sphere. But by the properties of the
map, $j_{\pm}K\subset \Omega$ and $\Sigma \subset j_{\pm}K$. Thus, we can
conclude the homology arguments just as in the previous case.

Now, from the proof of Proposition \ref{P:d1 square*}, we have
$H_m(W)\cong H_{n-m-1}(V)$, $0<m<n-1$; we know that $H_{n-m-1}(V)\cong H_{n-m-1} (S^n-\Omega)$ by homotopy equivalence; and there is a perfect linking pairing between $H_{n-m-1}(S^n-\Omega)$ and $H_m(\Omega)$ which gives an isomorphism since each is a finitely generated vector space. Putting these together with the above identities establishes the proposition. 
 \end{proof}

We next study the maps $d^q_1$ and $d^q_2$ using $j_{\pm *}$. Recall that 
\begin{align}\label{E:d1 and d2*}
d_q^i (\alpha \otimes 1)&=j_{-*}(\alpha)\otimes t -j_{+*}(\alpha)\otimes 1 \\
 &= t(j_{-*}(\alpha)\otimes 1) -j_{+*}(\alpha)\otimes 1, \notag
\end{align}
where $\alpha\in H_q(S^n-\Omega;\Q)$ or $H_q(W, \Sigma;\Q)$ according to whether $i=1$ or $2$.

Let $\{\alpha^p_i\}$, $\{\beta^q_i\}$, $\{\gamma^p_i\}$, and $\{\delta^q_i\}$ represent dual
bases for $H_p(W;\Z)$, $H_q(S^n-W;\Z)$, $H_p(\Omega;\Z)$, and $H_q(S^n-\Omega;\Z)$, all modulo
torsion, so that
\begin{equation}\label{E:dualbasis*}
L'(\alpha^p_i \otimes \beta^q_j)=L''(\gamma^p_i \otimes \delta^q_j)=\delta_{ij},
\end{equation} 
where $\delta_{ij}$ is here the delta function (i.e. $1$ if $i=j$ and $0$ otherwise). These collections then also form bases for the rational homology groups that result by tensoring with $\Q$, and the relations \eqref{E:dualbasis*} hold under the induced perfect rational pairing.

Let $\{\xi^p_i\}$ be a basis for $\td H_p(\Sigma)$. Then, letting
$\{\bar\alpha^p_i\}$ and $\{\bar\xi^{p-1}_i\}$ represent the bases $\{\alpha^p_i\}$ and $\{\xi^{p-1}_i\}$ under their isomorphic images as direct summands in $H_p(W, \Sigma)$ (see the proof of Proposition \ref{square2*}), $\{\bar\alpha^p_i\}$ and $\{\bar\xi^{p-1}_i\}$ taken together form a basis for $H_p(W, \Sigma)$. Similarly, we define  $\{\hat \gamma^p_i\}$ together with $\{\hat \xi^{p-1}_i\}$ forming a basis for $H_p(\Omega, \Sigma)$.

Let
\begin{align*}
j_{+*}(\delta^q_j)&=\sum_i \lambda_{ij}^q \beta_i^q\\
j_{-*}(\delta^q_j)&=\sum_i \sigma_{ij}^q \beta_i^q\\
j_{+*}(\bar\alpha^q_j)&=\sum_i \mu_{ij}^q \hat\gamma_i^q+ \sum_i e_{ij}^q \hat\xi_i^{q-1} \\
j_{-*}(\bar\alpha^q_j)&=\sum_i \tau_{ij}^q \hat\gamma_i^q+\sum_i f_{ij}^q \hat\xi_i^{q-1}\\
j_{+*}(\bar\xi^{q-1}_j)&=\sum_i \phi_{ij}^q \hat\gamma_i^q+ \sum_i g_{ij}^q \hat\xi_i^{q-1} \\
j_{-*}(\bar\xi^{q-1}_j)&=\sum_i \psi_{ij}^q \hat\gamma_i^q+\sum_i h_{ij}^q \hat\xi_i^{q-1},
\end{align*}
where the first two equations are maps $H_q(S^n-\Omega)\to H_q(S^n-W)$ and the
rest are maps $H_q(W,\Sigma)\to H_q(\Omega,\Sigma)$. Note that the $\lambda$,
$\sigma$, $\mu$, and $\tau$ will all be integers (by the chain map
interpretation of $j_{\pm}$ and the fact that the $\alpha$ and $\delta$ were
initially chosen as generators of the torsion free parts of the appropriate
integral homology groups).

\begin{lemma}
In the above equations, all of the $e_{ij}$ and $f_{ij}$ are $0$ and each $g_{ij}=h_{ij}=\delta_{ij}$ (i.e. $1$ if $i=j$ and $0$ otherwise).
\end{lemma}
\begin{proof}
The proof comes from studying the action of $j_{\pm}$ on chain representatives
of the $\bar\alpha$ and the $\bar \xi$.

First, since the $\bar \alpha$ come from $H_q(W)$ under the standard map
$H_q(W)\to H_q(W,\Sigma)$ induced by projection of chain complexes, each
$\bar\alpha$ can be represented by a chain $a$ mod $C_*(\Sigma)$ (where
$C_*(\Sigma)$ is the chain complex of $\Sigma$) and such that $a$ is a cycle
\emph{in $W$}. But since $j_{\pm }$ is induced by an isotopy  which
fixes $\Sigma$, the image of each $\bar \alpha$ should also have such a
representation, i.e. $j_{\pm}a$ is a cycle in $\Omega$. But by a similar
argument from the geometric underpinnings of the boundary map $\bd_*:
H_q(\Omega,\Sigma)\to \td H_{q-1}(\Sigma)$, each $\hat \xi$ is represented by a
chain $x$ mod $C_*(\Sigma)$ with $\bd x$ a non-zero cycle representing a basis
element of $H_{q-1}(\Sigma)$. Since none of these can occur in the image of
$a$ under $j_{\pm}$ (since $\bd j_{\pm} a= j_{\pm} \bd a=0$), each of the $e$
and $f$ must be zero.

For the last pair of maps, we observe similarly that each $\bar \xi$ is represented by a chain $y$
mod $C_*(\Sigma)$ with $\bd y$ a non-zero cycle representing a basis element of $H_{q-1}(\Sigma)$.
We have $\bd j_{\pm} y=j_{\pm} \bd y =\bd y$ since $j_{\pm}$ fixes $\Sigma$. Since $\bar \xi$ and
$\hat \xi$ are both induced by the same basis for $\td H_{q-1}(\Sigma)$ and since the $\hat \xi$
component of any element of $H_*(\Omega, \Sigma)$ is determined by its image under the boundary map
$\bd_*$, it is clear that each $\bar \xi$ maps to an element whose component in the summand $\td
H_{*-1}(\Sigma)$ of $H_*(\Omega,\Sigma)$ is the corresponding $\hat \xi$.
\end{proof}

Thus we have 
\begin{align*}
j_{+*}(\delta^q_j)&=\sum_i \lambda_{ij}^q \beta_i^q\\
j_{-*}(\delta^q_j)&=\sum_i \sigma_{ij}^q \beta_i^q\\
j_{+*}(\bar\alpha^q_j)&=\sum_i \mu_{ij}^q \hat\gamma_i^q \\
j_{-*}(\bar\alpha^q_j)&=\sum_i \tau_{ij}^q \hat\gamma_i^q\\
j_{+*}(\bar\xi^{q-1}_j)&= \sum_i \phi_{ij}^q \hat\gamma_i^q+\hat \xi^{q-1}_j\\
j_{-*}(\bar\xi^{q-1}_j)&= \sum_i \psi_{ij}^q \hat\gamma_i^q+\hat \xi^{q-1}_j.
\end{align*}

Using \eqref{E:d1 and d2*} we have therefore that the matrices $d^1_q$ can be written as 
\begin{equation*}
d^1_q(t)=(t\sigma_{ij}^q-\lambda_{ij}^q),
\end{equation*}
while $d^2_q$ has the form
\begin{equation*}
\left(
\begin{matrix}
P^q & 0\\
Q^q & R^q
\end{matrix}
\right),
\end{equation*}
where $P^q$ is the matrix given by
\begin{equation*}
P^q=(t\tau_{ij}^q-\mu_{ij}^q)
\end{equation*}
and $R^q$ is the block sum of $\td B_{q-1}$ copies of $t-1$ along the diagonal, where $\td B_{q-1}$ is the reduced $(q-1)$-Betti number of $\Sigma$. 

Since the Alexander polynomials in which we are interested are the
determinants of these presentation matrices, the $\phi_{ij}$ and
$\psi_{ij}$ of $Q$ are at present irrelevant (since by elementary linear
algebra,
\emph{det}$(d^2_q)=$\emph{det}$(P^q)$\emph{det}$(R^q)=$
\emph{det}$(P^q)(t-1)^{\td B_{q-1}}$). To
determine the relationships amongst the $\lambda$, $\mu$, $\sigma$, and
$\tau$, we identify the chains representing the $\bar\alpha$ with those
representing $\alpha$ as above, and similarly for the $\hat\gamma$ and the
$\gamma$. Then we can apply the linking pairings to get:
\begin{align*}
L'(\alpha_k^p\otimes j_{+*}(\delta^q_j))&=\sum_i \lambda_{ij}^q L'(\alpha_k^p\otimes \beta_i^q)=\lambda_{kj}^q\\
L'(\alpha_k^p\otimes j_{-*}(\delta^q_j))&=\sum_i \sigma_{ij}^q L'(\alpha_k^p\otimes \beta_i^q)=\sigma_{kj}^q\\
L''(j_{+*}(\alpha^q_j)\otimes \delta_k^p)& =\sum_i \mu_{ij}^q
L''(\gamma_j^q\otimes \delta_k^p)=\mu_{kj}^q\\
L''(j_{-*}(\alpha^q_j)\otimes\delta_k^p )&=\sum_i \tau_{ij}^q
L''(\gamma_j^q\otimes \delta_k^p)=\tau_{kj}^q.
\end{align*}
We can now use \eqref{E:commutivity*} and \eqref{E:commutivity2*} to obtain $\sigma^q_{jk}=\mu^p_{kj}$ and $\lambda^q_{jk}=\tau^p_{kj}$. This implies that $P^q(t)= -t d^p_1(t^{-1})'$, where $'$ indicates transpose.

It remains only to prove that the $d^i_q$, $0<q<n-1$, are non-singular and \emph{det}$(d^1_q(1))=\pm 1$. As noted above, the first will show that the $d^i_q$ are presentation matrices. The theorem will then follow by taking determinants. 

\begin{lemma}
The $d^i_q$, $0<q<n-1$, are non-singular, and, in particular,
\emph{det}$(d^1_q(1))=\pm 1$.
\end{lemma}
Certainly if $d^1_q$ is nonsingular, then $d^2_p$ is nonsingular, since, up to sign, the determinant
of $d^2_p$ will be $(t-1)^{\td B_{p-1}}$ times some power of $t$ times the determinant of $d^1_q(t^{-1})$, and the last will be nonsingular if $d^1_q(t)$ is. Therefore, it remains to show that the $d^1_q$, $0<q<n-1$, are non-singular. We will show that $d^1_q(1)$ has determinant $\pm 1$ which will establish the claim.
\begin{align*}
d^1_q(1)&= (\sigma_{ij}^q-\lambda_{ij}^q)\\
&=( L'(\alpha_i^p\otimes j_{-*}(\delta^q_j))- L'(\alpha_i^p\otimes j_{+*}(\delta^q_j)),
\end{align*}
but by the properties of the linking pairings above, this is the matrix of the
perfect intersection pairing between $H_p(W;\Q)$ and $H_q(S^n-\Omega;\Q)$, which
is equivalent to the perfect intersection pairing between $H_p(V,\bd V;\Q)$ and
$H_q(V;\Q)$. In fact, since we have chosen generators corresponding
to the generators of
the integral homology groups modulo torsion and since $H_p(W;\Z)\cong
H_p(W,\bd W;\Z)$ and $H_q(S^n-\Omega;\Z)\cong H_q(V;\Z)$ (the latter by
homotopy equivalence and the former as in the proof of Proposition \ref{P:d1
square*}), this is the matrix of
the \emph{integral} perfect intersection pairing between the torsion free parts
$H_p(V,\bd V;\Z)$ and $H_q(V;\Z)$. Therefore, this matrix is
unimodular over
$\Z$ and has determinant $\pm 1$. $\qedsymbol$

\subsubsection{Polynomials of the boundary ``knot''}\label{S: boundary
knot*}

We now wish to study the properties of the polynomials associated to the homology of the
infinite cyclic cover of $X=\bd
\overline{N(\Sigma)}-K$, or, in other words, the
complement of $K\cap S$ in $S$. Note that $(S, K\cap S)$ is a locally-flat
manifold pair. If $\Sigma$ were a point singularity, this would be the
boundary sphere knot of a slicing locally-flat disk knot (see Section
\ref{S: disk knots}). Note, however, that for the case of multiple point
singularities, we here diverge slightly from our previous treatment. Instead
of linking the point singularities with an arc and considering the regular
neighborhood of that arc, we instead consider the regular neighborhood of the
collection of points. This will consist of a collection of balls, and, in this
case, $(S, K\cap S)$ will be a collection of locally-flat sphere knots in
disjoint spheres. 

In order to study these ``boundary knots'' or ``link knots'', we begin by
examining the
kernels
of the boundary maps of the vertical exact sequences in the following
commutative diagram in which $0<k<n-1$:

{\small
\begin{diagram}[LaTeXeqno]\label{E:grid*}
& \\
& & \dTo<a & &\dTo<b && \dTo<c &\\
{} & \rTo^0 & H_k(V;\Q)\otimes_{\Q}\Gamma & \rTo^{d^1} &  H_k(Y;\Q)\otimes_{\Q}\Gamma & \rTo^{e^1} & H_k(\td C;\Q) & \rTo^0 &{}  \\
& & \dTo<r && \dTo<s && \dTo<u &\\
 {}& \rTo^0 & H_k(V,F;\Q)\otimes_{\Q}\Gamma & \rTo^{d^2} &  H_k(Y,Z;\Q)\otimes_{\Q}\Gamma & \rTo^{e^2} & H_k(\td C,\td X;\Q) & \rTo^0 & {} \\
& & \dTo<{\bd_*}  & & \dTo<{\bd'_*} & & \dTo<{\bd''_*} \\
 {}& \rTo^{f^3} & H_k(F;\Q)\otimes_{\Q}\Gamma & \rTo^{d^3} &  H_k(Z;\Q)\otimes_{\Q}\Gamma & \rTo^{e^3} & H_k(\td X;\Q) & \rTo^{f^3} & {}\\
& & \dTo<a & &\dTo<b && \dTo<c &\\
&&&&&&&\hfill .\\
\end{diagram}}

\noindent The top two rows are the exact rows of the Mayer-Vietoris sequences
constructed in Section \ref{S: duality*}. The bottom row is the Mayer-Vietoris
sequence of the restriction to $S$ of the
construction which gives us the top row. The columns are the usual
long exact sequences of the pairs in which the left two have been tensored
with
$\Gamma$ over $\Q$. This preserves exactness since the initial sequence
consists of free modules, in fact vector spaces. Commutativity of the diagrams
is obvious at the chain level. The $0$ maps are a consequence of the
non-singularity of $d^1$ and $d^2$ (see Section \ref{S: duality*}).

Using the exactness, this diagram induces the following commutative diagram:
\begin{equation}\label{E:ker pres}
\begin{diagram}
0&\rTo & \text{cok }a  & \rTo & \text{cok }b & \rTo & \text{cok }c &\rTo & 0\\
&& \dTo^{\cong}& & \dTo^{\cong}& & \dTo^{\cong}&& \\
0&\rTo &  \text{ker }\bd_* & \rTo &  \text{ker }\bd'_* & \rTo &  \text{ker }\bd''_*  &\rTo & 0.
\end{diagram}
\end{equation}
That the left vertical map is an isomorphism follows readily from the
isomorphisms of \emph{ker} $\bd_*$ to \emph{im} $r$ from exactness, the
canonical isomorphism of \emph{coim} $r$ to \emph{im} $r$ induced by $r$, and
the isomorphism \emph{coim} $r$ to \emph{cok} $a$ induced by \emph{ker}
$r=$\emph{im} $a$ from the exactness. The other vertical isomorphisms follow
similarly, and so the sequences are isomorphic. Furthermore, the long
exact sequences of the \emph{rational} homology of the pairs $(V, F)$ and
$(Y,Z)$, as exact sequences of vector spaces, must split at each term; in
other
words, each is isomorphic to an exact sequence of vector spaces of the form
$\to A\oplus B\to B\oplus C\to C\oplus D \to$. This
splitting and exactness is preserved under the tensor product with $\Gamma$
over $\Q$ so that
each of the left two kernels and cokernels in diagram \eqref{E:grid*} is
isomorphic to a
direct
rational vector space summand of the appropriate homology module tensored
with $\Gamma$. Thus each of the left four terms is a free $\Gamma$-module, and
once we show that the rows are exact and the rightmost maps are surjective,
the leftmost non-trivial maps will give us a presentation matrix for
\emph{ker}
$\bd''_*\cong$\emph{cok} $a$. For notational convenience, we relabel to get
the sequence
\begin{equation}
\begin{CD}
E @>>> G @>>> H
\end{CD}
\end{equation}
but leave ourselves free to think of these modules as kernels, images, cokernels, or coimages as the
proper contexts allow.

\begin{lemma}
The following sequence is exact:
\begin{equation}\label{E:kerpres}
\begin{CD}
E @>d>> G @>e>> H @>>> 0.
\end{CD}
\end{equation}
\end{lemma}
\begin{proof}
Thinking of $E$, $G$, and $H$ as the appropriate kernels,
\emph{im}$(d)\subset$\emph{ker}$(e)$ because $d$ and $e$ are induced by $d^2$ and $e^2$ and $e^2d^2=0$ by the exactness of the rows of \eqref{E:grid*}.

We next show that \emph{ker}$(e)\subset\,$\emph{im}$(d)$. Again we think of $E$, $G$,
and $H$ as the appropriate kernels. We will examine the following piece of \eqref{E:grid*}:
\begin{equation}
\begin{CD}
H_k(V, F;\Q)\otimes_{\Q}\Gamma @>d^2>> H_k(Y,Z;\Q)\otimes_{\Q}\Gamma @>e^2>> H_k(\td C,\td X;\Q) \\
@V \bd_*VV  @V \bd'_* VV @V \bd''_* VV \\
H_{k-1}(F;\Q)\otimes_{\Q}\Gamma @>d^3>> H_{k-1}(Z;\Q)\otimes_{\Q}\Gamma @>e^3>> H_{k-1}(\td X;\Q). 
\end{CD}
\end{equation}
Using the definitions of $E$, $G$, and $H$ and the splittings of the left two
vertical columns of \eqref{E:grid*}, we can write this isomorphically as
\begin{equation}
\begin{CD}
\stackrel{E}{\underset{A}{\scriptstyle{\oplus}}}@>d^2>>\stackrel{G}{\underset{B}{\scriptstyle{\oplus}}}@>e^2>>C\\
@V \bd_*VV  @V \bd'_* VV @V \bd''_* VV \\
\stackrel{A}{\underset{P}{\scriptstyle{\oplus}}}@>d^3>>\stackrel{B}{\underset{Q}{\scriptstyle{\oplus}}}@>e^3>>R,\\
\end{CD}
\end{equation}
where the label changes are the obvious ones, and, of course, $H\subset C$.
Note that $A$ and $B$ are the kernels of the maps $a$ and $b$ at the bottom of
\eqref{E:grid*}. Thus since they are all kernels of the appropriate vertical
maps, $d^2(E)\subset G$ and $d^3(A)\subset B$. Therefore $d^2$ and $d^3$ each
have the block forms $\left( \begin{smallmatrix} X & 0\\ Y & Z
\end{smallmatrix} \right) $ in the appropriate bases. By the commutativity,
the lower right submatrix of $d^2$ is the upper left submatrix of $d^3$, at
least up to equivalence under change of bases; the upper left submatrix of
$d^2$ is the matrix $d$ of \eqref{E:kerpres} in the statement of the lemma.  
Observe that $d$ must be represented by a square matrix: if it had more rows
than columns, $d^2$ would have determinant $0$ which is impossible since $d^2$
is nonsingular; if it had more columns than rows, then, using the fact that
$d$ is also the lower right submatrix of a similar block decomposition of
$d^1$ (by simply moving the whole argument up one level of the grid
\eqref{E:grid*}), $d^1$ would have determinant $0$ which is also impossible.
But $d^2$ is square and so its lower right submatrix, say $\delta$ (so that
$d^2=\left( \begin{smallmatrix} d & 0\\ Y & \delta \end{smallmatrix} \right)
$), is also square. In addition, $d$ and $\delta$ must each be nonsingular
since $d^2$
is.

Now, let $x\in G$ be also in \emph{ker}$(e)$. Since $e: G\to H$ is a restriction of $e^2$,
$e^2(x)=0$. Therefore, by the exactness of the Mayer-Vietoris sequence, $x\in$\emph{im}$(d^2)$ and
so $x=d^2(\epsilon +\alpha)$, where $\epsilon\in E$ and $\alpha\in A$. We need to show that $\alpha=0$ so that $x\in$\emph{im}$(d)=$\emph{im}$(d^2|_{E})$. Since $x\in$\emph{ker}$(\bd'_*)$ and $\epsilon\in$\emph{ker}$(\bd_*)$,
\[0=\bd'_* x=\bd'_*d^2(\epsilon+\alpha)=d^3\bd_*(\epsilon+\alpha)=d^3\bd_*\alpha.\]
But $\bd_*\alpha\in$\emph{ker}$(a)=A$ so that $ d^3\bd_*\alpha=\delta\bd_*\alpha$. Since $\delta$ is nonsingular, hence injective, $\bd_*\alpha$ must be $0$, but this is only possible if $\alpha=0$ since $\alpha\in A$ and $A$ is mapped injectively under $\bd_*$. This completes the proof that \emph{ker}$(e)\subset$\emph{im}$(d)$.

Lastly, we show that $e$ is surjective, this time treating $E$, $G$, and $H$ as
the appropriate cokernels. We can make use of the following fact of
homological algebra
\cite[p.3]{IV}: In any exact category, given the commutative diagram with
exact rows

\begin{equation}\label{E:coksnake}
\begin{CD}
W @>>> X @>>> Y @>>> Z \\
@V w VV @V x VV @V y VV @V z VV\\
W' @>>> X' @>>> Y' @>>> Z' \\
\end{CD}
\end{equation}
and such that $z: Z\to Z'$ is injective, the induced sequence  
\begin{equation*}
\begin{CD}
\text{cok }w @>>> \text{cok }x @>>> \text{cok }y 
\end{CD}
\end{equation*}
is exact. 

From \eqref{E:grid*}, we have a commutative diagram with exact rows
\begin{equation}
\begin{CD}
 H_{k+1}(Z;\Q)\otimes_{\Q}\Gamma @>e^3>> H_{k+1}(\td X;\Q) @>f^3>> H_{k}(F;\Q)\otimes_{\Q}\Gamma \\
 @V b VV @V c VV @VVV\\
 H_k(Y;\Q)\otimes_{\Q}\Gamma @>e^2>> H_k(\td C;\Q) @>0>> H_{k-1}(V;\Q)\otimes_{\Q}\Gamma. \\
\end{CD}
\end{equation}
We can truncate this to get a diagram
\begin{equation}
\begin{CD}
 H_{k+1}(Z;\Q)\otimes_{\Q}\Gamma @>e^3>> H_{k+1}(\td X;\Q) @>f^3>> \text{im }f^3 @>>>0\\
 @V b VV @V c VV @V f VV @VVV\\
 H_k(Y;\Q)\otimes_{\Q}\Gamma @>>> H_k(\td C;\Q) @>0>> 0 @>>>0, \\
\end{CD}
\end{equation}
which is still commutative (an easy verification) with exact rows. The map in the last column is an isomorphism, so the fact quoted above gives an exact sequence
\[
\begin{CD}
\text{cok }b  @>e>>\text{cok }c  @>>>\text{cok }f.   
\end{CD}
\]
In this case, \emph{cok}$(f)$ is clearly $0$ and so $e$ is surjective. 
\end{proof}

Thus $d$ gives a presentation matrix for $H$ which we will now study. From
here on, the $\Q$ in the homology notation will once again  be implied but 
not written. We will
also use $E_k$, $G_k$, $H_k$, and $d_k$ when we mean to think of the groups as
the appropriate kernels in the appropriate dimensions and $\mathfrak E_k$,
$\mathfrak G_k$, $\mathfrak H_k$, and $\mathfrak d_k$ when we think of them as
the appropriate cokernels. For $\mathfrak E_k$ and $\mathfrak G_k$, we also
sometimes make the identification of the cokernels with appropriate direct
summands of $H_k(V)\otimes \Gamma$ and $H_k(Y)\otimes\Gamma$ which maps onto
the cokernels under the projection to them (since $E_k$ and $G_k$ are
submodules, they can be automatically identified as summands). Note that $d_k$
is the restriction of $d^2$ to $E_k$, while $\mathfrak d_k$ can be thought of
as $d^1$ acting on the summand $\mathfrak E_k$ followed by the projection to
the summand $\mathfrak G_k$.

From the splitting of the leftmost column of \eqref{E:grid*} before
tensoring with $\Gamma$, there exist vectors space summands $\td E_k$ and
$\Td{\mathfrak E}_k$ in $H_k(V,F)$ and $H_k(V)$, respectively, such that
$E_k=\td E_k\otimes\Gamma$ and $\mathfrak E_k=\Td{\mathfrak E}_k\otimes
\Gamma$. Furthermore, $r$ can be written as $\td r\otimes
\text{\emph{id}}$, where $\td r: H_k(V)\to H_k(V,F)$ is the map of the
long exact sequence induced by inclusion (and induces an isomorphism of
the summands $\td E_k\cong \Td{\mathfrak E}_k$ ). We can make similar
conclusions about $G$ in the second column of \eqref{E:grid*} and carry
over all of the tilde notations.

For what follows, it is once again simpler to make the identifications of
Section \ref{S: duality*}:
$H_k(V;\Q)\cong H_k(S^n-\Omega;\Q)$, $H_k(Y;\Q)\cong H_k(S^n-W;\Q)$,
$H_k(V,F;\Q)\cong H_k(W,\Sigma;\Q)$, and $H_k(Y,Z;\Q)\cong H_k(\Omega,
\Sigma;\Q)$, but for convenience we maintain all of the other labels, both
of submodules and maps, making the suitable identifications.  We continue
to use $\{\alpha^k_i\}$, $\{\beta^k_i\}$, $\{\gamma^k_i\}$, and
$\{\delta^k_i\}$ as bases for $H_k(W)$, $H_k(S^n-W)$, $H_k(\Omega)$, and
$H_k(S^n-\Omega)$, respectively, appropriately dually paired, and
$\{\xi^k_i\}$ as a basis for $\td H_k(\Sigma)$. Recall that
$H_k(W,\Sigma)\cong H_k(W)\oplus \td H_{k-1}(\Sigma)$ and
$H_k(\Omega,\Sigma)\cong H_k(\Omega)\oplus \td H_{k-1}(\Sigma)$,
$0<k<n-1$,
from the proof of Proposition \ref{square2*}. We observe that $\td
E_k\subset H_k(W)\subset H_k(W,\Sigma)$ and $\td G_k\subset
H_k(\Omega)\subset H_k(\Omega,\Sigma)$. In fact, since we have the diagram
\begin{equation}
\begin{CD}
H_k(S^n-\Omega) @> \cong >> H_k(V)\\
@VVV @V \td r VV\\
H_k(W,\Sigma) @> \cong >> H_k(V,F),
\end{CD}
\end{equation}
where the top isomorphism is induced by homotopy equivalence and $r$ is
induced by the chain projection, and since $\td E_k\cong$\emph{im}$(r)$,
then any
element $\varepsilon\in \td E_k$ can be represented by the image of a
cycle (mod $C(F)$) in $V$ which is thus a cycle in $W$. Therefore, the
image of $\varepsilon$ under $\bd_*: H_k(W,\Sigma)\to H_{k-1}(\Sigma)$ is
$0$, which implies that $\varepsilon$ lies in the the $H_k(W)$ summand. The
argument for $\td G_k\subset H_k(\Omega)\subset H_*(\Omega,\Sigma)$ is the
same. Thus $\td E_k$ and $\td G_k$ are contained in the summands spanned
by the $\{\alpha^k_i\}$ and $\{\gamma^k_i\}$, respectively. We can now
prove the following lemma:

\begin{lemma}\label{L: orth pair}
$\td E_p\subset H_p(W,\Sigma)$ and $\Td {\mathfrak G}_q\subset H_q(S^n,W)$ are perfectly dually paired under $L'$, $p+q=n-1$;
$\Td {\mathfrak E}_p\subset H_p(S^n-\Omega)$ and $\td G_q\subset H_q(\Omega, \Sigma)$  are perfectly dually paired under $L''$, $p+q=n-1$. 
\end{lemma}
\begin{proof}
We begin with the latter:

By the preceding discussion and without loss of generality, let us assume
that the $\{\gamma^q_i\}$ are chosen so that the first $m$ form a basis
$\{\td g^q_i\}_{i=1}^m$ for $\td G$.  We claim that the sub-basis
$\{\delta^p_i\}_{i=1}^m$ in $H_p(S^n-\Omega)$, which is dual to the $\{\td
g^q_i\}_{i=1}^m=\{\gamma^q_i\}_{i=1}^m $ under $L''$, can be taken as the
basis for $\Td{\mathfrak{E}}_p$ under the projection from
$H_p(S^n-\Omega)$.

To see this, we first observe
that, up to sign, $L''([v],\td s([y]))=L'(\td r([v]),[y])$ for $[v]\in H_p(S^n-\Omega)$ and
$[y]\in H_q(S^n-V;\Q)$. This follows by considering the definition of the linking
pairings. If $v$ and $y$ are chains representing $[v]$ and $[y]$, then they
also represent $\td r[v]$ and $\td s[y]$ (as relative chains). Then $L''([v], \td s([y]))$ is the intersection
number of $y$ with a chain in $S^n$ whose boundary is $v$, while
$L'(\td r([v]),[y])$ is the intersection number of $v$ with a chain in $S^n$
whose
boundary is $y$. By the properties of intersection numbers, these agree. 

Now suppose that $v$ is an element of $H_p(S^n-\Omega)$ which lies in \emph{ker}$(\td r)$ and
that $\{\td{\mathfrak g}_i\}_{i=1}^m$ are basis elements of $\Td{\mathfrak G}_q$ which map
onto the $\td g_i$. Then
$0=L'(\td r(v),\td{\mathfrak g}_i)=L''(v,\td g_i)$. Therefore, the intersection of
\emph{ker}$(\td
r)$ and the dual space to $\td G_q$ is $0$. Thus, the dual subspace to $\td G_q$, spanned by
$\{\delta^p_i\}_{i=1}^m$, can be chosen as a sub-basis for  $\Td{\mathfrak E}_p$ under projection. In other words, $\{\delta_i^p+\text{\emph{ker}}(\td r)\}_{i=1}^m$ is a basis for a linear subspace of $\Td{\mathfrak E}_p$ and $L''(\delta_i^p+\text{\emph{ker}}(\td r), \gamma_j)=\delta_{ij}$. 
We will be done once we show that \emph{dim}$(\Td{\mathfrak E}_p)=m$. Since the above gives \emph{dim}$(\Td{\mathfrak E}_p)\geq m$, we need only show \emph{dim}$(\Td{\mathfrak E}_p)\leq m$

Observe that the same arguments, suitably but easily modified, apply to a basis
$\{\td \alpha^p_i\}_{i=1}^\mu$ for $\td E_p$ (where we have taken a subbasis of the $\{\alpha^p_i\}$ which span $H_p(W)$) to show that the duals $\{\beta^q_i\}_{i=1}^\mu$ span a subspace of $\Td{\mathfrak G}_q$ under the projection, and $\mu\leq \text{dim }\Td{\mathfrak G}_q$. But then we have 
\[\text{dim }\Td{\mathfrak E}_p=\text{dim } \td E_p =\mu\leq \text{dim
}\Td{\mathfrak G}_q=\text{dim }\td G_q=m,\]  
which is what we needed to show.

This establishes the duality of $\Td {\mathfrak E}_p$ and $\td G_q$. The other statement follows similarly. 
\end{proof}

Using this lemma and once again the fact that $L''([v],\td s([y]))=L'(\td
r([v]),[y])$ for $[v]\in H_p(S^n-\Omega)$ and $[y]\in H_q(S^n-V;\Q)$, we
can choose bases $\{\td\alpha^k_i\}$, $\{\td \beta^k_i\}$,
$\{\td\gamma^k_i\}$, and $\{\td\delta^k_i\}$ of $\td E_k$, $\Td{\mathfrak
G}_k$, $\td G_k$, and $\Td{\mathfrak E}_k$ such that:
\begin{enumerate}
\item $L'(\td\alpha^p_i\otimes \td\beta^q_j)= L''(\td\gamma^p_i\otimes \td\delta^q_j)=\delta_{ij}$, the Kronecker delta function, $p+q=n-1$; and

\item $\td r(\td\delta^k_i)= \td\alpha^k_i$ and $\td s(\td\beta^k_i)= \td\gamma^k_i$.
\end{enumerate}
In fact, we can, for example, start with a basis $\{\gamma_i^q\}$ for $\td
G_q$, dualize it to a basis for $\Td{\mathfrak E}_p$, push these to a
basis
for $\td E_p$ under $r$, and then dualize again to get a basis for
$\Td{\mathfrak G}_q$. That $s$ applied to these last basis elements
returns us to our initial basis is easy to check using the duality and
that $L''([v],\td s([y]))=L'(\td r([v]),[y])$.

With this choice of bases, the diagram
\begin{equation*}
\begin{CD}
\Td{\mathfrak E}_k \otimes \Gamma= \mathfrak E_k @>\mathfrak d_k>> \Td{\mathfrak G}_k \otimes \Gamma= \mathfrak G_k\\
@V\td r \otimes \text{id}=r V \cong V @V\td s \otimes \text{id}=s V \cong V\\ 
\td{E}_k \otimes \Gamma= E_k @>d_k>> \td{G}_k \otimes \Gamma= G_k
\end{CD}
\end{equation*}
makes it clear that as matrices $d_k=\mathfrak d_k$. 

We can now establish duality for the polynomials of the modules $H_k$.

\begin{proposition}\label{P:kersym}
$d_p(t)=-t d_q(t^{-1})'$, $p+q=n-1$, where $'$ indicates transpose.
\end{proposition}
\begin{proof}
By the immediately preceding comment, it suffices to show that $d_p(t)=-t
\mathfrak d_q(t^{-1})'$, where the bases of the modules have been chosen
as in the preceding discussion.

The proof is essentially the same as that of Theorem \ref{T:duality*}:
\begin{align*}
d_q (\tilde \alpha \otimes 1)&=\td j_{-*}(\td \alpha)\otimes t -\td j_{+*}(\td \alpha)\otimes 1 \\
 &= t(\td j_{-*}(\td \alpha)\otimes 1) -\td j_{+*}(\td \alpha)\otimes 1, 
\notag
\end{align*}
where $\td \alpha\in \td E_q$ and $\td j_{\pm}$ indicates the restriction of $j_{\pm}$ to $\td E_q$; and, similarly, 
\begin{align*}
\mathfrak d_q (\td{\delta}) \otimes 1)&=\Td{\mathfrak
j}_{-*}(\td{\delta}\otimes t -\Td{\mathfrak j}_{+*}(\td{\delta})\otimes 1 \\
 &= t(\Td{\mathfrak j}_{-*}(\td{\delta})\otimes 1) -\Td{\mathfrak
j}_{+*}(\td{\delta})\otimes 1, \notag
\end{align*}
where $\td{\delta}\in \Td{\mathfrak E}_q$ and $\Td {\mathfrak j}_{\pm}$ indicates $j_{\pm}$ restricted to the summand $\Td{\mathfrak E}_q$ followed by projection to $\Td{\mathfrak G}_q$ of $H_q(S^n-W)$.

Let
\begin{align*}
\td {j}_{+*}(\td {\alpha}^q_j)&=\sum_i \mu_{ij}^q \td \gamma_i^q\\
\td {j}_{-*}(\td \alpha^q_j)&=\sum_i \tau_{ij}^q \td\gamma_i^q\\
\Td{\mathfrak j}_{+*}(\td{\delta}^q_j)&=\sum_i \lambda_{ij}^q \td{\beta}_i^q\\
\Td{\mathfrak j}_{-*}(\td{\delta}^q_j)&=\sum_i \sigma_{ij}^q \td{\beta}_i^q.
\end{align*}
Then $d_q$ and $\mathfrak d_q$ have the forms
\begin{align*}
d_q(t) &=(t\tau_{ij}^q-\mu_{ij}^q)\\
\mathfrak d_q (t) &=(t\sigma_{ij}^q-\lambda_{ij}^q).
\end{align*}

To determine the relationships amongst the $\lambda$, $\mu$, $\sigma$, and $\tau$, we use the
linking pairings to get: 
\begin{align*}
L'(\td \alpha_k^p\otimes \td{\mf{j}}_{+*}(\td \delta^q_j))&=\sum_i
\lambda_{ij}^q L'(\td\alpha_k^p\otimes \td \beta_i^q)=\lambda_{kj}^q\\
L'(\td\alpha_k^p\otimes \td{\mf{j}}_{-*}(\delta^q_j))&=\sum_i
\sigma_{ij}^q L'(\td\alpha_k^p\otimes \td \beta_i^q)=\sigma_{kj}^q\\
L''(\td j_{+*}(\td\alpha^q_j)\otimes \td\delta_k^p)& =\sum_i \mu_{ij}^q
L''(\gamma_j^q\otimes \td\delta_k^p)=\mu_{kj}^q\\
L''(\td j_{-*}(\td\alpha^q_j)\otimes\td\delta_k^p )&=\sum_i \tau_{ij}^q
L''(\td\gamma_j^q\otimes \td\delta_k^p)=\tau_{kj}^q.
\end{align*}

Once we establish that our previous equations
\begin{align}
L'(\alpha \otimes j_{-*}(\delta))&= L''(j_{+*}(\alpha) \otimes \delta) \label{E:com1*}\\
L'(\alpha \otimes j_{+*}(\delta))&=L''(j_{-*}(\alpha) \otimes \delta), \label{E:com2*}
\end{align}
for $\alpha\in H_p(W;\Q)$ and $\delta\in H_q(S^n-\Omega;\Q)$, are still
applicable for the restricted pairings with $\alpha\in \td E_p$ and
$\delta\in \Td{\mathfrak E}_q$,
we can employ them to obtain $\sigma^q_{jk}=\,\mu^p_{kj}$ and
$\lambda^q_{jk}=\tau^p_{kj}$. This will imply that
$d_p(t)=-t\, \mathfrak d_q(t^{-1})'$, and the proposition will be proved.

We begin once again with the
observation that $j_{\pm *}$ takes elements of $\td E_k$ to elements of $\td G_k$ so that for
$\alpha\in \td E_p$ and $\delta\in \Td{\mathfrak E}_q$,
 \[ L''(j_{\pm*}(\alpha) \otimes 
\delta)= L''(\td j_{\pm*}(\alpha) \otimes \delta) \]
simply as a matter of making the obvious restrictions. On the other hand
$j_{\pm*}(\delta)$ might have components in both $\Td {\mathfrak G}_q$ and
its complementary summand. Since $\Td{\mathfrak j}_{\pm *}$ is $j_{\pm *}$
followed by projection to $\Td {\mathfrak G}_q$, we can write $
j_{\pm*}(\delta)= \Td{\mathfrak j}_{\pm *}(\delta)+x$, where $x$ lies in
the \emph{ker}$(s)$. But $E_p$ and \emph{ker}$(s)$ are orthogonal as in
the proof of Lemma \ref{L: orth pair}, so we have, for $\alpha$ and
$\delta$ as above,
\[ L'(\alpha \otimes j_{-*}(\delta))= L'(\alpha \otimes \Td{\mathfrak j}_{\pm *}(\delta)+x)=
L'(\alpha \otimes \Td{\mathfrak j}_{\pm *}(\delta)) + L'(\alpha \otimes x)= L'(\alpha \otimes \Td{\mathfrak j}_{\pm *}(\delta)).\]
Putting these together with \eqref{E:com1*} and \eqref{E:com2*} gives the desired
\begin{align*}
L'(\alpha \otimes \Td{\mathfrak j}_{-*}(\delta))&=L''(\td j_{+*}(\alpha) \otimes \delta) \\
L'(\alpha \otimes \Td{\mathfrak j}_{+*}(\delta))&=L''(\td j_{-*}(\alpha) \otimes \delta) 
\end{align*}
for $\alpha\in \td E_p$ and $\delta\in \Td{\mathfrak E_q}$. 
\end{proof}

\begin{corollary}\label{E:deltasim}
\emph{det}$(d_p(t))\sim\,$\emph{det}$(d_q(t^{-1}))$, $p+q=n-1$, where $\sim$ indicated the similarity relationship for polynomials in $\Gamma$. 
\end{corollary}
\begin{proof}
This follows immediately from Proposition \ref{P:kersym}  by taking
determinants.
\end{proof}

\begin{theorem}
Recall that $S= \bd \overline{N(\Sigma)}$, $X=S-(K\cap S)$, and $\td X$ is
the infinite cyclic covering of $X$. Let $\nu_i(t)$, $0<i<n-2$, be the
Alexander polynomials of $K\cap S$ in $S$. In other words, $\nu_i(t)$ is
the determinant of the presentation matrix of the $\Gamma$-module $H_i(\td
X; \Q)$. Then $\nu_i(t)=r_i(t)(t-1)^{\td B_i}$, where $\td B_i$ is the
\emph{i}th reduced Betti number of $\Sigma$; $r_P(t)\sim r_Q(t^{-1})$,
$P+Q=n-2$; and, if $\nu_i(t)$ is taken primitive in $\Lambda$, then
$r_i(1)=\pm 1$.
\end{theorem}

\begin{proof}
We will make us of the long exact sequence
\begin{equation}
\begin{CD}
@>>> H_i(\td X) @>u_{i}>> H_i(\td C) @>v_{i}>> H_i(\td C, \td
X) @>\bd_{i*}>>,
\end{CD}
\end{equation}
in which we continue to suppress the $\Q$'s which indicate \emph{rational}
homology.  Observe that the $\Gamma$-module structure is preserved
trivially at the chain level, by interpreting $t$ as the covering
transformation, so that this is an exact sequence of
$\Gamma-$\emph{modules}. Since the $H_i(\td C)$ and $H_i(\td C, \td X)$
are $\Gamma-$torsion modules for $i\leq n-2$ by Theorem
\ref{T:duality*}, the $H_i(\td X)$ must also be $\Gamma-$torsion modules
for $i<n-2$.
Thus, recalling that any module over a principal ideal domain can be given
a square presentation matrix, $\nu_i(t)$ will be well defined as the
determinant of that of $H_i(\td X)$. (Equivalently, we can think of
$\nu_i(t)$ as $\prod \nu_{i_j}
(t)$ where $H_i(\td X) =\bigoplus_j \Gamma/(\nu_{i_j}(t))$.)

Recall that, by Corollary \ref{C: kern}, we know that whenever we have an exact sequence of torsion $\Gamma$-modules, say
\[
\begin{CD}
M_1 @>f_1 >> M_2 @>f_2>> M_3 @>f_3>> M_4, 
\end{CD}
\]
then the determinant of the presentation matrix of $M_2$ is the product of the determinants of the presentation matrices of \emph{ker}$(f_2)$ and \emph{ker}$(f_3)=$\emph{im}$(f_2)$. 

Let $c_k$ be the determinant of the matrix $d_k$ above. With
$\lambda_k(t)$, $\mu_k(t)$, and $\nu_k(t)$ all as above, we have then that
$c_k|\mu_k$ and $c_k|\lambda_k$, each because $d_k$ is the
presentation matrix of $H_k$, the kernel of $\bd_{k*}$. Further, we must
then have that $\frac{\mu_k(t)}{c_k(t)}$ is the determinant of the
kernel of $u_{k-1}$ and $\frac{\lambda_k(t)}{ c_k(t)}$ is the
determinant of the kernel of $v_k$. Thus
\begin{equation*}
\nu_k(t)\sim \frac{\lambda_k(t)}{ c_k(t)} \frac{\mu_{k+1}(t)}{
c_{k+1}(t)}.
\end{equation*}

Recall that 
$\mu_q(t)\sim \lambda_p(t^{-1})(t-1)^{\td B_{q-1}}$, $p+q=n-1$, and we
have just proven in Corollary \ref{E:deltasim} that $c_q(t)\sim c_p
(t^{-1})$. Thus 
\begin{equation*}
\frac{\mu_q(t)}{ c_q(t)} \sim \frac{\lambda_p(t^{-1})(t-1)^{\td B_{q-1}}}{
c_p(t^{-1})}.
\end{equation*}
Further, since $(t-1)^{\td B_{k-1}}|\mu_k(t)$ but $(t-1)\nmid
\lambda_k(t)$ (because $\lambda(1)=\pm 1$), it follows from the above
formula for the decomposition of the determinants that $(t-1)^{\td
B_{k-1}}|\nu_{k-1}(t)$. Therefore, if we take $p+q=n-1$, $P=p$, $Q=q-1$,
and $r_k(t)=\frac{\nu_{k}(t)}{(t-1)^{\td B_{k}}}$, then

\begin{align*}
r_P(t^{-1})\sim \frac{\nu_P(t^{-1})}{(t^{-1}-1)^{\td B_P}}&\sim 
\frac{1}{(t^{-1}-1)^{\td B_p}}\cdot\frac{\lambda_p(t^{-1})}{ c_p(t^{-1})}
\cdot\frac{\mu_{p+1}(t^{-1})}{ c_{p+1}(t^{-1})}\\
&\sim \frac{1}{(t^{-1}-1)^{\td B_p}}\cdot\frac{\mu_q(t)}{
c_q(t) (t-1)^{\td B_{q-1}}}
	\cdot\frac{\lambda_{q-1}(t)(t^{-1}-1)^{\td B_{p}}}{ c_{q-1}(t)}\\
&\sim \frac{1}{(t-1)^{\td B_{q-1}}}\cdot\frac{\mu_q(t)}{ c_q(t) }
	\cdot\frac{\lambda_{q-1}(t)}{ c_{q-1}(t)}\\
&\sim \frac{\nu_Q(t)}{(t-1)^{\td B_Q}} \sim r_Q(t).
\end{align*}

Lastly, we know, again from Corollary \ref{E:deltasim}, that we can take
each $c_k$, $\lambda_k/c_k$, and $\mu_k/c_k$ to be primitive
in $\Lambda$, which will make $\lambda_k$, $\mu_k$, and $\nu_k$ primitive
in
$\Lambda$. Since, in that case, $\lambda_k(1)=\pm 1$, we must have each of
its factors $c_k(1)$ and $c_k(1)/c_k(1)$ equal to $\pm 1$. Also, since
$\frac{\mu_k(t)}{(t-1)^{\td
B_{k-1}}}=\lambda_{n-k-1}(t)$ and $c_k(t)$ are equal to $\pm1$ at $1$, so
must be
$\frac{\mu_k(t)}{(t-1)^{\td B_{k-1}}(t)\Delta_k(t)}$.

But then 
\[r_k(t)=\frac{\nu_k(t)}{(t-1)^{\td B_{k}}}=\frac{1}{(t-1)^{\td
B_k}}\frac{\lambda_k(t)}{ \Delta_k(t)} \frac{\mu_{k+1}(t)}{ \Delta_{k+1}(t)}\]
must also be primitive in $\Lambda$ and evaluate to $\pm 1$ at $1$.

\end{proof}

\begin{remark}
Note that when $n=2q+2$, our duality results and the proof of the theorem imply
that $r_q(t)$ is similar to a polynomial of the form $p(t)p(t^{-1})$,
$p\in \Lambda$.
\end{remark}
\begin{remark}
In the case where the singularity $\Sigma$ is a point, the results of this
section reduce to well-known facts about locally-flat sphere slice knots (see
\cite{L66}, \cite{M66}, \cite{L69}). 
\end{remark}

\subsubsection{The subpolynomials}\label{S: subpolys*}

The same algebraic considerations, which we applied in Sections
\ref{S: real disk knot} and \ref{S: conclusion alexa}
to split the
three sets of
Alexander polynomials of a disk knot into three sets of subpolynomials
and to show that these subpolynomials satisfy their own duality
relationships, readily 
generalize to the case of a knot with more general
singularities. Note
that all of the $(t-1)$ factors are
shared between the relative and boundary polynomials. 

Furthermore, if
$c_k(t)$ is the polynomial factor shared by $H_k(\td C)$ and $H_k(\td C,
\td X)$ (i.e. the polynomial of  \emph{ker}$\left[\bd_*'': H_k(\td C, \td
X;\Q)\to H_{k-1}(\td X;\Q)\right]$), then, for a knot $S^{2q-1}\subset
S^{2q+1}$, we can generalize the necessary conditions  we obtained
for the middle
dimension polynomial,
$c_q(t)$, of a disk knot in Section \ref{S: middim disc}. In fact, if
we replace integral homology and integral pairings with rational
homology and rational pairings, the computations of the presentation 
and pairing matrices goes through unchanged. It is only necessary to note
that, in this context, the pairings $L'$ and $L''$ again induce perfect
pairings between certain kernels and coimages (or cokernels) of diagram
\eqref{E:grid*},
but this is shown in Section \ref{S: boundary knot*}. Therefore, we have
the following proposition:

\begin{proposition}
$H_q=$ \emph{ker}$(\bd_*'': H_q(\td C, \td X;\Q)\to H_{q-1}(\td X;\Q))$ 
has a presentation matrix of the form $\tau t - (-1)^{q+1}
R'\tau'R^{-1}$, where $R$ is the matrix of the map $\td{\mathfrak E}\to
\td E$ induced by $\td r: H_q(V)\to H_q(V, F)$. $\mathfrak H_q=H_q(\td
C;\Q)/$\emph{ker}$(H_q(\td C;\Q)\to H_q(\td C,\td X;\Q )) $ has
presentation matrix $(-1)^{q+1}(R^{-1})'\tau R t-\tau'$. Furthermore,
there is $(-1)^{q+1}$-Hermitian pairing $\blm{\,}{}: \mathfrak H_q \times
\mathfrak H_q \to
Q(\Gamma)/\Gamma$ with matrix representative $\frac{t-1}{(R^{-1})'\tau
-(-1)^{q+1}\tau'tR^{-1}}$ with respect to the appropriate basis.
\end{proposition}

All of these necessary conditions on the polynomials can now be summarized
in the following theorem. The duality conditions on the Alexander
subpolynomials follows from that on the Alexander polynomials as in the
proof of Lemma \ref{L: subduality}. The only change, in fact, is the need
to keep special track of the $(t-1)$ factors, but, as already noted, we
know that these must all divide the $a_i$.

\begin{theorem} \label{T: sstrata}
Let $\nu_j(t)$, $\lambda_i(t)$, and $\mu_i(t)$, $0<j<n-2$ and $0<i<n-1$,
denote the Alexander polynomials corresponding to
$H_j(\td X)$, $H_i(\td C)$, and $H_i(\td C, \td X)$, respectively, of a
knotted $S^{n-2}\subset S^n$. We can assume these polynomials to be primitive in $\Lambda$. Then there exist polynomials $a_i(t)$, $b_i(t)$, and $c_i(t)$, primitive in $\Lambda$, such that
		\begin{enumerate} 
			
		\item $\nu_j(t)\sim a_j(t)b_j(t)$, 
					
		\item $\lambda_i(t)\sim b_i(t)c_i(t)$, 
			
		\item $\mu_i(t)\sim c_i(t)a_{i-1}(t)$,

\item	$a_i(t)\sim b_{n-2-i}(t^{-1})(t-1)^{\td B_i}$,	
\item $c_i(t)\sim c_{n-1-i}(t^{-1})$,

\item $b_i(1)=\pm 1$,
\item $c_i(1)=\pm 1$,

\item if $n=2q+1$, then $c_q(t)$ is the determinant of a matrix of the
form $(R^{-1})'\tau R t-(-1)^{q+1}\tau'$, where $\tau$ and $R$ are  
matrices such that $R$ has non-zero determinant.
 \end{enumerate}

Furthermore, if $n=2q+1$, there is a $(-1)^{q+1}$-Hermitian pairing
$\blm{\,}{}: \mathfrak H_q \times \mathfrak H_q \to Q(\Gamma)/\Gamma$ with
matrix representative $\frac{t-1}{(R^{-1})'\tau -(-1)^{q+1} \tau'tR^{-1}}$
with respect to an appropriate basis, where $\mathfrak
H_q=$\emph{cok}$(c)$ in diagram \eqref{E:grid*} above.
 $\qedsymbol$
\end{theorem}

\begin{remark}
Note that $H_0(\td C;\Q)\cong \Q\cong
\Gamma/(t-1)$, since $\td C$ is connected and the
action of $t$ on $H_0(\td
C;\Q)$ is trivial. Similarly, $H_0(\td X;\Q)\cong \oplus \Q\cong \oplus
\Gamma/(t-1)$, where the number of summands is equal to the number of
components of $\Sigma$. And of course, $H_0(\td C, \td X)=0$. Therefore,
by the long exact polynomial sequence of the knot, it is consistent in the
above theorem to take $a_0(t)\sim (t-1)^{\td B_0}$. 
\end{remark}

\subsubsection{High dimensions}\label{S: high dim *}

For completeness, we observe the following concerning the triviality of the  
knot homology modules in the dimensions above those which we have treated in
detail. We maintain the above notation.

\begin{proposition}\label{P: upper dims}$\,$\newline
\begin{enumerate}
\item $H_i(\td X; \Q)\cong 0$ for  $i\geq n-2$.
\item $H_i(\td C;\Q) \cong 0$ for $i\geq n-1$.
\item $H_i(\td C, \td X; \Q) \cong 0$ for  $i\geq n-1$.
\end{enumerate}
\end{proposition}
\begin{proof}
The assertion for $H_i(\td C,\td X;\Q)$ will follow from the other two
and the long exact sequence of the pair. 

The proposition is trivial for
$i\geq n-1$ in the case of $\td X$ and for $i\geq n$ in the case of $\td
C$
because $\td X$ and $\td C$ are noncompact manifolds of  
respective dimensions $n-1$ and $n$.

To show that $H_{n-1}(\td C;\Q)\cong 0$, we can employ Assertion 9 of Milnor
\cite{M68}, which states that for $M$ a compact triangulated $n-$manifold, $\td
M$ the infinite cyclic cover, there is a perfect orthogonal pairing to $\Q$ of
$H_{i-1}(\td M;\Q)\cong H^{i-1}(\td M;\Q)$ and $H_{n-i}(\td M, \bd\td M;\Q)\cong
H^{n-i}(\td M,\bd\td M; \Q)$, provided that $H_*(\td M;\Q)$ is finitely
generated over $\Q$. We can take $i=1$ and $M=C$ (replacing $C$ by the homotopy
equivalent knot exterior to get compactness).  Then $H_0(\td C, \bd
\td C;\Q)\cong 0$,
which implies $H_{n-1}(\td C;\Q)\cong 0$, provided $H_*(\td M;\Q)$ is finitely
generated over $\Q$. But Assertion 5 of the same paper states that this holds if 
$M$ is a homology circle over $\Q$, which we know to be true by Alexander
duality.

The same argument holds to show that $H_{n-2}(\td X;\Q)\cong 0$ provided $H_*(\td X;\Q)$ is
finitely generated over $\Q$.  We know that $H_i(\td X;\Q)$ is $0$ for $i\geq n-1$ and that it
is a torsion $\Gamma$-module for $0<i<n-2$, which implies that it is finite dimensional over
$\Q$ in this dimension range. $H_0(\td X;\Q)$ is also finite
dimensional, being equal
in dimension to the finite number of components of $\Sigma$. Therefore, it
only remains to show
that $H_{n-2}(\td X;\Q)$ is finite dimensional over $\Q$.  For this, we will show directly
that $H_{n-1}(\td C, \td X;\Q)=0$. Then, because $H_{n-2}(\td C;\Q)$ is finite dimensional (in
fact a torsion $\Gamma$-module), the result for $H_{n-2}(\td X;\Q)$ will follow from the long
exact sequence of the pair.

To prove that $H_{n-1}(\td C, \td X;\Q)=0$, we once again employ the
Mayer-Vietoris sequence of
the covering:
\begin{equation*}
\to H_{n-1}( V, F;\Q)\otimes_{\Q}\Gamma \to H_{n-1}( Y,  Z;\Q)\otimes_{\Q}
\Gamma \to
H_{n-1}(\td C;\td X;\Q) \overset{0}{\to}.
\end{equation*}
The last map can be taken as the zero map because we know that
\begin{equation*}
H_{n-2}( V,F;\Q) \to H_{n-2}(Y, Z;\Q)
\end{equation*}
is injective from
the proof of Theorem \ref{T:duality*}. So the proof will be complete if we
show that $H_{n-1}(Y, Z;\Q)=0$. But we saw in Section \ref{S:
sstrata} that $H_{n-1}(Y, Z;\Q)\cong H_{n-1}(\Omega,\Sigma;\Q)$, and,
since $\Sigma$ is a complex of dimension at most $n-4$, this is isomorphic
to $H_{n-1}(\Omega;\Q)$. Hence it suffices to show that this group is $0$.

Consider the long exact sequence (with rational coefficients suppressed in the notation) 
\begin{equation*}
\begin{CD}
H_n(\Omega)@>>> H_n(S^n)@>>> H_n(S^n,\Omega) @>>> H_{n-1}(\Omega) @>>> H_{n-1}(S^n).
\end{CD}
\end{equation*}
$H_{n-1}(S^n)=0$ and $H_n(S^n)\cong \Q$, trivially, and $H_n(\Omega)=0$
by Alexander duality.
Furthermore, by Lefschetz duality, $H_n(S^n,\Omega)\cong
H^0(S^n-\Omega)\cong \Q$ since $S^n-\Omega\sim_{h.e.} V$ is connected.
Therefore, the long exact sequence reduces to
\begin{equation*}
\begin{CD}
0 @>>> \Q @>>> \Q  @>>> H_{n-1}(\Omega;\Q) @>>> 0.
\end{CD}
\end{equation*}
Since any injective map $\Q\to \Q$ must be an isomorphism, the result
follows.  
\end{proof}

%% file: sstratab.tex

\subsection{Constructions}\label{S: construction*}

For the case of knots with  general singularities, realization of polynomials
is more
difficult than it was for the case of point singularities because the allowable set of
polynomials may depend subtly on the
properties of the singular set, its link pairs, and its embedding. However, in the
following sections, we
will employ several constructions to create knots with singularities and
to compute
their Alexander polynomials. These will provide at least partial realization results.

In Section \ref{S: frame spin}, we will use the \emph{frame spinning} of Roseman
\cite{Ro89} (generalized to spin non-locally-flat knots)  to construct knots with certain
kinds of manifold singularities. In Section \ref{S: twist spin}, we further generalize
this construction to  create \emph{frame twist-spinning}. Together, these
procedures include as
special cases the superspinning of Cappell \cite{C70} and the twist spinning of Zeeman    
\cite{Z65}. 
In Section \ref{S: suspend}, we construct knots by suspension.

\subsubsection{Frame spinning}\label{S: frame spin}

To construct some examples of knots with a given singular stratum, we will employ the technique of
\emph{frame spinning},
which was introduced by Roseman in \cite{Ro89} and studied further by Suciu \cite{Su92} and Klein and Suciu \cite{KS91}.
It generalizes the earlier techniques of spinning and the superspinning of Cappell \cite{C70}. We begin by describing
this procedure.

Let $K$ be a knot $S^{m-2}\subset S^m$, and let $M^k$ be a $k$-dimensional framed
submanifold of $S^{m+k-2}$ with framing $\phi$. Suppose that $S^{m+k-2}$ is embedded in
$S^{m+k}$ by the standard (unknotted) embedding. Roughly speaking, the frame spun knot
$\sigma_M^{\phi}(K)$ is formed by removing a standard disk pair  $(D^{m},
D^{m-2})$ at each point of $M$ and replacing it with the disk knot obtained by removing a
neighborhood of a nonsingular point of the knot $K$.

More specifically, let $(D_{-}^{m}, D_{-}^{m-2})$ be an unknotted open disk pair which is
the open neighborhood pair of a point which does not lie in the singularity  of the
embedding of the knot $K\subset S^{m}$. Let $(D_{+}^{m}, D_{+}^{m-2})=(S^{m}, K)-
(D_{-}^{m},
D_{-}^{m-2})$. This is a disk knot, possibly non-locally-flat, with the unknotted
locally-flat sphere
pair as boundary. Let $M^k\times D^{m-2}$ be the normal bundle of $M^k\in S^{m+k-2}$
determined by the trivialization $\phi$. Finally, writing $S^{m+k}$ as
$S^{m+k}=S^{m+k-2}\times D^2 \cup_{S^{m+k-2}\times S^1} D^{m+k-1}\times S^1$, we let
$S^{m+k-2}\times 0$ in the first factor represent the unknot in which $M$ is embedded.

Now define $\sigma_M^{\phi}(K)$  to be the $(m+k-2)$-sphere 
\begin{equation*}
(S^{m+k-2}-M^k\times \text{ int } D^{m-2}) \cup_{M^k\times S^{m-3}} M^k \times D_{+}^{m-2}
\end{equation*}
embedded in the $(m+k)$-sphere
\begin{equation*}
(S^{m+k}-M^k\times \text{ int } (D^{m-2}\times D^2)) \cup_{M^k\times S^{m-1}} M^k \times D_{+}^{m}.
\end{equation*}
This construction corresponds to removing, for each point of $M$, the trivial disk pair
$(D^{m},
D^{m-2})$, which is the fiber of the normal bundle of $M$, and replacing it
with the knotted disk pair 
$(D_{+}^{m}, D_{+}^{m-2})$.  In the above references, $K$ is always assumed to be a
locally-flat knot, but there is nothing to prevent us from applying this construction
to
a non-locally-flat knot so long as we are careful to embed $(D_{-}^{m}, D_{-}^{m-2})$
away from the singularity. Observe that, in the case where $M^k$ is the sphere $S^k$ with
the standard unknotted embedding and bundle framing, $\sigma^{\phi}_M(K)$ is the
superspin of $K$ (see \cite{C70}).

Let $n=m+k$. We obtain the following:

\begin{proposition}
Let $M^k$ be a manifold which can be embedded in $S^{n-2}$ with trivial normal bundle. Then there exists a knot $S^{n-2}\subset S^{n}$ with $M$ as its only singular stratum.
\end{proposition}
\begin{proof}
Let $K$ be any knot $S^{m-2}\subset S^m$ whose singular set constists of a single point. Let $\phi$ be a trivialization for the normal bundle of the embedding of $M\in S^{n-2}$. Then $\sigma_M^{\phi}(K)$ provides an example.
\end{proof}

To study the Alexander polynomials that occur from such constructions, we first need a
geometric formula for the exterior of a frame-spun knot. This is provided, as follows, by
Suciu in \cite{Su92}, although we adopt our own notations. Throughout this section, let
$X(\cdot)$ denote the exterior of a knot and $\td X(\cdot)$ the corresponding infinite
cyclic covering. $X(K)$, the exterior of the knot $K$, is homeomorphic to the exterior of
the induced disk knot $(D_{+}^{m}, D_{+}^{m-2})$. Its intersection with the exterior of
the induced boundary sphere pair is $D^{m-2} \times S^1$ because the induced boundary
sphere pair of the disk knot is unknotted. Let $ M^k \times \text{int}(D^{m-1})$ represent
the intersection of the tubular neighborhood of $M^k$ in $S^n$ with $D^{n-1}\times 0\in
D^{n-1}\times S^1$, the exterior of the trivial knot $S^{n-2}\subset S^{n}$. It can be
seen that
\begin{equation*}
X(\sigma_M^{\phi}(K)) =(D^{n-1}-M^k \times \text{int}(D^{m-1}))\times S^1
\cup_{M^k\times D^{m-2} \times S^1} M^k\times X(K). 
\end{equation*}

In the following lemma, we use Cov$[\cdot]$ to denote the infinite cyclic covering where the tilde notation would be unwieldly.

\begin{lemma}
{\footnotesize
\begin{align*}
\td X(\sigma_M^{\phi}(K))&
=\text{Cov}[(D^{n-1}-M^k \times \text{int }(D^{m-1}))\times S^1]
\cup_{\text{Cov}[M^k\times D^{m-2} \times S^1]} \text{Cov}[ M^k\times X(K)]\\
&\sim_{h.e.} (D^{n-1} -M^k \times \text{int }(D^{m-1}))\times \R \cup_{M^k\times
D^{m-2}\times \R} M^k\times \td X(K)
\end{align*}}
\end{lemma}
\begin{proof}
As observed in \cite{KS91}, if $V$ is a Seifert surface for the knot $K$ and we define
\begin{equation*}
\sigma_M^{\phi}(V)= (D^{n-1}-M^k \times \text{int}D^{m-1})  \cup_{M^k\times D^{m-2} }
M^k\times V,
\end{equation*}
then $\sigma_M^{\phi}(V)$ is a Seifert surface for $\sigma_M^{\phi}(K)$ (Note that if $K$
is not locally-flat then we mean the knot exteriors and Seifert surfaces in the sense of
Section \ref{S: sstrata}). Using this Seifert surface we can form the infinite cyclic
cover $\td X(\sigma_M^{\phi}(K))$ by the usual ``cut and past'' construction. From this,
the first equation follows by considering what the construction does on each piece. The
second equation follows from the observation that the covering space can be obtained by
``unwrapping'''' $S^1$ to $\R$.
\end{proof}

\begin{remark}
The ability to create the Seifert surface in this manner relies heavily on the following
fact: While the particular framing of $M$ may serve to ``spin'' the knots $K$
tangentially to $S^{n-2}$, the knots are never ``twisted''. No
rotation takes place along the meridians circling $S^{n-2}$ in $S^n$. Thus, contrary to a
remark of Roseman \cite{Ro89}, frame spinning can not yield instances of Zeeman's
\emph{twist spinning} \cite{Z65}. In cases involving twisting, it is not always possible
to get the Seifert surfaces to ``line up'' so that they may be connected by a disk
(although this can happen in special cases, particularly with fibered knots where the
Seifert surfaces can be forced to align by ``rotating them around the fibration'').
However, see the following section (Section \ref{S: twist spin}), in which we introduce a
method to obtain such twisting.
\end{remark}

We can now use a Mayer-Vietoris sequence to study the Alexander modules of 
$\sigma_M^{\phi}(K)$. In particular, we have the rational exact sequence (in
which we suppress the $\Q$'s from the notation)
\begin{equation*}
\to \td H_i(M^k\times D^{m-2}\times \R)\to\td H_i(D^{n-1}\times
\R)\oplus \td H_i(M^k\times \td X(K))\to \td H_i(\td X
(\sigma_M^{\phi}(K)))\to ,
\end{equation*}
in which we have used the homotopy equivalence of $D^{n-1}-M^k \times \text{int
}(D^{m-1})$ and $D^{n-1}$ to replace $\td H_i((D^{n-1}-M^k \times \text{int
}(D^{m-1}))\times \R;\Q)$ with $\td H_i(D^{n-1}\times \R;\Q)$.
This sequence simplifies in the obvious manner to 
\begin{equation}\label{E: MV}
\to \td H_i(M^k;\Q) \overset{i_*}{\to} \td H_i(M^k\times \td X(K);\Q)\to \td  H_i(\td X
(\sigma_M^{\phi}(K));\Q)\to.
\end{equation}
From this, we will prove:
\begin{proposition}
Let $B_i$ be the $i$th Betti number of $M^k$. Let $\lambda_j(t)$ be the $j$th Alexander polynomial of $K$ and $\lambda_i^{\sigma}(t)$ the $i$th Alexander polynomial of $\sigma_M^{\phi}(K)$. Then, for $0<i<n-1$, 
\begin{equation*}
\lambda_i^{\sigma}(t)=\prod_{l=1}^{m-2} [\lambda_l(t)]^{B_{i-l}}.
\end{equation*}

\end{proposition}
\begin{proof}
We first study $H_i(M^k\times \td X(K);\Q)$, which, by the K\"{u}nneth theorem, is 

\begin{equation*}
H_i(M^k\times \td X(K);\Q)\cong \bigoplus_{j+l=i} H_j(M^k;\Q)\otimes_{\Q} H_l(\td X (K);\Q). 
\end{equation*}

 Now let $B_i$ be the $i$th
Betti number of $M^k$, and, for a rational vector space $A$, let $A^i$ denote the direct
sum of $i$ copies of $A$. Since
$A\otimes_{\Q} \Q^i=A^i$,
\begin{equation}\label{E: spin kun}
H_i(M^k\times \td X(K);\Q)\cong \bigoplus_{j+l=i}  [H_l(\td X (K);\Q)]^{B_{j}}.
\end{equation}

To establish the $\Gamma$-module structure, it is evident from the geometry that if
$\alpha\otimes_{\Q} \beta \in H_j(M^k;\Q)\otimes_{\Q} H_l(\td X (K))$, then
$t(\alpha\otimes_{\Q} \beta)= \alpha\otimes_{\Q} (t\beta)$.  Therefore, equation
\eqref{E: spin kun} is a $\Gamma$-module isomorphism.

For the knot $K: S^{m-2}\into S^m$, recall that
\begin{enumerate}
\item $H_0(\td X(K);\Q)=\Q$,
\item $H_i(\td X(K);\Q)=0$, $i\geq m-1$,
\end{enumerate}
by Proposition \ref{P: upper dims}.
Taking this into account,
\begin{equation*}
H_i(M^k\times \td X(K);\Q)\cong \bigoplus_{j+l=i, l<m-1}  [H_l(\td X (K);\Q)]^{B_{j}}.
\end{equation*}

We next study the map $i_*: H_i(M^k;\Q) \to H_i(M^k\times \td X(K);\Q)$ in the Mayer-Vietoris
sequence \eqref{E: MV}. From the geometric constructions above  and consideration of the chain maps used to define the Mayer-Vietoris sequence and K\"{u}nneth theorem, $i_*$ is the map which takes an element $\alpha\in H_i(M^k;\Q)$ to $\alpha \otimes \{*\} $ in the submodule $H_i(M^k;\Q)\otimes H_0(\td X(K);\Q)$ of $ H_i(M^k\times \td X(K);\Q)$, where $\{*\}$ is a point representing the generator of $ H_0(\td X(K);\Q)$. It follows that $i_*$ is injective, and thus the Mayer-Vietoris sequence gives
\begin{equation*}
H_i(\td X (\sigma_M^{\phi}(K));\Q)\cong \bigoplus_{j+l=i, 0<l<m-1}  [H_l(\td X (K);\Q)]^{B_{j}}.
\end{equation*}

So, if $\lambda_j(t)$ is the $j$th Alexander polynomial of $K$, then the $i$th Alexander polynomial of $\sigma_M^{\phi}(K)$ is 
\begin{equation*}
\lambda_i^{\sigma}(t)=\prod_{l=1}^{m-2} [\lambda_l(t)]^{B_{i-l}}
\end{equation*}
since the polynomial associated to a direct sum of torsion $\Gamma$-modules is the product of the polynomials associated to the summands.
\end{proof}

Now assume that the knot $K$ has singular set, $\Sigma$. Then $\sigma_M^{\phi}(K)$ will
have singular set $\Sigma \times M$ stratified by $(\Sigma\times M)_i=\Sigma_{i-k}\times
M$. We can use our previous duality results (Theorem \ref{T:duality*}) to calculate
the relative Alexander polynomials of the pair given by the spun knot complement and the link
pair complement of $\Sigma \times M$. In particular, let $B_i$ continue to denote the
$i$th Betti number of $M$, let $\td {\mathfrak b}_i$ denote the $i$th reduced Betti
number of $\Sigma$, and let $\td \beta_i$ denote the $i$th reduced Betti number of
$M\times \Sigma$. Then, for $i>0$,
\begin{align*}
\mu_i^{\sigma}(t)&\sim (t-1)^{\td \beta_{i-1}}\lambda^{\sigma}_{n-1-i}(t^{-1}) \\
&= (t-1)^{\td \beta_{i-1}}\prod_{l=1}^{m-2} [\lambda_l(t^{-1})]^{B_{n-1-i-l}} \\
&\sim (t-1)^{\td \beta_{i-1}}\prod_{l=1}^{m-2} \left[\frac{\mu_{m-1-l}(t)}{(t-1)^{\td{\mathfrak b}_{m-2-l}}}\right]^{B_{n-1-i-l}}.
\end{align*}
Rather than explore the relations among these Betti numbers directly, we can simplify
this formula by alternatively studying the relative Alexander polynomial directly using a
Mayer-Vietoris sequence, as we did for the $\lambda^{\sigma}_i(t)$. The one significant
difference is that the relative homology module for $K$ in dimension $0$ is $H_0(\td
X(K), \text{Cov}(X\cap \overline{N(\Sigma)}))=0$ so that instead of all of the maps
$i_*$ of the Mayer-Vietoris sequence being injective, they are all $0$ instead. Since
$\td  H_i(M^k;\Q)\cong(\Gamma/(t-1))^{\td B_{i-1}}$, because the $\Gamma$ action of $t$ 
on the cover $M^k\times \R$ induces the identity map on the homology,
we can 
conclude by
polynomial algebra that
\begin{equation*}
\mu_i^{\sigma}(t)= (t-1)^{\td B_{i-1}}\prod_{l=0}^{m-2} [\mu_l(t)]^{B_{i-l}}
\end{equation*}
for $i>0$. Note that $\mu_0= 1$.

Letting $\nu_i(t)$ denote the $i$th Alexander polynomial of the link pair of $\Sigma$
for the knot $K$,
the $i$th Alexander polynomial of the link pair of $\Sigma \times M$ for the knot
$\sigma_M^{\phi}(K)$
is easily derived from the K\"{u}nneth theorem to be
\begin{equation*}
\nu_i^{\sigma}(t)=\prod_{l=0}^{m-3}[\nu_l(t)]^{B_{i-l}}.
\end{equation*} 
Note that $\nu_0(t)=(t-1)^{\td{\mathfrak b}_0+1}$.

As we know, the Alexander polynomials of $\sigma_M^{\phi}(K)$ factor into
Alexander subpolynomials $a^{\sigma}_i(t)$, $b^{\sigma}_i(t)$, and
$c^{\sigma}_i(t)$. It is an exercise with the long exact sequences to show that this
factorization is
preserved under the spinning, modulo some minor extra complication in the $t-1$
factors. In other words,
\begin{align}
b_i^{\sigma}(t)&=\prod_{l=1}^{m-2} [b_l(t)]^{B_{i-l}}\\
c_i^{\sigma}(t)&=\prod_{l=1}^{m-2} [c_l(t)]^{B_{i-l}}\\
a_i^{\sigma}(t)&=(t-1)^{\td \beta_{i}}\prod_{l=1}^{m-2}
\left[\frac{a_l(t)}{(t-1)^{\td{\mathfrak
b}_{l}}}\right]^{B_{i-l}}\label{E: spin a}.
\end{align}
For example, to perform the calculation for the $c_i$, let $\td L(K)$ represent the the infinite cyclic covering of the intersection of the knot exterior $X(K)$ with the closed  neighborhood of the singularity $\bar N(\Sigma)$ (i.e. the ``link exterior'', the usual subset for our relative homology modules) and similarly for $\sigma_M^{\phi}(K)$. The above
calculations show that the kernel module of the map of the long exact sequence
\begin{equation*}
\begin{CD}
H_i(\td X(\sigma_M^{\phi}(K));\Q)@>>> H_i(\td X(\sigma_M^{\phi}(K)), \td L(\sigma_M^{\phi}(K));\Q)
\end{CD}
\end{equation*}
is isomorphic to the kernel of the map 
\begin{equation*}
\begin{CD}
\bigoplus_{j+l=i, 0<l<m-1}[H_l(\td X (K);\Q)]^{B_{j}} @>f>>\bigoplus_{j+l=i,
l<m-1}[H_l(\td X (K),\td L(K);\Q)]^{B_{j}} 
\end{CD}
\end{equation*}
because we know that $H_i(\td X(\sigma_M^{\phi}(K));\Q)$ maps trivially to
the other summand of 
\begin{equation*}
H_i(\td X(\sigma_M^{\phi}(K)), \td
L(\sigma_M^{\phi}(K));\Q),
\end{equation*}
which consists of a sum of $\Q$'s with trivial
$\Gamma$-action, i.e. $\Gamma/(t-1)$'s. (The triviality of this part of
the map, $f$, is a result of the splitting of maps of torsion modules into
their $p$-primary summands (see the proof of Proposition \ref{P:alt.
poly.}).) But from the K\"{u}nneth theorem, this map is induced by the
usual map $p_*: H_l(\td X (K);\Q)\to H_l(\td X (K),\td L(K);\Q)$ in the
long exact homology sequence tensored with the identity map on the
homology of $M$. Since we are working with rational homology, this tensor
product is an exact functor and so the kernel of the map as a rational
vector space is $\oplus H_j(M)\otimes \text{ker }p_*$. But the
$\Gamma$-module structure is also evidently preserved, acting trivially on
the $H_j(M)$ factors and with the action on $\text{ker }p_*$ induced by
that on $H_l(\td X (K);\Q)$. Passing from the modules to the polynomials
gives the above equation for $ c_i^{\sigma}(t)$. The other equations are
handled similarly modulo their $(t-1)$ factors, these terms being
accounted for separately by consideration of what the $t-1$ factors must
be according to the duality formulas of Theorem \ref{T:duality*}.

Perhaps a simpler way to look at  what happens to the Alexander polynomials of a frame spun knot is the following interpretation which follows readily (with a little checking) from the calculations above. Take, for example, the polynomial $\lambda_i(t)$ of the knot $K$ for some $i$. This polynomial will be a factor of $\lambda^{\sigma}_j(t)$ a number of times equal to the the Betti number $B_{j-i}$ of $M$. So if, for example, we take $M=S^k$, each $\lambda_i(t)$ will appear exactly twice as a factor of the $\lambda_j^{\sigma}(t)$, once in its ``native'' dimension $i$ and once $k$ dimensions higher. 
Similar consideration apply for all of the other polynomials and subpolynomials modulo the $t-1$ terms which can be computed at the end by tallying the reduced Betti numbers. 

By taking $\Sigma$ to be a point we can therefore construct a knot with singularity $M$ and certain
specified polynomials as follows:

\begin{proposition}\label{P: frame spin real}
Let $M^k$ be a manifold which embeds in $S^{n-2}$ with trivial normal bundle with framing $\phi$ and
such that $n-k>3$. Let $\Sigma$ be a single point. Let $B_i$ denote the $i$th Betti number of $M$,
and let $\td{\mathfrak b}_i$ and
$\td \beta_i$ denote the $i$th reduced Betti numbers of $\Sigma$ and $M\times \Sigma$, respectively.
Suppose that we are given any set of polynomials, $a_i(t)$, $b_i(t)$, $c_j(t)$ and
$c'_l(t)$,
$0<i<n-k-2$, $0<j<n-k-1$, and $0<l<n-1$, which satisfy:
\begin{enumerate}		
\item	$a_i(t)\sim b_{n-k-2-i}(t^{-1})$,	
\item $c_i(t)\sim c_{n-k-1-i}(t^{-1})$,
\item $c'_i(t)\sim c'_{n-1-i}(t^{-1})$,
\item $b_i(1)=\pm 1$,
\item $c_i(1)=\pm 1$,
\item $c'_i(1)=\pm 1$,
\item if $n-k=2p+1$, $p$ even, $p\neq 2$, then $c_p(t)$ is the determinant of a matrix
of the form $(R^{-1})'\tau R
t-(-1)^{q+1}\tau'$, where $\tau$ and $R$ are integer matrices such that $R$ has non-zero determinant
and
$(R^{-1})'\tau R$ is an integer matrix; if $n-k=2p+1$, $p$ even, $p= 2$, then 
$|c_p(-1)|$ is an odd square,
\item if $n=2q+1$, $q$ even, then $|c'_q(-1)|$ is an odd square. 
\end{enumerate}
Then there exists a knotted $S^{n-2}\subset S^n$ with singular set $M$ and
Alexander subpolynomials $a_i^{\sigma}(t)$,  $b_i^{\sigma}(t)$,
and $c_i^{\sigma}(t)$ satisfying
\begin{align*}
a_i^{\sigma}(t)&\sim (t-1)^{\td \beta_{i}} \prod_{l=1}^{m-2}
\left[a_l(t)\right]^{B_{i-l}}\\
b_i^{\sigma}(t)&\sim\prod_{l=1}^{m-2} [b_l(t)]^{B_{i-l}}\\
c_i^{\sigma}(t)&\sim c'_i(t) \prod_{l=1}^{m-2} [c_l(t)]^{B_{i-l}}. \\
\end{align*}
(The first equation comes from equation \eqref{E: spin a} by taking into account that $\Sigma$ is
being taken as a point. Note that this also implies that
$\td \beta_j$ is just the reduced jth Betti number of $M$.)

\end{proposition}

Fist, we need one lemma regarding disk knots which was not available before our discussion of frame spinning. It will be proven below.
\begin{lemma}\label{L: triv poly}
For any $n\geq 4$, there exists a locally-flat disk knot, $ D^{n-2}\subset D^n$, with
non-trivial
boundary knot and with all Alexander polynomials equal to $1$. Equivalently, there exists
a sphere knot, $S^{n-2}\subset S^n$, with a point singularity and with all Alexander
polynomials equal to $1$.
\end{lemma}

\begin{proof}[Proof of the Proposition]
By the results of Section \ref{S: disk knots} and Lemma \ref{L: triv
poly}, there exists
a sphere
knot, $K$, with a single point singularity and with the desired polynomials
$a(t)$, $b(t)$, and $c(t)$ (if the construction of Section \ref{S: disk
knots}
yields a locally-flat knot, we can take the knot sum with the knot of the lemma).  Then,
by the
calculations above, $\sigma_M^{\phi}(K)$ has the desired $a_i^{\sigma}(t)$,
$b_i^{\sigma}(t)$, and $c_i^{\sigma}(t)$ except for the $c'_i(t)$ factors.  To get the
latter, we can take the knot sum with a locally-flat knot $S^{n-2}\subset S^n$ which has
the $c'_i(t)$ as its Alexander polynomials. Such a knot exists by \cite{L66}.
\end{proof}

\begin{proof}[Proof of Lemma \ref{L: triv poly}]
It is well known that there exist nontrivial locally-flat knots $S^1\subset S^3$ whose
Alexander polynomials are trivial but whose knot groups are nontrivial (see, for example,
\cite{Rolf}). For $n> 4$, we can now frame spin one of these knots about the sphere $S^{n-4}$
with the trivial framing and embedding in $S^{n-3}$ to obtain a locally-flat knot $K:
S^{n-3}\into S^{n-1}$. For $n=4$, we can simply choose $K$ to be the knot with which we
started. In each case the knot is still nontrivial because superspinning preserves knot
groups by Cappell \cite{C70}. Next, we convert $K$ to a disk knot, $L:D^{n-3}\into
D^{n-1}$, by removing a trivial disk pair neighborhood of a point on the knot, just as in
the first step of the frame spinning construction. Lastly, we take as our desired disk
knot the product of the disk knot, $L$, with an interval $I$, $L\times I: D^{n-3}\times I
\into D^{n-1} \times I$. Since its exterior is homotopy equivalent to the exterior of
$L$, all of its Alexander polynomials $\lambda_i(t)$ are trivial. By duality, the
$\mu_i$ are also trivial. The $\nu_i$ are then trivial by polynomial algebra from the
long exact sequence of the knot pair. The
boundary knot is the knot sum $K\#(-K)$, where $-K$ denotes the reflection of
$K$.  Therefore, by the Van Kampen
theorem, the group of the boundary knot is nontrivial and so $L\times I$ is, in fact,
knotted.

We have thus produced a disk knot, $L\times I$, with the desired properties. To obtain
the desired sphere knot with point singularity, we simply take the cone on the boundary
of the disk knot.
\end{proof}

%% file: sstratac.tex

\subsubsection{Frame twist-spinning}\label{S: twist spin}

We now slightly generalize the frame-spinning construction to include
``twisting''. In the special case where we frame twist-spin about a circle,
$S^1$, embedded with standard framing in $S^{n-2}$, we will obtain the twist-spun knots
of Zeeman \cite{Z65}.

Before beginning the construction, we recall that one alternative way to compute
Alexander modules and hence Alexander polynomials is the following: Rather than
considering the infinite cyclic cover of the knot complement and its homology
with rational coefficients, we can instead consider the homology of the knot
complement with a  certain local coefficient system with $\Gamma$ as the stalk.
If $\alpha$ is an element of the fundamental group of the knot complement
and $\ell(\alpha)$ denotes the linking number of $\alpha$ with the knot, 
then the action of the fundamental group on the stalk module is given
by $\alpha(\gamma)=t^{\ell(\alpha)}\gamma$, and this completely determines the
coefficient system which we shall call $\varGamma$. It is not hard to see
that the (simplicial or
singular) chain complex of $\Gamma$-modules determined by this coefficient
system on the knot complement is equivalent to the chain complex of the infinite
cyclic cover with rational coefficients. Thus, if $X$ stands for the knot
complement, the homology modules
$H_*(X; \vg)$ and $H_*(\td X; \Q)$
are isomorphic. (See, for example, \cite{Ha} for a related discussion of the
relationship between homology with local coefficients and homology of covering
spaces).  

The procedure for forming a frame twist-spun knot from a lower dimensional
knot is similar to the procedure for frame spinning except that we add a
``longitudinal twist'' to the gluing. To set up the proper language, we adapt
some notation from Section 6 of Zeeman's paper, \cite{Z65}, in which he
introduces twist spinning. Following Zeeman, if we consider the unit sphere
$S^{m-1}$ in the Euclidean space $\R^m=\R^{m-2}\times \R^2$, then we can
define the latitude for a point $y\in S^{m-1}$ as its projection onto
$\R^{m-2}$ and its longitude as the angular polar coordinate of the
projection of $y$ onto the $\R^2$ term. Hence the latitude is always
well-defined, while the longitude is either undefined or a unique point of
$S^1$ dependent on whether or not $y$ lies in the sphere $S^{m-3}$ that is
the intersection of $S^{m-1}$ with $\R^{m-2}\times 0$. Notice that in the
case where the longitude in undefined, the point on the sphere is uniquely
determined by its latitude (just as on a standard globe). As in Zeeman's
paper, to simplify the notation in abstract cases, we will simply refer to the
latitude-longitude coordinates, $(z, \theta)$, in either case.

Now, just as for frame spinning, we choose a knot $K\subset S^m$ and form the
pair $(D_{+}^{m}, D_{+}^{m-2})=(S^{m}, K)- (D_{-}^{m}, D_{-}^{m-2})$ by
removing a trivial (unknotted) disk pair. We can then identify the trivial boundary
sphere pair $(S^{m-1}, S^{m-3})$ with the unit sphere of the preceding
paragraph and its intersection with $\R^{m-2}\times 0$. Thus, each boundary
point in $(S^{m-1}, S^{m-3})$ can be described by its latitude and longitude
coordinates $(z,\theta)\in D^{m-2}\times S^1$. Then $M^k\times (D_{+}^{m},
D_{+}^{m-2})$ gives a bundle of knots, and the points in $\bd[M^k\times
(D_{+}^{m},D_{+}^{m-2})]$ have coordinates $(x,z,\theta)$, where $x\in M$ and
$(z,\theta)$ are the latitude-longitude coordinates of $\bd (D^m, D^{m-2})$. 

Similarly, given an embedding of $M^k\subset S^{m+k-2}$ with framing $\phi$,
where $S^{m+k-2}$ is the $(m+k-2)$-sphere embedded in $S^{m+k}$ with the
standard normal bundle, we form
\begin{equation*}
(S^{m+k}, S^{m+k-2})-M^k\times
\text{int}(D^{m-2}\times D^2, D^{m-2})
\end{equation*}
as in the frame spinning
construction (Section \ref{S: frame spin}).  Again the boundary can be identified as $M^k \times (S^{m-1},
S^{m-3})$, and the framing $\phi$, together with the trivial framing of
$S^{m+k-2}$ in $S^{m+k}$, allows us to assign to this boundary the same
$(x,z,\theta)$-coordinates.

Given a map $\tau: M^k\to S^1$, we can form the frame twist-spun knot
$\sigma_M^{\phi,\tau}(K)$ as
\[
[(S^{m+k}, S^{m+k-2})-M^k\times \text{int}(D^{m-2}\times D^2, D^{m-2})]
\cup_f [M^k\times (D_{+}^{m}, D_{+}^{m-2})],
\]
where $f$ is the attaching homeomorphism of the boundaries 
\[f:\bd[M^k\times (D_{+}^{m}, D_{+}^{m-2})]\to\bd[(S^{m+k}, S^{m+k-2})-M^k\times \text{int}(D^{m-2}\times D^2, D^{m-2})]\]
  which, identifying each with $M^k \times (S^{m-1},
S^{m-3})$ as above, takes
$(x,z,\theta)\to (x, z, \theta + \tau(x))$, where we define the addition in
the last coordinate as the usual addition on $S^1$. The map $f$ is clearly
well-defined on $M^{k}\times (S^{m-1}-S^{m-3})$ and also on $M^k\times
S^{m-3}$, if we ignore the undefined longitude coordinate. To see that this is
a well-defined continuous map overall, simply observe that on each sphere
$*\times (S^{m-1}, S^{m-3})$, the map is just the rotation by angle $\tau(x)$
of the longitude coordinate induced by the rotation in the second factor of
$\R^{m-2}\times \R^2$. Considered along with the continuity of $\tau$, $f$ is
obviously a homeomorphism.
 
Roughly speaking, we are once again removing a bundle of trivial knots over
$M$ and replacing it with a bundle of non-trivial knots. The new element is
the longitudinal twist determined by $\tau$. The framing $\phi$ employed in
non-twist frame spinning dictates how the trivial bundle of knots over $M$ is
attached ``latitudinally'', while the addition of ``twist'' allows us to
alter the attachment ``longitudinally''. As an example, if $M$ is taken as
the standard circle $S^1\subset S^{m+k-2}$ with $\phi$ the trivial framing, then
$\sigma_{S^1}^{\phi}(K)$ gives us the spun knot of Artin, but if $\tau:
S^1\to S^1$ is a map of degree $k$, then $\sigma_{S^1}^{\phi,\tau}(K)$ is the
$k$-twist spun knot of Zeeman \cite{Z65}. Note also that if $\tau$ is the trivial map, which will always be the case if $M$ is simply-connected, then the frame twist-spin $\sigma_{S^1}^{\phi,\tau}(K)$ is simply the standard frame-spin
$\sigma_{S^1}^{\phi}(K)$.

We next wish to compute the Alexander polynomials of the frame twist-spun
knots. First, we recall the following basic facts of algebra, some of which
we have used before: Since $\Gamma$ is a principle ideal domain, any
$\Gamma$-module can be written as $\Gamma^{\mathfrak B} \oplus (\oplus_i
\Gamma/(p_i))$ for some $\mathfrak B \geq 0$ and $p_i\in \Gamma$, $p_i\neq 0$. It is
also sometimes assumed that the $p_i$ satisfy $p_i|p_{i+1}$ or some
other similar formula simply to provide a normalization, but we will not
impose that condition here. If $b=0$, then $\prod_i p_i$ is what we have been calling the polynomial associated to the module. We will sometimes refer to the $p_i$ as the invariants or torsion invariants of the module.

We also recall that for $p,q\neq0$, $\Gamma/(p)\otimes_{\Gamma}
\Gamma/(q)\cong \Gamma/(p)*_{\Gamma}\Gamma/(q)\cong\Gamma/(d(p,q))$, where
$\otimes_{\Gamma}$ and $*_{\Gamma}$ represent the tensor and torsion products
over the ring $\Gamma$, respectively, and $d(p,q)$ is the greatest common
divisor of $p$ and $q$ in $\Gamma$. In addition, for any $\Gamma$-module $A$,
$\Gamma\otimes_{\Gamma}A\cong A$ and $\Gamma*_{\Gamma}A\cong0$.   
Since we will always assume $\Gamma$ as our ground ring in the following, we
will often simply use $\otimes$ and $*$ to mean the respective products over
$\Gamma$. Observe that the distributivity of $\otimes$ and $*$ over $\oplus$
allow us to calculate the tensor and torsion products of any two
$\Gamma$-modules $A\cong \Gamma^{\mathfrak B_A}\bigoplus (\oplus_i
\Gamma/(p_i))$ and $B\cong \Gamma^{\mathfrak B_B}\bigoplus (\oplus_i
\Gamma/(q_i))$. If we let $A^i$ stand for the direct sum of $i$ copies of the
module $A$ and $T(A)$ stand for the torsion summand of the module $A$
(i.e. $T(A)\cong \oplus_i \Gamma/(p_i)$), then we obtain the following
formulas:
 \begin{align*}
A\otimes B &\cong \Gamma^{\mathfrak B_A+\mathfrak B_B}\oplus T(A)^{\mathfrak B_B}\oplus T(B)^{\mathfrak B_A}
\oplus (\oplus_{i,j}\Gamma/(d(p_i,q_j)))\\
A*B &\cong \oplus_{i,j}\Gamma/(d(p_i,q_j)).
\end{align*}

We will also be using the fact that an
exact sequence of $\Gamma$-torsion modules can be split into the direct sum of exact
sequences of the $p$-primary summands of the modules (see the proof of Proposition
\ref{P:alt. poly.}).

With these formulas in hand, we can compute the Alexander modules of 
frame twist-spun knots. Suppose that $K$ is the knot $S^{m-2}\subset S^m$
which is to be spun and that its Alexander modules are
$H_j(S^{m}-S^{m-2};\vg)\cong \oplus_l \Gamma/(p_{jl})$ (recall that these will
always be $\Gamma$-torsion modules). We wish to compute the homology modules
$H_j(S^n-S^{n-2};\vg)$, where $n=m+k$, $S^{n-2}\subset S^n$ is the frame
twist-spun knot $\sigma_M^{\phi,\tau}(K)$, and $\vg$ is the local coefficient
system as discussed above. Using the above description of the spun knot, let
\begin{align*}
Y&= S^n-(S^{n-2}\cup\text{int} (M^k\times D^m))\\
Z&= M^k\times (D_+^m-D_+^{m-2}),
\end{align*}
so that
\begin{align*}
Y\cap Z&= M^k\times (S^{m-1}-S^{m-3}) \\
Y\cup Z&= S^n-S^{n-2}.
\end{align*}
Then we can employ the Mayer Vietoris sequence
\begin{equation}\label{E: MV tw}
\to H_j(Y\cap Z;\vg|_{Y\cap Z})\overset{i_*}{\to} H_j(Y;\vg|_Y)\oplus
H_j(Z;\vg|_Z)\to
H_j(Y\cup Z;\vg)\to
\end{equation}
to computer the Alexander modules. 

We now examine the terms and maps of this sequence. 

First, we observe that $Y\sim_{h.e.} D^{n-1}\times S^1\sim_{h.e.} S^1$, just as it is
for the corresponding piece of the exterior of the non-twist frame-spun knot in the
previous section. The $S^1$ here can be viewed as a meridian of the knot outside of a neighborhood of the surgery. Therefore, $H_j(Y;\vg|_Y)\cong H_j(\td Y; \Q)\cong H_j(\R;\Q)$, so 
\begin{equation*}
H_j(Y;\vg|_Y)=
\begin{cases}
 \Q\cong \Gamma/(t-1),& j=0\\
0,& j\neq 0.
\end{cases}
\end{equation*}

For the $Z$ component, we need to investigate the coefficient system $\vg|_Z$. Since $Z$
is a product space, its fundamental group is
$\pi_1(M^k\times(D_+^m-D_+^{m-2}))=\pi_1(M^k)\times \pi_1(D_+^m-D_+^{m-2})$. Therefore,
the action of an element $\alpha\times\beta=(\alpha\times 1)\cdot(1\times \beta)=
(1\times \beta)\cdot (\alpha\times 1)$ of the fundamental group on the stalk $\Gamma$
over the basepoint is determined by the product of the actions of $\alpha$ and $\beta$,
which we can take to be loops in $M^k\times *$ and $*\times(D_+^m-D_+^{m-2})$. But this
means that $\vg|_Z$ is equivalent to the product system $\vg|_{M^k\times *}\boxtimes
\vg|_{*\times(D_+^m-D_+^{m-2})}$. Therefore, we can compute $H_j(Z;\vg|_Z)$ via the
K\"{u}nneth theorem (see \cite{Sp}) to be
\begin{align*}
H_j(Z;\vg|_Z)\cong &\bigoplus_{r+s=j}H_r(M^k; \vg|_{M^k})\otimes H_s(D_+^m-D_+^{m-2};\vg|_{D_+^m-D_+^{m-2}})\\
&\bigoplus \bigoplus_{r+s=j-1}H_r(M^k; \vg|_{M^k})*
H_s(D_+^m-D_+^{m-2};\vg|_{D_+^m-D_+^{m-2}}),
\end{align*}
where we have written $\vg|_{M^k}$ to mean $\vg|_{M^k\times *}$ and similarly for the other term. 

We may observe that the terms $ H_s(D_+^m-D_+^{m-2};\vg|_{D_+^m-D_+^{m-2}})$ are none
other than the Alexander modules for the knot $K$. To see this, we need only show that the action of an element  $\alpha\in \pi_1(D_+^m-D_+^{m-2})$ on the stalk $\Gamma$ is given by multiplication by $t^{\ell_K(\alpha)}$, $ \ell_K(\alpha)$ being the linking number of $\alpha$ with the knot $K$. But $\vg|_{D_+^m-D_+^{m-2}}$ is the restriction of the system $\vg$ on $S^n-S^{n-2}$, so the action on $\Gamma$ of a curve representing $\alpha$ is multiplication by $t^{\ell_{\sigma(K)}(\alpha)}$, where the exponent is the linking numer of the loop $\alpha\subset D_+^m-D_+^{m-2}\subset S^n-S^{n-2}$ with the spun knot. So we need only show that the two linking numbers are equivalent. As an element of $H_1(D_+^m-D_+^{m-2})$ (or $H_1(S^n-S^{n-2})$) under the Hurewicz map, $\alpha$ bounds in $D_+^m\subset S^n$.  If $\alpha=\bd c$ and we use $a\cap b$ to denote the intersection number of the chains $a$ and $b$, then 
\[ \ell_K(\alpha)=c\cap D_+^{m-2}=c\cap S^{n-2}= \ell_{\sigma(K)}(\alpha), \]
where the leftmost and rightmost equalities are taken from the definitions of linking and intersection numbers and the central equality is due to  $ D_+^{m-2}=S^{n-2}\cap D_+^m$.
 
Observe that, because the knot modules are all torsion $\Gamma$-modules
then $H_j(Z;\vg|_Z)$ will also be a torsion $\Gamma$-module. 

The homology modules $H_r(M^k; \vg|_{M^k})$ depend, of course, on $M$ so that we cannot
give a general formula for their structure. However, we can obtain a little more
information about the structure of the coefficient system $\vg|_{M^k}$. In fact, we claim
that the action of $\alpha\in\pi_1(M)$ on the stalk $\Gamma$ is given by multiplication
by $t^{\text{deg}(\tau(\alpha))}$, where $ \text{deg}(\tau(\alpha))$ is the degree of the
map $S^1\to S^1$ given by the image of the loop $\alpha$ under $\tau$. To see this, it is
simplest if we choose basepoints so that $M\times *\subset M\times (D_+^m-D_+^{m-2})$
lies in the boundary $M\times (S^{m-1}- S^{m-3})$. This allows us to consider the loop
which represents $\alpha\in \pi_1(M\times *)$ as lying in the component $Y$ via the
attaching homeomorphism $f$. We need to compute the linking number of $\alpha$ with the
spun knot. As remarked, $Y\sim_{h.e.} D^{n-1}\times S^1\sim_{h.e.} S^1$, where $S^1$
gives a meridian of the knot. So if we let $h:Y\to S^1$ be the homotopy equivalence, we
need only compute the degree of $h\circ f(\alpha)$. But by considering the construction,
we can choose $h$ so that its restriction to $ M\times (S^{m-1}-S^{m-3})\subset Y$ is the
projection to the third coordinate in the $(x, z, \theta)$ coordinate system. In other
words, the map is given by projection to the longitude coordinate. So if $M\times
*=M\times (0,0)$ in the coordinate system, then it is clear from the above description of
$f$ that $ h\circ f\circ\alpha(t)=\tau (\alpha(t))$, so the degree of $h\circ f(\alpha)$
is equal to the degree of $\tau(\alpha)$.

The homology of $Y\cap Z\cong M^k\times (S^{m-1}-S^{m-3})$ can also be computed by the
K\"{u}nneth theorem, but here the result is much simpler becuse $S^{m-1}-S^{m-3}$ is an
unknotted sphere pair. Since $m$ must be $\geq 3$, the same linking number argument
applies to show that the homology modules of $S^{m-1}-S^{m-3}$ with coefficient system
$\vg|_{ S^{m-1}-S^{m-3}}$ are the Alexander modules of a trivial knot. In other words,
$H_0(S^{m-1}-S^{m-3};\vg|_{ S^{m-1}-S^{m-3}})\cong\Gamma/(t-1)$, and the homology is
trivial in all other dimensions. Thus,
\begin{multline*}
H_j(Y\cap Z;\vg|_{Y\cap Z})\cong\\ 
\left(H_j(M^k; \vg|_{M^k})\otimes
\Gamma/(t-1)\right) 
\oplus \left(H_{j-1}(M^k; \vg|_{M^k})* \Gamma/(t-1)\right).
\end{multline*}

Notice that $t-1$ is prime in $\Gamma$, so for any $p\in \Gamma$, $d(t-1,p)$ is
$(t-1)$ or $1$. Therefore, $ H_j(Y\cap Z;\vg|_{Y\cap Z})$ is a direct sum of
$\Gamma/(t-1)$'s.

Next, we claim that the map $i_*$ of the Mayer-Vietoris sequence \eqref{E: MV tw} is
injective. We have computed that all the terms of the sequence are torsion
$\Gamma$-modules except for the $H_j(S^n-S^{n-2};\vg)$, but these must also be torsion
modules because the other terms are or simply because we know that these are knot
modules. We know that exact sequences of $\Gamma$-torsion modules can be broken up into
the direct sum of the exact sequences of their $p$-primary components (see the proof of
Proposition \ref{P:alt. poly.}). We also know that $H_j(S^n-S^{n-2};\vg)$ has no
$(t-1)$-primary component for $j>0$ because $t-1$ does not divide the Alexander
polynomials
(which we know up to similarity must evaluate to $\pm 1$ at $1$). Therefore, on the exact
sequence summand corresponding to the $(t-1)$-primary components, the
$H_j(S^n-S^{n-2};\vg)$ terms are $0$, and $i_*$ must be injective. But $ H_j(Y\cap
Z;\vg|_{Y\cap Z})$ consists entirely of its $t-1$ primary component as noted in the
previous paragraph. Therefore, $i_*$ is injective for all $j>0$. It is also injective for
$j=0$ by standard arguments.

Therefore, we obtain short exact sequences  
\begin{equation*}
0\to H_j(Y\cap Z;\vg|_{Y\cap Z})\overset{i_*}{\to} H_j(Y;\vg|_Y)\oplus
H_j(Z;\vg|_Z)\to
H_j(Y\cup Z;\vg)\to 0,
\end{equation*}
and based upon the previous calculations, we can compute $H_j(Y\cup Z;\vg)\cong H_j(S^n-S^{n-2};\vg)$ to be
\begin{align*}
H_j(S^n-S^{n-2};\vg)\cong &\bigoplus_{\overset{r+s=j}{s>0}}H_r(M^k; \vg|_{M^k})\otimes H_s(D_+^m-D_+^{m-2};\vg|_{D_+^m-D_+^{m-2}})\\
&\bigoplus \bigoplus_{\overset{r+s=j-1}{s>0}}H_r(M^k; \vg|_{M^k})* H_s(D_+^m-D_+^{m-2};\vg|_{D_+^m-D_+^{m-2}})
\end{align*}
for $j>0$.

Supposing that the Alexander modules of the knot $K$ are given as
$H_j(D_+^m-D_+^{m-2};\vg)\cong \oplus_l \Gamma/(\lambda_{jl})$ and $H_j(M^k;
\vg|_{M^k})\cong \Gamma^{\mathfrak B_j}\oplus\oplus_l\Gamma/(\zeta_{jl})$, we can then
compute the Alexander polynomial $\lambda_j^{\tau}(K)$, $j>0$, of the frame twist-spun
knot to be
\begin{equation*}
\lambda_j^{\tau}(t)= \prod_{\overset{r+s=j}{s>0}}\left((\prod_l\lambda_{sl}^{\mathfrak
B_r}\cdot\prod_{i,l}d(\zeta_{ri},\lambda_{sl})\right) \cdot
\prod_{\overset{r+s=j-1}{s>0}}\left(\prod_{i,l}d(\zeta_{ri},\lambda_{sl})\right).
\end{equation*}

We next calculate the ``relative'' and ``boundary'' polynomials $\mu^{\tau}_j(t)$ and $\nu^{\tau}_j(t)$ of the spun knot $\sigma^{\phi,\tau}_M(K)$. Let $\bar p (t)=p(t^{-1})$  for any $p\in \Gamma$, suppose $\mathfrak B_i$ continues to denote the rank of the free $\Gamma$ component of $H_i(M;\vg|_M)$, and let $\td \beta_i$ denote the $i$th reduced Betti
number of $M\times \Sigma$. 
Then, for $j>0$, we can calculate $\mu^{\tau}_j(t)$ using the duality of Alexander
polynomials:

{\small
\begin{align*}
\mu_j^{\tau}(t)&\sim (t-1)^{\td \beta_{j-1}}\bar\lambda^{\tau}_{n-1-i} \\
&= (t-1)^{\td \beta_{j-1}}
\prod_{\overset{r+s=n-j-1}{s>0}}\left((\prod_l\bar\lambda_{sl}^{\mathfrak B_r}\cdot
\prod_{i,l}\overline{d(\zeta_{ri},\lambda_{sl})}\right) \cdot 
\prod_{\overset{r+s=n-j-2}{s>0}}\left(\prod_{i,l}\overline{d(\zeta_{ri},\lambda_{sl})}\right)\\
&= (t-1)^{\td \beta_{j-1}}
\prod_{\overset{r+s=n-j-1}{s>0}}\left((\prod_l\bar\lambda_{sl}^{\mathfrak B_r}\cdot
\prod_{i,l}d(\bar\zeta_{ri},\bar\lambda_{sl})\right) \cdot 
\prod_{\overset{r+s=n-j-2}{s>0}}\left(\prod_{i,l}d(\bar\zeta_{ri},\bar\lambda_{sl})\right).
\end{align*}}

\noindent One can also go a step further and calculate $\mu_i^{\tau}$ in terms of the
$\mu_{sl}$,
which we define as follows: Let $X(K)$ denote the exterior of the knot $K$, and let
$L(K)$ represent intersection of the knot exterior $X(K)$ with the closed  neighborhood
of the singularity $\overline{N(\Sigma)}$ (i.e. the ``link exterior''). Then we know that $H_i(X(K),L(K);\vg)$ has the form $T_i\oplus_l \Gamma/(\mu_{il})$, where $T_i$ is the $(t-1)$-primary summand of $H_i(X(K),L(K);\vg)$, $t-1\nmid\mu_{il}$, and $\mu_{il}\neq 0$. Applying Theorem \ref{T:duality*}, we may assume that each $\mu_{il}=\bar\lambda_{m-i-1,l}$. Thus  
\begin{multline*}
\mu_j^{\tau}(t)\sim
(t-1)^{\td \beta_{j-1}}
\prod_{\overset{r+s=n-j-1}{s>0}}\left((\prod_l\mu_{m-s-1,l}^{\mathfrak B_r}\cdot
\prod_{i,l}d(\bar\zeta_{ri},\mu_{m-s-1,l})\right)\\
\cdot 
\prod_{\overset{r+s=n-j-2}{s>0}}\left(\prod_{i,l}d(\bar\zeta_{ri},\mu_{m-s-1,l})\right).
\end{multline*}

Lastly, to calculate $\nu^{\tau}_j(t)$, we can once again employ the K\"{unneth} theorem since $L(\sigma_M^{\phi,\tau}(K))=M\times L(K)$. We have
\begin{align*}
H_j(L(\sigma_M^{\phi,\tau}(K));\vg)\cong &H_j(M\times L(K);\vg)\\
\cong &\bigoplus_{r+s=j}H_r(M^k; \vg|_{M^k})\otimes H_s(L(K);\vg|_{L(K)})\\
&\bigoplus \bigoplus_{r+s=j-1}H_r(M^k; \vg|_{M^k})* H_s(L(K);\vg|_{L(K)}).
\end{align*}
Based on our previous calculations in Section \ref{S: boundary
knot*}, we know that if we let ${\mathfrak b}_i$ stand for the $i$th  Betti number of
$\Sigma$, the singular set of $K$, then the $(t-1)$-primary summand, $T_j$, of
$H_j(L(K);\vg)$ will be isomorphic to $[\Gamma/(t-1)]^{{\mathfrak b}_{j}}$. (For $j>0$,
we
showed that it
was $[\Gamma/(t-1)]^{\td{\mathfrak b}_{j}}$ for reduced Betti number $\td {\mathfrak
b}_j$, but $\td {\mathfrak b}_j=\mf
b_j$ in this range and clearly
$H_0(\widetilde{L(K)};\Q)\cong [\Gamma/(t-1)]^{\mf b_0}\cong \Q^{\mf b_0}$.)
So we can set $ H_{j}(L(K);\vg|_{L(K)})\cong T_{j} \oplus_l
\Gamma/(\nu_{jl})$, where $T_j\cong [\Gamma/(t-1)]^{\mf b_j}$, $t-1\nmid \nu_{jl}$,
and $\nu_{jl}\neq 0$. Then we can
use the above equation to calculate the Alexander polynomial of  
$H_j(L(\sigma_M^{\phi,\tau}(K));\vg)$:\pagebreak
\begin{align*}
\nu^{\tau}_j(t)\sim &
\prod_{r+s=j}\left((t-1)^{\mathfrak B_r\cdot {\mathfrak
b}_s}\prod_l\nu_{sl}^{\mathfrak B_r}\cdot\prod_i d(\zeta_{ri},t-1)^{{\mathfrak
b}_s}\cdot\prod_{i,l}d(\zeta_{ri},\nu_{sl})\right)\\
\qquad\qquad&\cdot  \prod_{r+s=j-1}\left(\prod_i d(\zeta_{ri},t-1)^{{\mathfrak
b}_s}\cdot
\prod_{i,l}d(\zeta_{ri},\nu_{sl})\right)\\
\sim &(t-1)^{ \beta_j}\prod_{r+s=j}\left(\prod_l\nu_{sl}^{\mathfrak B_r} \cdot 
\prod_{i,l}d(\zeta_{ri},\nu_{sl})\right)
\cdot  \prod_{r+s=j-1}\left(\prod_{i,l}d(\zeta_{ri},\nu_{sl})\right),
\end{align*}
where, for the last line, we have used our knowledge of to what power the $t-1$ factor
should occur, based upon some polynomial algebra and our calculations for 
$\lambda_j^{\tau}$ and $\mu_j^{\tau}$  

\begin{remark} As a special case, we can take $M=S^1$ with the standard trivialization
and $\tau:S^1
\to S^1$ to be a map of degree $k\neq 0$. Then $\sigma_{S^1}^{\phi,\tau}(K)$ is the Zeeman $k$-twist spin of $K$. Since the action of a generator of $\alpha\in\pi_1(S^1)$ on the stalk $\Gamma$ is multiplication by $t^k$, we have
\begin{align*}
H_i(S^1; \vg|_{S^1})\cong
\begin{cases}
\Gamma/(t^k-1), & i=0\\
0,&i\neq 0.
\end{cases}
\end{align*}
Therefore, $\mathfrak B_r=0$ for all $r$, $\zeta_{0,1}=t^k-1$, and all other torsion
invariants $\zeta_{ri}$ are trivially equal to $1$. Thus, for $j>0$, we get the
polynomials: 
\begin{align*}
\lambda_j^{\tau}(t)\sim&\prod_{l}d(t^k-1,\lambda_{jl}) \cdot 
\prod_{l}d(t^k-1,\lambda_{j-1,l})\\
\mu_j^{\tau}(t)\sim&
(t-1)^{\td \beta_{j-1}}
\prod_{l}d(t^{-k}-1,\mu_{m-n+j,l}) \cdot 
\prod_{l}d(t^{-k}-1,\mu_{m-n+j+1,l})\\
\nu^{\tau}_j(t)\sim&
(t-1)^{{\mathfrak b}_j+{\mathfrak b}_{j-1}}\prod_{l}d(t^k-1,\nu_{jl})\cdot 
\prod_{l}d(t^k-1,\nu_{j-1,l}).
\end{align*}
\end{remark}

\begin{remark} If $k=0$, we can check that we obtain the polynomials of the non-twist
frame-spun knots
as in the last section. For $k=1$, note that all of the $\lambda^{\tau}_i$, $0<i<n-1$, are trivial (i.e. similar to $1$), while the $\mu^{\tau}_i$ and $\nu^{\tau}_i$ are all powers of $t-1$. 
\end{remark}

As for the subpolynomials $a^{\tau}_i(t)$, $b^{\tau}_i(t)$, and $c^{\tau}_i(t)$, the
existence of $\Gamma$-torsion terms in $H_i(M;\vg|_M)$, the lack of naturality in the
splitting of the K\"{u}nneth theorem, and the lack of exactness of the tensor and torsion
products make it impossible to derive simple formulae in terms of the subpolynomials of
the knot being spun as we did in Section \ref{S: frame spin} for frame spun knots. This
is not a great loss,
however, since we can always calculate the
subpolynomials from $\lambda^{\tau}_i(t)$, $\mu^{\tau}_i(t)$, and $\nu^{\tau}_i(t)$ by
``dividing in'' from the outside of the exact sequence. In other words, recall that we
can calculate $a^{\tau}_i(t)$, $b^{\tau}_i(t)$, and $c^{\tau}_i(t)$ by 
$c^{\tau}_{n-2}(t)= \lambda^{\tau}_{n-2}(t)$ and then
\begin{align*}
a^{\tau}_{n-3}(t)&= \frac{\mu^{\tau}_{n-2}(t)}{ c^{\tau}_{n-2}(t)}\\
b^{\tau}_{n-3}(t)&= \frac{\nu^{\tau}_{n-3}(t)}{ a^{\tau}_{n-3}(t)}\\
c^{\tau}_{n-3}(t)&= \frac{\lambda^{\tau}_{n-3}(t)}{b^{\tau}_{n-3}(t)}\\
&\vdots&\hfill.
\end{align*}
Of course we could also begin from the other side with $c^{\tau}_1(t)$ equal to $\mu^{\tau}_1(t)$ divided by its $t-1$ terms, and so on.

Lastly, we summarize the above calculations as the following realization theorem:

\begin{theorem}\label{T: twist spin real}
Let $M^k$, $n-k>3$, be a manifold which embeds in $S^{n-2}$ with trivial normal bundle
with framing $\phi$. Given a map $\tau: M\to S^1$, let $\mathfrak B_i$ be the rank of the
free part and $\zeta_{il}$ be the torsion invariants of the $\Gamma$-modules
$H_i(M;\vg|_M)$, where the coefficient system is given as above. (Note that these modules
are independent of the knot being spun in the construction.) Then, if $K$ is a knot
$S^{m-2}\subset S^m$ with Alexander invariants $\lambda_{il}$, $\mu_{il}$, and
$\nu_{il}$
and with singular set $\Sigma$ with reduced Betti numbers $\td{\mathfrak b}_i$, then
there exists a frame twist-spun knot $\sigma^{\phi,\tau}_M(K)$ with singular set $M\times
\Sigma$ (whose reduced Betti numbers we denote $\td \beta_i$) and with Alexander
polynomials given for $j>0$ by:
\begin{align*}
\lambda_j^{\tau}(K)\sim& 
\prod_{\overset{r+s=j}{s>0}}\left((\prod_l\lambda_{sl}^{\mathfrak
B_r}\cdot\prod_{i,l}d(\zeta_{ri},\lambda_{sl})\right) \cdot 
\prod_{\overset{r+s=j-1}{s>0}}\left(\prod_{i,l}d(\zeta_{ri},\lambda_{sl})\right)\\
\mu_j^{\tau}(K)\sim&
(t-1)^{\td \beta_{j-1}}
\prod_{\overset{r+s=n-j-1}{s>0}}\left((\prod_l\mu_{m-s-1,l}^{\mathfrak B_r}\cdot
\prod_{i,l}d(\bar\zeta_{ri},\mu_{m-s-1,l})\right)\\
& \qquad\qquad\cdot 
\prod_{\overset{r+s=n-j-2}{s>0}}\left(\prod_{i,l}d(\bar\zeta_{ri},\mu_{m-s-1,l})\right)\\ 
\nu^{\tau}_j(t)\sim &
(t-1)^{\td \beta_j}\prod_{r+s=j}\left(\prod_l\nu_{sl}^{\mathfrak B_r} \cdot
\prod_{i,l}d(\zeta_{ri},\nu_{sl})\right)
\cdot  \prod_{r+s=j-1}\left(\prod_{i,l}d(\zeta_{ri},\nu_{sl})\right).
\end{align*}
In particular, by frame twist-spinning knots with a single point as their singular set,
we obtain knots with $M$ as their singular sets. 
\end{theorem}

\begin{remark}
Although we have focused on realizing given Alexander polynomials in our previous
constructions of knots with point singularities, observe that the methods of proof
actually allow us to create knots with given \emph{invariants}. In fact, we can create
disk knots with specific single invariants and string these together using the disk knot
sum. Then coning on the boundary gives us a knot with the same invariants and a point
singularity. Putting this together with the above theorem, we know exactly what kinds of
polynomials can be realized as those of frame twist-spun knots with singular set $M$,
modulo our ability to compute the homology $H_j(M;\vg|_M)$ and our previous difficulty
with the polynomial $c_2(t)$ of a disk knot $D^3\subset D^5$.
\end{remark}

\begin{remark}
We can, of course,  further enrich the class of polynomials we can realize as polynomials of a knot with singular set $M$ by attaching locally-flat knots to our frame twist-spun knots using ordinary knot sums away from the singularities.  
\end{remark}

%% file: sstratad.tex
10/11/01 \today

\subsubsection{Suspensions}\label{S: suspend}

Another method for obtaining new knots from old ones is by suspension, which, in
some sense, constitutes the extreme case, as the number of singular strata will
always increase. In particular, if we begin with the knotted sphere pair
$S^{n-3}\subset S^{n-1}$ with singular set $\Sigma$, filtered by the nested
subsets $\Sigma_i$ and with ``pure strata''
$U_k=\Sigma_{n-k+1}-\Sigma_{n-k}$, then the suspension, thought of as
$(S^{n-1},S^{n-3})\times I/\{x\times 0\sim *_-, x\times 1 \sim *_+\}$, is a
sphere pair $S^{n-2}\subset S^n$. Its singular set is given by the
suspension of $\Sigma$, and it is filtered by the suspensions points,
$\{*_{\pm}\}$, and the suspensions of the $\Sigma_i$.

We will employ the italicized
$\varSigma$ to denote suspensions. Thus, the suspension of the knot $K$ with 
Alexander polynomials  $\lambda_i\sim b_ic_i$, $\mu_i\sim c_ia_{i-1}$, and
$\nu_i\sim a_ib_i$ will be
denoted by $\vs K$ with polynomials $\lambda_i^{\vs}$, $\mu_i^{\vs}$, and
$\nu_i^{\vs}$. We will compute these polynomials.

\begin{proposition}\label{P: suspend real}
With the notation as above,
\begin{enumerate}
\item $\lambda_i^{\vs}\sim \lambda_i\sim b_ic_i$
\item $\mu_i^{\vs}\sim \mu_{i-1}\sim c_{i-1}a_{i-2}$
\item $\nu_i^{\vs}\sim a_{i-1}b_{i}c_i^2$.
\end{enumerate}
\end{proposition}
\begin{proof}
We first observe that $\lambda_i^{\vs}=\lambda_i$. This follows immediately from
the fact that the suspension points $*_{\pm}$ lie in the knot $\vs K$ so that
$S^n-\vs K\cong (S^{n-1}-K)\times (0,1)$. Therefore, $S^n-\vs K\sim_{h.e.}
S^{n-1}-K$, and $H_i(S^n-\vs K;\vg)\cong H_i(S^{n-1}-K;\vg)$. The claim
follows because  $\lambda_i^{\vs}$ and $\lambda_i$ are the polynomials
associated to these modules, respectively. (Note that the local coefficient
system, $\vg$, on $S^n-\vs K$ is simply that induced by the homotopy equivalence
with $S^{n-1}- K$).

We next turn to the computation of $\nu_i^{\vs}$. We will use $N$ to denote 
open regular neighborhoods and $\bar N$ to denote
closed regular neighborhoods, letting the context in each case determine
the
ambient space. Then $\nu_i^{\vs}$ will be the 
polynomial associated to the homology of the ``link compliment''
$\bd \bar N(\vs(\Sigma))-\vs K\cap \bd \bar N(\vs(\Sigma))$ or,
equivalently, the ``link exterior'' 
$\bd \bar N(\vs(\Sigma))- N(\vs K\cap \bd \bar N(\vs(\Sigma)))$. For the current
argument, it is simplest
if we think of the suspended knot pair as $\left[(S^{n-1},
S^{n-3})\times I\right]\cup
\amalg_{\pm} \bar c_{\pm} (S^{n-1}, S^{n-3})$, where $\amalg_{\pm} \bar c_{\pm}
(S^{n-1},
S^{n-3})$ indicates the disjoint union of the ``northern'' ($+$) and
``southern'' ($-$) closed cones on the
original knot pair. The cones attach to the product with the unit
interval in the
obvious manner. In this case, it is clear that $\bar N(\vs(\Sigma))\cong
\bar N(\Sigma)\times I\cup \amalg_{\pm} \bar c_{\pm}S^{n-1}$ and
$\bd \bar N(\vs(\Sigma))\cong (\bd \bar N(\Sigma)\times I)\cup \amalg_{\pm}
(S^{n-1}- N(\Sigma))$. Finally, since $ N(\Sigma)\subset 
N(K)\subset S^{n-1}$, when we remove the neighborhood around $\vs K$ in $\bd
\bar N(\vs(\Sigma))$, we
see that 
\begin{multline*}
\bd \bar N(\vs(\Sigma))- N(\vs K\cap \bd \bar N(\vs(\Sigma)))\cong\\
\left[(\bd \bar N(\Sigma)-N(K\cap \bd \bar N(\Sigma)))\times I\right]\cup
\amalg_{\pm}
(S^{n-1}_{\pm}-  N(K_{\pm}))
\end{multline*}
glued together at $(\bd \bar N(\Sigma)-N(K\cap \bd \bar N(\Sigma)))\times \{0\}$
and $(\bd \bar N(\Sigma)-N(K\cap \bd \bar N(\Sigma)))\times \{1\}$.

For simplicity, let us readopt some of the notation of Section \ref{S:
disk
knots}. Let us set $X=\bd
\bar
N(\Sigma)-N(K\cap \bd \bar N(\Sigma))$, the link complement of $K$,  and
$C=S^{n-1}- N(K)$, the knot complement of 
$K$. Let us also define
$X_{\vs}=\bd \bar N(\vs(\Sigma))-N(\vs K\cap \bd \bar N(\vs(\Sigma)))$, the
link
complement of $\vs K$, and $C_{\vs}=S^n-N(\vs K)$, the knot complement of
$\vs K$. We will also continue to use $+$ and $-$ in the subscript as
indicators in the cases
where there are multiple copies.
From the preceding paragraph, we can form a Mayer-Vietoris sequence (in which
the coefficient system $\vg$ or its restriction is implied):
\begin{equation*}
\to H_i(X_+)\oplus H_i(X_-) \overset{i_*}{\to} H_i(X\times I)\oplus
H_i(C_+)\oplus
H_i(C_-)\to H_i(X_{\vs}) \to.
\end{equation*}

Now, to study the polynomials, we already know that each  $H_i(X)$ has
associated polynomial $\nu_i\sim a_ib_i$ and each $ H_i(C)$ has associated
polynomial $\lambda_i\sim b_ic_i$. So, according to the results of Section
\ref{S: poly alg}, we can determine the polynomial associated to
$H_i(X_{\vs})$
by determining the polynomial associated to the kernel of $i_*$. 
But, by the definition of the Mayer-Vietoris sequence, the map $i_*$ is induced
by inclusion so that each induced  map $H_i(X_\pm)\to H_i(X\times I)$ is the
identity and each map $H_i(X_\pm)\to H_i(C_\pm)$ is the standard map, say $j_*$, 
induced by
inclusion. Form this, and identifying $H_i(X_+)\cong H_i(X_-)\cong
H_i(X\times I)$ and $H_i(C_+)\cong H_i(C_-)$,  we can conclude that the
kernel of $i_*$ consists of pairs $(\alpha,-\alpha)$, where $\alpha\in H_i(X)$
and furthermore $\alpha\in\text{ker}(j_*)$. This implies that
ker$(i_*)\cong$ker$(j_*)$ and that the polynomial associated to the kernel of
$i_*$ is
$a_i$, as, by definition, this is the polynomial of the kernel of $j_*$. Hence, 
in
the exact sequence, 
the natural factorization of the polynomial associated to
$H_i(X_+)\oplus H_i(X_-)$ is as $a_i$ times $a_ib_i^2$, and the natural
factorization of the polynomial associated to $H_i(X\times I)\oplus
H_i(C_+)\oplus H_i(C_-)$ is as $a_ib_i^2$ times $b_ic_i^2$. Applying this
argument
in all dimensions, we see then that the polynomial  $\nu_i^{\vs}$ associated to
$H_i(X_{\vs})$ must be $a_{i-1}b_ic_i^2$. 

Finally, note that in the long exact sequence of the pair $(C_{\vs},X_{\vs})\cong
(C\times I,X_{\vs})$, 
\begin{equation*}
\begin{CD}
@>>> H_i(X_{\vs}) @>j_*>> H_i(C\times I) @>>> H_i(C\times I, X_{\vs}) @>>>,
\end{CD}
\end{equation*}
the map $j_*:H_i(X_{\vs}) \to H_i(C\times I)$ is an epimorphism. This is due to
the fact that any cycle in $C\times I$ can be homotoped to a cycle in $C\times
[0]\subset X_{\vs}$. Since the polynomial of $H_{i-1}(C\times I)$ is
$\lambda_{i-1}\sim
\lambda_{i-1}^{\vs}\sim b_{i-1}c_{i-1}$ and the polynomial of $ H_{i-1}(
X_{\vs})$
is
$a_{i-2}b_{i-1}c_{i-1}^2$, this implies that the polynomial of
$H_i(C\times I,X_{\vs})$ is $\mu_i^{\vs}\sim a_{i-2}c_{i-1}$.
\end{proof}